\documentclass[a4paper,10pt]{article}
\usepackage[utf8]{inputenc}
\usepackage{biblatex}
\bibliography{bibliography.bib}
\usepackage{hyperref}
\usepackage{float}
\usepackage{amsmath}
\usepackage{amsthm}
\usepackage{amsfonts}
\usepackage[bbgreekl]{mathbbol}
\usepackage{bbm}
\usepackage{geometry}
\usepackage{graphicx}
\usepackage{xcolor}
\usepackage{tikz}
\usepackage{multirow}
\usepackage{siunitx}
\usepackage{listings}
\usepackage{xypic}
\usepackage{subfiles}
\usepackage{authblk}
\usepackage{caption}
\usepackage{subcaption}
\usepackage{comment}

\graphicspath{{./pics/}}

\theoremstyle{definition}

\theoremstyle{remark}
\newtheorem{remark}{Remark}[section]

\newcommand{\bs}{\boldsymbol}
\newcommand{\nun}{\boldsymbol{\nu}}
\newcommand{\dieltens}{\bs{\doubleunderline{\varepsilon}}}
\newcommand{\tens}[1]{\bs{\doubleunderline{#1}}}

\newcommand{\curl}{\textrm{curl}}

\newcommand{\nomegap}{\hat{\omega}_p}
\newcommand{\nomegac}{\hat{\omega}_c}
\newcommand{\nb}{\boldsymbol{b}}

\newcommand{\bx}{\bs{x}}

\newcommand{\gb}{\mathrm{gb}}
\newcommand{\inc}{\mathrm{inc}}
\newcommand{\scat}{\mathrm{sca}}
\newcommand{\ex}{\mathrm{ex}}
\newcommand{\CFL}{\mathrm{CFL}}

\newcommand{\bE}{\bs{E}}
\newcommand{\bB}{\bs{B}}
\newcommand{\bY}{\bs{Y}}

\newcommand{\arr}[1]{{\bs{\mathsf{#1}}}}
\newcommand{\arrE}{\arr{E}}
\newcommand{\arrB}{\arr{B}}
\newcommand{\arrY}{\arr{Y}}

\newcommand{\mat}[1]{\mathbb{#1}}
\def\matM{\mat{M}}
\def\matG{\mat{G}}
\def\matC{\mat{C}}
\def\matD{\mat{D}}
\def\matR{\mat{R}}
\def\matB{\mat{B}}
\def\matA{\mat{A}}

\def\matI{\mat{I}}
\def\matZ{\mat{0}}

\def\RR{\mathbbm{R}}

\newcommand{\psydac}{\textit{Psydac} }

\def\doubleunderline#1{\underline{\underline{#1}}}

\title{Time-splitting methods for the cold-plasma model using Finite Element Exterior Calculus}

\author[a,b]{Elena Moral Sánchez \thanks{Corresponding author. Email address: \href{mailto:Elena.Moral.Sanchez@ipp.mpg.de}{Elena.Moral.Sanchez@ipp.mpg.de}}}
\author[a,b]{Martin Campos Pinto}
\author[b]{Yaman Güçlü}
\author[a,b]{Omar Maj}
\affil[a]{\small{\textit{Department of Mathematics, Technical University of Munich, Boltzmannstraße 3, 85748 Garching, Germany}}}
\affil[b]{\small{\textit{Max-Planck Institute for Plasma Physics, Boltzmannstraße 2, 85748 Garching, Germany}}}

\begin{document}
\newgeometry{left=2.5cm, right=1.5cm, top=2cm, bottom=2cm}

\maketitle

\begin{abstract}
In this work we propose a high-order structure-preserving discretization
of the cold plasma model which describes the propagation of electromagnetic waves in magnetized plasmas.
By utilizing B-Splines Finite Elements Exterior Calculus, we derive a space discretization that preserves the underlying Hamiltonian structure of the model, and we study two stable time-splitting geometrical integrators.
We approximate an incoming wave boundary condition in such a way that the resulting schemes are compatible with a time-harmonic / transient decomposition of the solution, which allows us to establish their long-time stability. This approach readily applies to curvilinear and complex domains. We perform a numerical study of these schemes which compares their cost and accuracy against a standard Crank-Nicolson time integrator, and we run realistic simulations where the long-term behaviour is assessed using frequency-domain solutions.
Our solvers are implemented in the Python library \psydac which makes them memory-efficient, parallel and essentially three-dimensional.

\vspace{0.5cm}

\end{abstract}

\section{Introduction}
The cold-plasma model describes the propagation of a high-frequency electromagnetic wave in a magnetized plasma, by coupling Maxwell's system with a linearized current density equation corresponding to the response of the plasma close to an equilibrium state. Simulation codes constitute an important tool to study the various phenomena involved in the plasma-wave interaction. Turbulent magnetized plasma is an inhomogeneous medium, and highly anisotropic due to the strong background magnetic field. For references, see the classical books \cite{Stix1992WavesIP, swanson2020plasma, Brambilla1998}.

Despite its apparent simplicity, the model is able to accurately simulate physical scenarios relevant to reflectometry diagnostics on fusion plasmas \cite{reflectometry_FDTD_Silva_2019, ITER_position_reflectometers} and plasma wave heating \cite{jacquot_2d_2013}, which consists of using electromagnetic waves to heat up the plasma in magnetic-confinement fusion devices. For instance, Electron Cyclotron Resonance Heating (ECRH) is one of the heating systems of the ITER Tokamak\footnote{\href{https://www.iter.org/machine/supporting-systems/external-heating-systems}{https://www.iter.org/machine/supporting-systems/external-heating-systems}}. For the simulations of plasma wave heating, thermal effects are relevant in the core of the plasma where the wave interacts with the plasma particles, but the cold-plasma model remains valid for most of the domain and can be used to study the interaction of the wave with turbulent plasma fluctuations \cite{waveheating_Heuraux_2015}.

In simplified regimes, reduced models such as the Helmholtz equation or ray-tracing methods have been used, however some applications require the use of so-called full-wave codes, which solve the full vector-valued electromagnetic field \cite{waveheating_Heuraux_2015,DaSilva2024_overviewFDTD}. Certain full-wave codes solve the problem in the frequency-domain \cite{dumont_heating_2012}, i.e. they solve for time-harmonic solutions. Frequency domain codes usually rely on iterative solvers \cite{iterative_methods} which may have convergence issues due to the bad conditionning
of the operators, especially for high frequencies \cite{why_difficult_Helmholtz_Ernst2012} or close to plasma resonances \cite{LiseMarie_phdthesis}. Finding good preconditioners is an ongoing research topic. Domain Decomposition Methods (DDM) \cite{bookDDM_Dolean} have been the main line of research \cite{optimized_schwarz_maxwell_Dolean}. However, as the number of subdomains increases, the convergence deteriorates and, unlike for elliptic problems, the construction of coarse solvers is challenging. Another preconditioning approach has been the shifted-Laplacian method \cite{shiftedlaplacian_preconditioner}, but the trade-off between dissipation and preconditioner accuracy is difficult to assess.

In recent years, the Finite-Difference Time-Domain method (FDTD) has been typically used to run the time-domain problem \cite{tierens_unconditionally_2012, DaSilva2024_overviewFDTD}. A priori this method based on finite differences can be very fast since the computations only involve a small stencil. However, FDTD codes often require a very fine mesh and very small time-steps. The cost is critical for three-dimensional realistic scenarios \cite{2D_vs_3D_FDTD_daSilva_2022}, especially for high frequency. 
Furthermore, long-time stability is often an issue which can lead to the development of unphysical effects \cite{Heuraux2014_differentplasmas_FDTD}. Long-time stable FDTD schemes exist, but at the price of being either first order  \cite{DASILVA201524} or implicit \cite{tierens_unconditionally_2012}, which mitigates their intrinsic efficiency. Another drawback of FDTD codes is that their application to complex geometries is not straightforward \cite{Heuraux2014_differentplasmas_FDTD}. However, the interaction of the wave with the complex structures in the Tokamak vessel and the wave launcher can be essential for reflectometry diagnostics \cite{DaSilva2024_overviewFDTD, Santos2021}.

Our contribution is a full wave code for the time-domain cold-plasma model based on high-order B-spline structure-preserving Finite Element Exterior Calculus (FEEC) \cite{buffa_isogeometric_2011, FEEC_YamanSaidMartin2022} that applies to general domains of curvilinear geometry and is intrinsically stable on very long time ranges. In particular, we present three methods based on geometric integrators \cite{geomIntegration2006}. Two of the methods combine time-splitting techniques \cite{McLachlan_Quispel_2002} with a treatment of the non-Hamiltonian terms, motivated by the preservation of incoming wave boundary conditions.
The third method, used for the purpose of comparison, is a standard energy preserving Crank-Nicolson integrator which does not involve time-splitting.

The presented methods are implemented in the \textit{Psydac} library\footnote{\href{https://github.com/pyccel/psydac}{https://github.com/pyccel/psydac}}, which provides its own sparse stencil matrix format. As a result, our code is memory-efficient, it runs in parallel and is essentially three-dimensional. Besides, block Kronecker mass matrices are used to precondition the implicit steps. Complex geometries can be directly included in the simulation using mappings, 
and physical parameters can be either interpolated from general data or set symbolically thanks to the \textit{SymPDE} library\footnote{\href{https://github.com/pyccel/sympde}{https://github.com/pyccel/sympde}}, so that the code easily adapts to different plasma profiles and background magnetic fields.

The paper is organized as follows. Firstly we present a theoretical study of the properties of the problem and the derivation of the three methods. In particular, we study the long-time stability of two of the proposed schemes. The second part involves a numerical study, where we check the properties of the implemented methods (convergence order, numerical stability and preserved quantities) and compare their performance in terms of cost and accuracy.
The third part includes three realistic two-dimensional numerical simulations. We compute the solutions for long times and track the convergence to the frequency-domain solution.

\section{The cold-plasma model}
In a domain $\Omega \subset \mathbb{R}^3$, the cold-plasma model 
\cite{Brambilla1998} describes the time-evolution of electric and magnetic fields $\bE(\bs{x}, t)$, $\bB(\bs{x}, t)$, together with the current density $\bs{J}(\bs{x}, t)$, and it amounts to the system
\begin{subequations}
\label{dim_model}
\begin{align}
 \partial_t \bE - c \ \curl \bB &= - 4 \pi \bs{J} &\text{ in } \Omega, \label{dim_model:ampere_eq}\\
 \partial_t \bB + c \ \curl \bE &= \bs{0} &\text{ in } \Omega, \label{dim_model:faraday_eq} \\
  \partial_t \bs{J} + \omega_c \bs{J} \times \nb_0 &= (\omega_p^2/4 \pi) \bE - \nu_e \bs{J} &\text{ in } \Omega, \label{dim_model:current_density_eq}
\end{align}
\end{subequations}
where
\begin{equation*}
\omega_p(\bs{x}) := \sqrt{\cfrac{4\pi e^2 n_e(\bs{x})}{m_e}}, \qquad  \omega_c(\bs{x}) := \cfrac{e |\bB_0(\bs{x})|}{m_e c} > 0 \qquad \text{(in CGS units)}
\end{equation*}
are the electron plasma and cyclotron frequencies respectively.
Here, $c$ is the speed of light in vacuum,
$e>0$ the elementary charge, $m_e$ the electron mass,
$n_e(\bs{x})$ denotes the background electron plasma density 
and $\bB_0(\bs{x})$, the background magnetic field with point-wise modulus $|\bB_0| = \sqrt{(B_0)_1^2 + (B_0)_2^2 + (B_0)_3^2}$ and unitary vector field $\nb_0 = \bB_0/|\bB_0|$. The background electron plasma density and magnetic field are assumed to be bounded and constant in time.
Lastly, $\nu_e \ge 0$ denotes the electron-collision frequency, which introduces damping in the system.

Observe that the model is written in Gaussian units, so that $\bE$ and $\bB$ are of the same order, and only the electron species have been considered. In addition the initial condition must satisfy the constraint $\mathrm{div} (\bB) = 0$, which is preserved by 
Faraday's law \eqref{dim_model:faraday_eq}. 

In this article we will complement this model with an impedance
boundary condition corresponding to an incoming time-harmonic 
field of angular frequency $\omega_0 > 0$.

\subsection{Normalization}
In practice, the time-scale of the problem is normalized with the angular frequency $\omega_0$, yielding $\tilde{t} = \omega_0 t$. The spatial scale is normalized with $\omega_0 / c = 2\pi/\lambda_0$, where $\lambda_0$ is the wave-length in vacuum. This yields the normalized variable $\bs{\tilde{x}} = (\omega_0 / c) \bs{x}$.
For simplicity of notation, the tildes are dropped in the rest of the paper, and the quantities correspond to normalized quantities. The normalized cold plasma model reads then
\begin{subequations}
\label{model}
\begin{align}
 \partial_t \bE - \curl \bB &= - \hat{\omega}_p \bY &\text{ in } \Omega, \label{model:ampere_eq}\\
 \partial_t \bB + \curl \bE &= \bs{0} &\text{ in } \Omega, \label{model:faraday_eq} \\
 \partial_t \bY + \hat{\omega}_c \bY \times \nb_0 &= \hat{\omega}_p \bE - \hat{\nu}_e \bY &\text{ in } \Omega,
 \label{model:current_density_eq}
\end{align}
\end{subequations}
where the current density has been rescaled as $\bs{J} = (\omega_p/4 \pi) \bY$ and the normalized parameters are given by $\hat{\omega}_p = \omega_p / \omega_0$,  $\hat{\omega}_c = \omega_c / \omega_0$ and $\hat{\nu}_e = \nu_e / \omega_0$.

The well-posedness of \eqref{model} is usually studied with
the Hille-Yosida theorem, which relies on the fact that the evolution operator is maximal monotone.

\subsection{Boundary conditions}
For simplicity we consider a cubic domain,
with a boundary partitioned into
$\partial \Omega = \Gamma_A \cup \Gamma_P$.
On $\Gamma_P$ we gather faces with periodic boundary conditions 
while $\Gamma_A$ is an artificial boundary where we impose a type of 
impedance conditions 
known as Silver-Müller boundary conditions in time-domain 
\cite{assous2018mathematical}. In normalized units
they take the form 
\begin{equation}
\nun(\bx) \times \left( \bE(\bx, t) - \bB(\bx, t) \times \nun(\bx) \right) 
    = \nun(\bx) \times \bs{s}^\inc(\bx, t)
\qquad \text{ for } \bx \in \Gamma_A,
 \label{SM_BC_general}
\end{equation}
where $\nun$ is the outward unit normal vector on $\Gamma_A$, and
$\bs{s}^\inc$ is a boundary source which may be used to prescribe a given incoming field.
In this paper we consider an incoming wave on a face
$\Gamma_A^\inc := \Gamma_A \cap \{x=x^\inc\}$,
$\bx = (x,y,z)$, and an absorbing boundary condition
on $\Gamma_A \setminus \Gamma_A^\inc$.
Specifically, our incoming wave corresponds to a
time-harmonic circular cross-section Gaussian beam, launched from $x^\inc=0$ 
and propagating in the $x$-axis, with expression
\begin{equation}
    \hat \bE^\gb(\bx) = \bs{e} \frac{w_0}{w(x)} \exp \left[ -\frac{r^2}{w(x)^2} +  i \left(x + \frac{\kappa(x)}{2} r^2 - \psi(x)\right) \right]
    \label{gaussian_beam_3d}
\end{equation}
where $r^2 = (y-y_0)^2 + (z-z_0)^2$ is the radial distance
to $(0,y_0,z_0)$ the position of the focus, $w_0 >0$ is the waist radius,  $\kappa(x) = x / (x^2 + x_R^2)$ is the curvature of the wavefronts, $\psi(x) = \arctan(x/x_R)$ is the Gouy phase and $w^2(x) = w_0^2 \big(1 + (x/x_R)^2\big)$, where $x_R = w_0^2 / 2 > 0$ is the normalized Rayleigh range in vacuum.
The unit vector $\bs{e}$ indicates the polarization, 
for which we will consider different configurations.
Thus, our boundary source reads 
\begin{equation} \label{hs_BC}
  \bs{s}^\inc(\bx,t) = \Re \{ \hat{\bs{s}}^\inc(\bx) e^{- i t}\}
  \qquad \text{ with } \qquad
  \hat{\bs{s}}^\inc(\bx) = 
  \begin{cases}
    \hat{\bE}^\gb(\bx) - \hat{\bB}^\gb(\bx) \times  \nun(\bx) 
      & \text{on } \Gamma_A^\inc
    \\  
    \bs{0} 
      & \text{on } \Gamma_A \setminus \Gamma_A^\inc
  \end{cases}  
\end{equation}
where we have denoted $\hat{\bB}^\gb = -i\, \curl \hat{\bE}^\gb$. 
We remind that the the associated field 
$\bE^\gb(\bx, t) = \Re \{ \hat{\bs{E}}^\gb(\bx) e^{- i t}\}$ is an exact solution of the paraxial wave equation in vacuum, but only an
approximate solution of Maxwell's equations \cite{masterElena}.

\subsection{The frequency-domain problem}
\label{section_freq_domain}
Assuming that the solution of \eqref{model} with the boundary conditions \eqref{SM_BC_general} is time-harmonic of the form
\begin{equation}
\bE^{\mathrm{th}}(t) = \Re \{ \hat{\bE} e^{-it}\}, \quad \bB^{\mathrm{th}}(t) = \Re \{ \hat{\bB} e^{-it}\}, \quad \bY^{\mathrm{th}}(t) = \Re \{ \hat{\bY} e^{-it}\},
\label{cont_timeharmonicsol}
\end{equation}
then $\hat{\bB} = -i \curl \hat{\bE}$, $\hat{\bY} = i/\nomegap (\hat{\bE} - \dieltens \hat{\bE})$ and $\hat{\bE}$ must be the solution of the associated frequency-domain problem
\begin{subequations}\label{freq_problem}
\begin{align}
 \curl \curl \hat{\bE} - \dieltens \hat{\bE} &= \bs{0} \qquad \text{ on } \Omega, \label{freq_problem:eq}\\
 \nun \times (\hat{\bE} + i \curl \hat{\bE} \times \nun )  &=
 \nun \times \hat{\bs{s}}^\inc
 \qquad \text{ on } \Gamma_A, \label{freq_problem:bc}
\end{align}
\end{subequations}
where $\dieltens \hat{\bE}  = S \hat{\bE} - i D \nb_0 \times \hat{\bE} + (P-S) \nb_0 (\nb_0 \cdot \hat{\bE})$ is the dielectric tensor defined with the Stix parameters
\begin{equation*}
S(\bs{x}) = 1 - \cfrac{(1+i \hat{\nu}_e(\bs{x})) \hat{\omega}_p^2(\bs{x})}{(1+i \hat{\nu}_e(\bs{x}))^2 - \hat{\omega}_c^2(\bs{x})}, \qquad
 D(\bs{x}) =  \cfrac{\hat{\omega}_c(\bs{x}) \hat{\omega}_p^2(\bs{x})}{(1+i \hat{\nu}_e(\bs{x}))^2 - \hat{\omega}_c^2(\bs{x})}, \qquad
 P(\bs{x}) = 1 - \cfrac{\hat{\omega}_p^2(\bs{x})}{1+i \hat{\nu}_e(\bs{x})}.
\end{equation*}
In the following we proceed with the dispersion relation analysis, 
for which we assume that $\hat{\nu}_e=0$. In the case of homogeneous plasmas where $\dieltens$ is constant, the analysis resorts to the dispersion relation \cite{Brambilla1998} for plane waves of 
the form $e^{i \bs{k} \cdot \bx} \bs{A}$. In our case, the plasma is inhomogeneous and we resort to the local dispersion relation in the WKB approximation
\cite{Kravtsov1969,Bernstein1975,MCDONALD1988337}.
This relation can be obtained from the symbol 
\cite[Eq.~(3.15)]{MCDONALD1988337}
of the partial differential operator in \eqref{freq_problem:eq}, 
which reads
\begin{align}
  \tens{p}(\bx, \bs{k}) =  - \bs{k} \otimes \bs{k} + \|\bs{k}\|^2 \tens{I} - \dieltens(\bx) \quad \text{ and } \quad \mathrm{Char}(\tens{p}) = \{(\bx, \bs{k}) \in \Omega \times \mathbb{R}^3 | \det(\tens{p}(\bx, \bs{k}) = 0\}
  \label{def_symbol}
\end{align}
is its characteristic variety. Note that since $\hat{\nu}_e = 0$, $\dieltens$ is Hermitian and $\det(\tens{p})$ is real. Following the Stix frame of reference \cite{Stix1992WavesIP}, we align the $z$-axis with $\bs{b}_0$. By symmetry we can assume that $\bs{k}$ is contained in the $x$-$z$ plane, thus $\bs{k} = n(\sin \theta, 0 , \cos \theta)$, where $\theta$ is the angle from $\nb_0$ to the $x$-axis and $n$ is the refractive index (since $\bs{k}$ is normalized to $\omega_0/c$). In this work we are only interested in perpendicular propagation, i.e. $\theta = \pi/2$, for which there are two eigenmodes corresponding to
\begin{itemize}
 \item the O-mode polarization, also known as Transverse Magnetic (TM), where $\hat{\bE} \parallel \nb_0$,
 \item the X-mode polarization, also known as Transverse Electric (TE), where $\hat{\bE} \perp \nb_0$.
\end{itemize}
In this reference frame the dielectric tensor amounts to
\begin{equation}
\dieltens = \begin{bmatrix}
             S & -i D & 0 \\
             iD & S & 0 \\
             0 & 0 & P
            \end{bmatrix} \label{dieltens_Stixframe}
\end{equation}
and the symbol defined in \eqref{def_symbol} becomes
\begin{align}
\tens{p}(\bx, \bs{k}) =  \begin{bmatrix}
                            - S & iD & 0 \\
                            -iD & n^2 - S  & 0 \\
                            0 & 0 & n^2 - P
                           \end{bmatrix} \qquad (\theta=\pi/2).
  \label{def_symbol_Stix}
\end{align}
The symbol consists of two decoupled blocks, corresponding to the X-mode and O-mode symbols, respectively,
\begin{align*}
 \tens{p}_X(\bx, \bs{k}) =  \begin{bmatrix}
                            - S & iD  \\
                            -iD & n^2 - S
                           \end{bmatrix} \qquad \text{and} \qquad
 p_O(\bx, \bs{k}) = n^2 - P.
\end{align*}
Each symbol defines a different branch of the characteristic variety, $\mathrm{Char}(\tens{p}) = \mathrm{Char}(\tens{p}_X) \cup \mathrm{Char}(p_O)$. A simple calculation shows that the points $(\bx, \bs{k})$ within $\mathrm{Char}(\tens{p}_X)$ and $\mathrm{Char}(p_O)$ satisfy
$$
n^2 = S(\bs{x})-D^2(\bs{x})/S(\bs{x}) \quad \text{(X-mode)} \qquad \text{and}  \qquad n^2 = P(\bs{x}) \quad \text{(O-mode)},
$$
where we remind that $n = |\bs{k}|$. These equations are known as the local dispersion relations. Regions with real solutions $n$ are called propagative, whereas regions with complex solutions are called evanescent. Points $\bx \in \Omega$ where $n=0$ are called cutoffs and  points where $n^2 \to \pm \infty$ are called resonances.

From the local dispersion relation it is clear that the O-mode has a cutoff at $P=0$, i.e. when $\nomegap(\bx) = 1$, and does not have any resonances. In the case of the X-mode, there is a cutoff at $S^2=D^2$ and a resonance when $S=0$, namely,
$$
\nomegap^2(\bx) + \nomegac^2(\bx) \pm \nomegac(\bx) \nomegap^2(\bx) = 1 \text{ (X-mode cutoffs)} \qquad \text{and} \qquad
\nomegap^2(\bx) + \nomegac^2(\bx) = 1 \text{ (upper hybrid resonance)}.
$$
Since in our setup the wave is launched from vacuum, we only consider the lower density cutoff (plus sign).

\begin{remark}(Cutoffs and the spectrum of $\dieltens$)
From \eqref{dieltens_Stixframe}, it is clear that the eigenvalues of $\dieltens$ at every $\bx \in \Omega$ are
\begin{equation}
\lambda_1(\bs{x}) = S(\bs{x}) + D(\bs{x}), \quad \lambda_2(\bs{x}) = S(\bs{x}) - D(\bs{x}) \quad \text{and} \quad \lambda_3(\bs{x}) = P(\bs{x}).
 \label{eigenvals_dieltens}
\end{equation}
Therefore, for both polarizations, at least one eigenvalue is zero at the cutoff position.
\label{remark_cutoff_spectrum_dieltens}
\end{remark}

We conclude this section by discussing the well-posedness of the frequency-domain problem \eqref{freq_problem}.
References \cite{sebelin1997}, \cite{wellposedness_PEC_Back2015} show, under the assumption of $\hat{\nu}_e>0$ (this assumption regularizes the problem), the well-posedness of \eqref{freq_problem:eq} with different boundary conditions. In \cite{sebelin1997}, the result is shown for perfect electrical conductor boundary conditions, i.e. $\nun \times \hat{\bE} = \bs{0}$ on $\partial \Omega$. In \cite{wellposedness_PEC_Back2015}, it is proved in the case of Dirichlet boundary conditions $\nun \times \hat{\bE} = \nun \times \hat{\bE}_A$ for a prescribed field $\hat{\bE}_A$ and the case of Neumann boundary conditions $\nun \times \curl \hat{\bE} = \bs{j}_A$ for a prescribed $\bs{j}_A$.

In \cite{LiseMarie_phdthesis} the well-posedness of \eqref{freq_problem} with $\hat{\nu}_e=0$ is shown under certain assumptions.
Among others, all the eigenvalues of $\dieltens$ must be either strictly positive or strictly negative in the domain. As it is shown in remark \ref{remark_cutoff_spectrum_dieltens}, this excludes the cutoffs and the resonance. Another assumption is that $\dieltens$ must be uniformly bounded, which excludes the points where $\nomegac^2 = 1$ which corresponds to the fundamental cyclotron resonance.

In the present work, we assume that there are no resonances in the domain, but allow the existence of cutoffs.

\subsection{Weak form} 
Throughout this paper we denote by 
$\langle \bs{F}, \bs{G} \rangle_\Omega = \int_\Omega \bs{F} \cdot \bs{G} d\bx$ 
the $L^2$-inner product on the domain $\Omega$,
and by
$\langle \bs{F}, \bs{G} \rangle_{\Gamma_A} = \int_{\Gamma_A} \bs{F} \cdot \bs{G} d\bs{S}$
the product on the artificial boundary $\Gamma_A$.
We want to enforce the Silver-Müller boundary conditions in \eqref{SM_BC_general} to the time-domain problem \eqref{model}. For that we impose the Ampère-Maxwell equation \eqref{model:ampere_eq} weakly, while Faraday \eqref{model:faraday_eq} and
the current density \eqref{model:current_density_eq} equations are required to hold strongly.
We consider 
$\bY \in H(\curl, \Omega)$ and $\bB \in H(\mathrm{div}, \Omega)$. 
Given the impedance boundary condition, 
as discussed in \cite{FE_Maxwell_Monk}, a specific space is required for $\bE$,
\begin{align*}
 H_{\mathrm{imp}}(\curl, \Omega) &= \{\bs{U} \in H(\curl, \Omega) | \nun \times\bs{U} \in L^2_t(\partial \Omega) \}
 \\ \text{where } \quad L^2_t(\partial \Omega) &= \{ \bs{U} \in (L^2(\partial \Omega))^3 | \nun \cdot \bs{U} = \bs{0} \text{ on } \partial \Omega \}.
\end{align*}
Consequently the phase space is $V := H_{\mathrm{imp}}(\curl, \Omega) \times H(\mathrm{div}, \Omega) \times H(\curl, \Omega)$,
and the solution 
$(\bE, \bB, \bY)$ in $C^1(\RR_+;V)$ satisfies 
for every $(\bs{F}, \bs{C}, \bs{G}) \in V$ and every $t>0$,
\begin{subequations}
\label{weak_form}
\begin{align}
 \langle \bs{F}, \partial_t \bE \rangle_{\Omega} 
  - \langle \curl \bs{F}, \bB \rangle_{\Omega} 
  + \langle \nun \times \bs{F}, \nun \times \bE \rangle_{\Gamma_A} 
    + \langle \bs{F}, \nomegap \bY \rangle_{\Omega} 
  &= 
  \langle \nun \times \bs{F}, \nun \times \bs{s}^\inc \rangle_{\Gamma_A},
  \label{weak_form:ampere_eq}\\
\langle \bs{C}, \partial_t \bB \rangle_{\Omega} 
  + \langle \bs{C}, \curl \bE \rangle_{\Omega} &=0, 
  \label{weak_form:faraday_eq}\\
\langle \bs{G}, \partial_t \bY \rangle_{\Omega} 
  + \langle \bs{G}\times \bY , \nomegac \nb_0 \rangle_{\Omega} 
  - \langle \bs{G}, \nomegap \bE \rangle_{\Omega} 
  + \langle \bs{G}, \hat{\nu}_e \bY \rangle_{\Omega} &=0. 
   \label{weak_form:current_density_eq}
\end{align}
\end{subequations}
Using proper test functions and integration by parts,
one can check that equation~\eqref{weak_form:ampere_eq} is consistent
with the boundary conditions \eqref{SM_BC_general}.

\subsection{Incoming and scattered fields}

Let $\hat{\bE}^\inc, \hat{\bB}^\inc$ be the complex-valued
solution to the time-harmonic Maxwell equations in vacuum
with Silver-Müller boundary conditions,
\begin{equation}
  \label{model_inc}
  \begin{aligned}
 &\left\{\begin{aligned}
  -i \hat{\bE}^\inc - \curl \hat{\bB}^\inc = 0
  \\
  -i \hat{\bB}^\inc + \curl \hat{\bE}^\inc = 0
 \end{aligned}\right.
 \qquad \text{ on } \Omega 
 \\
 &~~~ \nun \times (\hat{\bE}^\inc - \hat{\bB}^\inc \times \nun )  = 
 \nun \times \hat{\bs{s}}^\inc
 \qquad \text{ on } \Gamma_A 
\end{aligned}
\end{equation}
where the boundary source term $\hat{\bs{s}}^\inc$ 
is defined in \eqref{hs_BC}.
The corresponding (real-valued) time-dependent incoming field
$(\bE, \bB)^\inc(\bx,t) = \Re \{ (\hat{\bE}, \hat{\bB})^\inc(\bx) e^{-it}\}$
satisfies 
\begin{subequations}
  \label{weak_form_inc}
 \begin{align}
  \langle \bs{F} , \partial_t \bE^\inc\rangle_\Omega - \langle \curl \bs{F}, \bB^\inc \rangle_\Omega 
    + \langle \nun \times \bs{F}, \nun \times \bE^\inc \rangle_{\Gamma_A} 
    &= \langle \nun \times \bs{F} , \nun \times \bs{s}^\inc\rangle_{\Gamma_A}, 
  \label{weak_form_inc:ampere_eq} \\
  \langle \partial_t \bB^\inc, \bs{C} \rangle_\Omega 
    + \langle \curl \bE^\inc, \bs{C} \rangle_\Omega &= 0
  \label{weak_form_inc:faraday_eq}
 \end{align}
 \end{subequations}
for all $\bs{F} \in H_{\mathrm{imp}}(\curl, \Omega)$, 
$\bs{C} \in H(\mathrm{div}, \Omega)$.
We next decompose the solution as
\begin{equation}
 \bE = \bE^\inc + \bE^\scat, \qquad \bB = \bB^\inc + \bB^\scat, \qquad \bY = \bY^\scat, \label{scat_field_decomposition}
\end{equation}
where the scattered field is characterized by the system  
\begin{subequations}
 \label{weak_form_scat}
\begin{align}
   \langle \bs{F}, \partial_t \bE^\scat \rangle_{\Omega} - \langle \curl \bs{F} , \bB^\scat\rangle_{\Omega} 
    + \langle \nun \times \bs{F} , \nun \times \bE^\scat \rangle_{\Gamma_A} 
    + \langle \bs{F}, \nomegap \bY^\scat  \rangle_{\Omega}
    &= 0,
   \label{weak_form_scat:ampere_eq}  \\
\langle \bs{C}, \partial_t \bB^\scat \rangle_{\Omega} + \langle \bs{C}, \curl \bE^\scat \rangle_{\Omega} &= 0,
\label{weak_form_scat:faraday_eq} \\
\langle \bs{G}, \partial_t \bY^\scat \rangle_{\Omega} 
  + \langle \bs{G} \times \bY^\scat , \nomegac \nb_0 \rangle_{\Omega}
  - \langle \bs{G}, \nomegap \bE^\scat \rangle_{\Omega} 
  + \langle \bs{G}, \hat{\nu}_e \bY^\scat \rangle_{\Omega} 
  &= \langle \bs{G}, \nomegap \bE^\inc \rangle_{\Omega} 
  \label{weak_form_scat:current_density_eq}
\end{align}
\end{subequations}
for $(\bs{F}, \bs{C}, \bs{G})$ arbitrary in $V$. 
In particular we observe that the scattered field solves the cold-plasma model with purely absorbing boundary conditions, its only source being the interaction of the incoming electric field with the electron plasma density.

\begin{remark}(boundary data for the incoming field)
  In \eqref{weak_form} the boundary condition has been defined using the Gaussian beam profile 
  \eqref{gaussian_beam_3d} for the sake of simplicity. If specific values are available for the 
  incoming wave, e.g., as provided by an experiment or a separate numerical solver, these values
  could be used as well to define the boundary source term $\bs{s}^\inc$.
\end{remark}

\subsection{Hamiltonian structure of the ideal case}

In the case of the \textit{ideal} model, that is, when $\Gamma_A = \emptyset$ and $\hat{\nu}_e=0$, the problem has a non-canonical Hamiltonian structure \cite{morrison1998} in the sense that it can be rewritten in the form
\begin{equation}
 \cfrac{d}{dt} \mathcal{F} (\bE,\bB,\bY) = \{ \mathcal{F}, \mathcal{H} \} (\bE,\bB,\bY) 
 \label{Poisson_system}
\end{equation}
for an arbitrary functional $\mathcal{F} \in C^\infty(V)$.
Here, $\mathcal{H}$ is the Hamiltonian
\begin{equation}
  \mathcal{H}(\bE,\bB,\bY) = \cfrac{1}{2} \left( \| \bE \|^2_{L^2(\Omega)} + \| \bB \|^2_{L^2(\Omega)} + \| \bY \|^2_{L^2(\Omega)} \right)
  \label{Hamiltonian_continuous}
\end{equation}
and the bracket $\{\cdot, \cdot \}$ 
is defined for any smooth functionals $\mathcal{F}, \mathcal{G}$ 
on $V$ such that the functional derivatives are also in $V$ by
\begin{equation}
 \{ \mathcal{F}, \mathcal{G} \} := \{ \mathcal{F}, \mathcal{G} \}_{\mathrm{Maxwell}} + \{ \mathcal{F}, \mathcal{G} \}_{p} + \{ \mathcal{F}, \mathcal{G} \}_{c} \in C^\infty(V),
 \label{Poisson_bracket}
\end{equation}
where
\begin{equation}
 \{ \mathcal{F}, \mathcal{G} \}_{\mathrm{Maxwell}} 
  = \Big \langle \curl \left( \frac{\delta \mathcal{F}}{\delta \bE} \right), \frac{\delta \mathcal{G}}{\delta \bB} \Big \rangle_{\Omega} 
  - \Big \langle \frac{\delta \mathcal{F}}{\delta \bB} , \curl \left( \frac{\delta \mathcal{G}}{\delta \bE} \right)\Big \rangle_{\Omega}
 \label{Maxwell_bracket}
\end{equation}
is the Poisson bracket for Maxwell's equations \cite{abraham2012manifolds}
and
\begin{equation}
 \{ \mathcal{F}, \mathcal{G} \}_{p} 
 = \Big \langle \frac{\delta \mathcal{F}}{\delta \bY} ,  \nomegap \frac{\delta \mathcal{G}}{\delta \bE} \Big \rangle_{\Omega} 
 - \Big \langle \frac{\delta \mathcal{F}}{\delta \bE} , \nomegap \frac{\delta \mathcal{G}}{\delta \bY} \Big \rangle_{\Omega}, \qquad
  \{ \mathcal{F}, \mathcal{G} \}_{c} 
  = \Big \langle  \frac{\delta \mathcal{F}}{\delta \bY} \times \frac{\delta \mathcal{G}}{\delta \bY},  \nomegac \nb_0 \Big \rangle_{\Omega}
  \label{cold_plasma_brackets}
\end{equation}
account for the plasma response and the cyclotron rotation.
The bracket in \eqref{Poisson_bracket} is skew-symmetric and verifies the Leibniz rule and the Jacobi identity, hence it is a Poisson bracket 
\cite{maj_hamiltonian_nodate}.

From the skew-symmetry, it is clear that the Hamiltonian \eqref{Hamiltonian_continuous} is preserved by the ideal system.
An advantage of the form \eqref{Poisson_system} is that it simplifies 
the study of invariants. In particular, one easily sees that the Casimirs
of the bracket, i.e., the functionals $\mathcal{C} \in C^\infty(V)$ such that
\begin{equation}
 \{ \mathcal{C}, \mathcal{G} \} = 0 \qquad \forall 
  \mathcal{G} \in C^\infty(V),
 \label{Casimir_invariants}
\end{equation}
are also invariants of the ideal system \cite{morrison1998}, indeed
taking $\mathcal{G}=\mathcal{H}$ in \eqref{Casimir_invariants} yields
$\frac{d}{dt} \mathcal{C} (\bE,\bB,\bY) = 0$.

\subsection{Non-ideal case and energy balance}
In the general case, \eqref{Poisson_system} becomes  
\begin{equation}
  \cfrac{d}{dt} \mathcal{F} (\bE,\bB,\bY) = \{ \mathcal{F}, \mathcal{H} \} (\bE,\bB,\bY)
  + (\mathcal{F}, \mathcal{H})(\bE,\bB,\bY) 
  + \left \langle \nun \times \frac{\delta \mathcal{F}}{\delta \bE},  \nun \times \bs{s}^\inc \right \rangle_{\Gamma_A}
  \qquad
  \forall \mathcal{F} \in C^\infty(V),
 \label{metriplectic_system}
\end{equation}
 where the Poisson bracket $\{\cdot, \cdot\}$ and Hamiltonian are as above, and where 
 \begin{equation}
  (\mathcal{F}, \mathcal{G}) := 
    -\left \langle \nun \times \frac{\delta \mathcal{F}}{\delta \bE},  \nun \times \frac{\delta \mathcal{G}}{\delta \bE} \right \rangle_{\Gamma_A}
    - 
    \left \langle \hat{\nu}_e \frac{\delta \mathcal{F}}{\delta \bY},  \frac{\delta \mathcal{G}}{\delta \bY} \right \rangle_{\Omega}
\end{equation}
is a symmetric, non-positive bilinear form (sometimes called a metric bracket),
see \cite{maj_hamiltonian_nodate}. 
In the homogeneous case ($\nu \times \bs{s}^\inc = 0$), 
system \eqref{weak_form} has the form of a metriplectic system, however 
we note that this metriplectic form dissipates the Hamiltonian.

In the presence of an incoming source, 
we compute the variation in time of the Hamiltonian by considering $\mathcal{F} = \mathcal{H}$
(i.e., by testing against the solution itself in \eqref{weak_form}), which yields
\begin{equation}
  \partial_t \mathcal{H} = 
  (\mathcal{H}, \mathcal{H})
  + \left \langle \nun \times \frac{\delta \mathcal{H}}{\delta \bE},  \nun \times \bs{s}^\inc \right \rangle_{\Gamma_A}
  =
  -\|\nun \times \bE\|_{L^2(\Gamma_A)}^2 
  - \langle \hat{\nu}_e(\bx) \bY, \bY \rangle_{\Omega} 
  + \langle \nun \times \bE,  \nun \times \bs{s}^\inc \rangle_{\Gamma_A}.
\label{deriv_Hamiltonian_continuous}
\end{equation}
In the case of the ideal model, this shows again that the Hamiltonian is constant. Otherwise, since $\hat{\nu}_e(\bx) \geq 0$, it holds that $\|\nun \times \bE\|_{L^2(\Gamma_A)}^2 + \langle \hat{\nu}_e(\bx) \bY, \bY \rangle_{\Omega} \geq 0$, hence the first two terms can only dissipate the energy. 
In addition, the third term may introduce fluctuations of energy.

Finally, as a direct consequence of Faraday's law \eqref{model:faraday_eq}, we have the invariant
\begin{equation}
 \partial_t (\mathrm{div} \bB) = 0. \label{invariant_divB}
\end{equation}
In addition, Ampère's equation \eqref{weak_form:ampere_eq} implies that the charge density 
$\rho = \mathrm{div} \bE $ satisfies a continuity equation
\begin{equation}
\partial_t \rho + \mathrm{div} (\nomegap \bY) = 0  \qquad \forall t \ge 0.
\end{equation}

Using Gauss' law, the total charge is $Q(t) := \int_{\Omega} \rho d\bx = \int_{\Omega} \mathrm{div} \bE d\bx$ 
By integration by parts we thus have
\begin{equation}
\frac{d}{dt} Q(t) = \int_{\Gamma_A} \partial_t \bE(t) \cdot \nun d \bs{S} 
  = \int_{\Gamma_A} \nomegap \bY  \cdot \nun d \bs{S}
\qquad \forall t \ge 0.
\end{equation}
In the case where $\Gamma_A = \emptyset$ the total charge is trivially preserved, but not in the presence
of an artificial boundary.

\section{Space discretization}
In order to discretize our fields, we use the Finite Element Exterior Calculus (FEEC) framework with
B-splines.

\subsection{FEEC framework}
Following \cite{buffa_isogeometric_2011}, we consider the three-dimensional de Rham diagram:
\begin{equation}
\begin{aligned}
\xymatrix@+1pc{
H^1(\Omega) \ar[d]^{\Pi_0} \ar[r]^-{\mathrm{grad}} & H(\mathrm{curl}, \Omega) \ar[d]^{\Pi_1} \ar[r]^-{\mathrm{curl}} & H(\mathrm{div}, \Omega) \ar[d]^{\Pi_2} \ar[r]^-{\mathrm{div}} & L^{2}(\Omega)  \ar[d]^{\Pi_3} \\
V^0_h \ar[d]^{\bs{\sigma}^0} \ar[r]^{\mathrm{grad}} & V^1_h \ar[d]^{\bs{\sigma}^1} \ar[r]^{\mathrm{curl}} & V^2_h \ar[d]^{\bs{\sigma}^2}  \ar[r]^{\mathrm{div}} & V^3_h
\ar[d]^{\bs{\sigma}^3}
\\
C_h^0 \ar[r]^{\matG} & C_h^1  \ar[r]^{\matC} & C_h^2 \ar[r]^{\matD} & C_h^3
}
\label{deRham_diagram}
\end{aligned}
\end{equation}
where $V^k_h$ refer to the discretized spaces, $\Pi_k$ are the commuting projections, 
$C_h^k$ are the B-spline coefficient spaces, 
$\bs{\sigma}^k$ are the discrete degrees of freedom corresponding to the spline coefficients 
and $\matG, \matC, \matD$ are the gradient, curl and divergence matrices respectively.
In particular, the discretized spaces are finite-dimensional and spanned by the B-spline basis constructed by the tensor-product of one-dimensional B-splines, that is,
\begin{subequations}
\label{Bspline_bases}
\begin{align}
 V^0_h &= \mathrm{span}_{\bs{i}} \{ \Lambda_{\bs{i}}^0 \} \subset L^2(\Omega), &\text{ with } \Lambda_{\bs{i}}^0(\bx) &= N_{i_1}^{p_1} (x) N_{i_2}^{p_2} (y) N_{i_3}^{p_3} (z), \label{Bspline_bases:V0}
 \\[1em]
V^1_h &= \mathrm{span}_{\bs{i}} \{ \bs{\Lambda}_{\bs{i}_x}^1, \bs{\Lambda}_{\bs{i}_y}^1, \bs{\Lambda}_{\bs{i}_z}^1 \}\subset (L^2(\Omega))^3 , &\text{ with }
\bs{\Lambda}_{\bs{i}_x}^1 &=
 \begin{bmatrix}
  \Lambda_{\bs{i}_x}^1\\
  0\\
  0
 \end{bmatrix},
 \bs{\Lambda}_{\bs{i}_y}^1 =
 \begin{bmatrix}
  0\\
  \Lambda_{\bs{i}_y}^1\\
  0
 \end{bmatrix},
 \bs{\Lambda}_{\bs{i}_z}^1 =
 \begin{bmatrix}
  0\\
  0\\
  \Lambda_{\bs{i}_z}^1
 \end{bmatrix}
, \label{Bspline_bases:V1}
\\[1em]
V^2_h &= \mathrm{span}_{\bs{i}} \{ \bs{\Lambda}_{\bs{i}_x}^2, \bs{\Lambda}_{\bs{i}_y}^2, \bs{\Lambda}_{\bs{i}_z}^2 \} \subset (L^2(\Omega))^3 , &\text{ with }
\bs{\Lambda}_{\bs{i}_x}^2 &=
 \begin{bmatrix}
  \Lambda_{\bs{i}_x}^2\\
  0\\
  0
 \end{bmatrix},
 \bs{\Lambda}_{\bs{i}_y}^2 =
 \begin{bmatrix}
  0\\
  \Lambda_{\bs{i}_y}^2\\
  0
 \end{bmatrix},
 \bs{\Lambda}_{\bs{i}_z}^2 =
 \begin{bmatrix}
  0\\
  0\\
  \Lambda_{\bs{i}_z}^2
 \end{bmatrix}
, \label{Bspline_bases:V2}
\\[1em]
 V^3_h &= \mathrm{span}_{\bs{i}} \{ \Lambda_{\bs{i}}^3 \} \subset L^2(\Omega) , &\text{ with } \Lambda_{\bs{i}}^3(\bx) &= D_{i_1}^{p_1} (x) D_{i_2}^{p_2} (y) D_{i_3}^{p_3} (z), \label{Bspline_bases:V3}
\end{align}
\end{subequations}
where $N_{i_k}^{p_k}$ is the $i_k$-th B-spline of degree $p_k$ in the $k$-th direction, $D_{i_k}^{p_k}$ is the Curry-Schoenberg B-spline \cite{de_boor_practical_2001} and
\begin{align*}
 \Lambda_{\bs{i}_x}^1(\bx) &= D_{i_{x,1}}^{p_1} (x) N_{i_{x,2}}^{p_2} (y) N_{i_{x,3}}^{p_3} (z),
 & \Lambda_{\bs{i}_y}^1(\bx) &= N_{i_{y,1}}^{p_1} (x) D_{i_{y,2}}^{p_2} (y) N_{i_{y,3}}^{p_3} (z),
 & \Lambda_{\bs{i}_z}^1(\bx) &= N_{i_{z,1}}^{p_1} (x) N_{i_{z,2}}^{p_2} (y) D_{i_{z,3}}^{p_3} (z),\\
 \Lambda_{\bs{i}_x}^2(\bx) &= N_{i_{x,1}}^{p_1} (x) D_{i_{x,2}}^{p_2} (y) D_{i_{x,3}}^{p_3} (z),
 & \Lambda_{\bs{i}_y}^2(\bx) &= D_{i_{y,1}}^{p_1} (x) N_{i_{y,2}}^{p_2} (y) D_{i_{y,3}}^{p_3} (z),
 & \Lambda_{\bs{i}_z}^2(\bx) &= D_{i_{z,1}}^{p_1} (x) D_{i_{z,2}}^{p_2} (y) N_{i_{z,3}}^{p_3} (z).
\end{align*}
In the above tensor-product spline bases, the differential matrices $\matG, \matC, \matD$ are
essentially connectivity matrices made of 0 and $\pm 1$, thanks to the univariate derivative
formula $(N_{i}^{p})' = D_{i}^{p} - D_{i+1}^{p}$ for any degree $p$ and index $i$.

\subsection{Discretization of the equation}

For our discrete phase space we take $V_h := V_h^1 \times V_h^2 \times V_h^1$. Since $V_h^1 \subset H_{\mathrm{imp}}(\mathrm{curl}, \Omega) \subset H(\mathrm{curl}, \Omega)$, this discretization is conforming in the sense that $V_h \subset V$.
In particular, we can define a semi-discrete solution
\begin{equation}
  \bE_h(t) \in V_h^1, 
  \quad 
  \bB_h(t) \in V_h^2 
  \quad \text{and} \quad 
  \bY_h(t) \in V_h^1
\end{equation}
characterized by the same equations as \eqref{weak_form} but 
with discrete test functions $(\bs{F}, \bs{C}, \bs{G}) \in V_h$.
In terms of (column) vectors of spline coefficients % the numerical solution 
$\arrE = \bs{\sigma}^1 (\bE_h)$, $\arrB = \bs{\sigma}^2 (\bB_h)$, $\arrY = \bs{\sigma}^1 (\bY_h)$, these equations read
\begin{align}
 \begin{bmatrix}
  \partial_t \arrE\\
  \partial_t \arrB\\
  \partial_t \arrY
 \end{bmatrix}
 &=
 \begin{bmatrix}
  - \matM_1^{-1} \mathbb{A}_{1} & \matM_1^{-1}\matC^T \matM_2 & -\matM_1^{-1} \matM_{1, \hat{\omega}_p}\\
  - \matC & \matZ & \matZ \\
  \matM_1^{-1} \matM_{1, \hat{\omega}_p} & \matZ &  - \matM_1^{-1} \left( \matR_{1, \hat{\omega}_c}  + \matM_{1, \hat{\nu}_e} \right)
 \end{bmatrix}
 \begin{bmatrix}
  \arrE\\
  \arrB\\
  \arrY
 \end{bmatrix}
 +
  \begin{bmatrix}
  \matM_1^{-1} \arr{S}^{\inc}(t) \\
  \bs{0}\\
  \bs{0}
 \end{bmatrix},
 \label{semidiscrete_form}
\end{align}
where $\matZ$ denotes the zero matrix, $\matM_k$ is the mass matrix of $V_h^k$
and
\begin{align*}
(\mathbb{A}_{1})_{\bs{i},\bs{j}} 
    &:= \langle \nun \times \bs{\Lambda}^1_{\bs{i}}, \nun \times \bs{\Lambda}^1_{\bs{j}} \rangle_{\Gamma_A}, 
  &  
  (\matM_{1, \hat{\omega}_p})_{\bs{i},\bs{j}} 
    &:= \langle \bs{\Lambda}^1_{\bs{i}}, \hat{\omega}_p \bs{\Lambda}^1_{\bs{j}} \rangle_{\Omega}, 
    \\[0.2em]
  (\matM_{1, \hat{\nu}_e})_{\bs{i},\bs{j}} 
    &:= \langle \bs{\Lambda}^1_{\bs{i}} , \hat{\nu}_e \bs{\Lambda}^1_{\bs{j}}\rangle_{\Omega},
  &
  (\matR_{1, \hat{\omega}_c})_{\bs{i},\bs{j}} 
    &:= \langle \bs{\Lambda}^1_{\bs{i}} \times \bs{\Lambda}^1_{\bs{j}}, \nomegac \nb_0 \rangle_{\Omega}.
\end{align*}
The source term is 
$$
(\arr{S}^\inc)_{\bs{i}} = \langle \nun \times  \bs{\Lambda}^1_{\bs{i}} , \nun \times \bs{s}^\inc \rangle_{\Gamma_A}
= (\cos (t) \arr{S}^{\inc}_R + \sin (t) \arr{S}^\inc_I)_{\bs{i}}
$$
with 
time-independent boundary source arrays given by
\begin{equation} \label{SincRI}
  \left\{
  \begin{aligned}
    (\arr{S}^\inc_R)_{\bs{i}} &:= 
    \langle \nun \times \bs{\Lambda}^1_{\bs{i}}, \nun \times \Re\{\hat{\bs{s}}^{\inc}\} \rangle_{\Gamma_A}
    \\ 
    (\arr{S}^\inc_I)_{\bs{i}} &:= 
    \langle \nun \times \bs{\Lambda}^1_{\bs{i}}, \nun \times \Im\{\hat{\bs{s}}^{\inc}\} \rangle_{\Gamma_A}.
  \end{aligned} 
  \right.
\end{equation}
where $\hat{\bs{s}}^{\inc}$ is given in \eqref{hs_BC}. 
Lastly, we remark that $\matR_{1, \hat{\omega}_c}$ is skew-symmetric and $\mathbb{A}_{1}$ is symmetric positive semi-definite.

\subsection{Discrete Hamiltonian structure}

Gathering the coefficients in a large column array
of the form
$\arr{U} = [\arrE, \arrB, \arrY]^T$ in 
$C^1(\RR_+;C_h)$ with $C_h = C_h^1 \times C_h^2 \times C_h^1$,
system~\eqref{semidiscrete_form} rewrites as
 \begin{equation}
 \partial_t \arr{U} = \mathbb{P} \nabla_{\arr{U}} \mathsf{H} + \mat{N} \nabla_{\arr{U}} \mathsf{H} + \arr{f}(t), \label{coeffs_semidiscrete_form_separated_Hamiltonian_part}
 \end{equation}
where $\mathbb{P}$ and $\mat{N}$ are the Poisson and ``metric'' matrices corresponding to 
the continuous metriplectic system \eqref{metriplectic_system}, 
namely
\begin{align}
 \mathbb{P} &:=
  \begin{bmatrix}
  \matZ & \matM_1^{-1} \matC^T & - \matM_1^{-1} \matM_{1, \hat{\omega}_p} \matM_1^{-1}\\
  - \matC \matM_1^{-1} & \matZ & \matZ \\
  \matM_1^{-1} \matM_{1, \hat{\omega}_p} \matM_1^{-1} & \matZ & - \matM_1^{-1} \matR_{1, \hat{\omega}_c}  \matM_1^{-1}
 \end{bmatrix}, &
 \mat{N} &:=
  \begin{bmatrix}
  - \matM_1^{-1} \matA_{1} \matM_1^{-1}& \matZ & \matZ \\
    \matZ & \matZ & \matZ \\
    \matZ & \matZ & - \matM_1^{-1} \matM_{1, \hat{\nu}_e} \matM_1^{-1}
 \end{bmatrix}.
\label{PN}
\end{align}
Here,
the source is $ \arr{f}(t) := [\matM_1^{-1} \arr{S}^{\inc}(t), \bs{0}, \bs{0}]^T $ 
and the discrete Hamiltonian is
\begin{align}
 \mathsf{H}(\arr{U}) &= \mathcal{H}(\bE_h, \bB_h, \bY_h) = \cfrac{1}{2} ( \arr{E}^T \matM_1 \arr{E} + \arr{B}^T \matM_2 \arr{B} + \arr{Y}^T \matM_1 \arr{Y}),
  \label{Hamiltonian_semidiscrete}
\end{align}
so that $\nabla_{\arr{U}} \mathsf{H} = [\matM_1 \arrE, \matM_2 \arrB, \matM_1 \arrY]^T$.
Another way of writing \eqref{coeffs_semidiscrete_form_separated_Hamiltonian_part} is to 
consider a general functional $\mathsf{F} \in C^\infty(C_h)$ and apply a chain rule:
\begin{equation}
  \cfrac{d}{dt} \mathsf{F} (\arr{U}) = \{ \mathsf{F}, \mathsf{H} \} (\arr{U})
  + (\mathsf{F}, \mathsf{H})(\arr{U}) 
  + (\nabla_{\arr{U}} \mathsf{F})^T \arr{f}(t).
 \label{metriplectic_system_h}
\end{equation}
Since $\mat{P}$ is constant and skew-symmetric it represents a Poisson bracket, moreover 
$\mat{N}$ is symmetric and negative semi-definite matrix. 
Therefore the homogeneous ($\arr{f} = \arr{0}$) equation admits 
a metriplectic structure with discrete Poisson and metric brackets
\begin{align}
 \{\mathsf{F}, \mathsf{G} \} := (\nabla_{\arr{U}} \mathsf{F})^T \mathbb{P} \nabla_{\arr{U}} \mathsf{G},
 \qquad
 ( \mathsf{F}, \mathsf{G} ) := (\nabla_{\arr{U}} \mathsf{F})^T \mathbb{N} \nabla_{\arr{U}} \mathsf{G}
  \qquad \text{ for } 
    \mathsf{F}, \mathsf{G} \in C^\infty(C_h).
\end{align} 
In the ideal case ($\mat{N} = \matZ$ and $\arr{f} = \arr{0}$), system \eqref{metriplectic_system_h}
preserves $\mathsf{H}$ (by antisymmetry of $\mat{P}$) and the Casimirs \cite{geomIntegration2006}, 
i.e. every functional $\mathsf{C} \in C^\infty(C_h)$ such that 
$(\nabla_{\arr{U}} \mathsf{C})^T \mathbb{P} = \bs{0}$.
In the general case, the semi-discrete energy balance involves both non-Hamiltonian terms of 
\eqref{coeffs_semidiscrete_form_separated_Hamiltonian_part}, namely
\begin{equation}
 \partial_t \mathsf{H} 
 = (\mathsf{H}, \mathsf{H}) 
 + (\nabla_{\arr{U}} \mathsf{H})^T \arr{f}(t)
 = -(\arrE^T \matA_1 \arrE + \arrY^T \matM_{1,\hat{\nu}_e} \arrY )
    + \arrE^T \arr{S}^{\mathrm{inc}}(t),
 \label{deriv_Hamiltonian_discrete}
\end{equation}
where $\arrE^T \matA_1 \arrE + \arrY^T \matM_{1,\hat{\nu}_e} \arrY \geq 0$ only contributes to energy dissipation, while the last term may introduce energy oscillations.

\section{Time discretization}
Geometric integrators are methods that preserve certain geometric quantities of the solutions, such as the Hamiltonian or the Casimirs of the Poisson bracket \cite{geomIntegration2006}. Thus, they are of special interest for long times simulations. Splitting methods are a type of geometric integrators consisting of solving certain subsystems sequentially, instead of the whole system at once. Usually this approach leads to smaller subsystems which are easier to solve.
For Hamiltonian systems, geometric splitting methods have been thoroughly studied in \cite{McLachlan_Quispel_2002}.
Here, we design two methods based on such time-splitting schemes where the non-Hamiltonian part 
is handled in a stable way.

The resulting splitting schemes will then be compared with a Crank-Nicolson scheme, which is a standard second-order implicit, stable (energy preserving) time integrator based on the trapezoidal rule.

\subsection{Time splittings for Hamiltonian systems} \label{sec:splittings}
In the ideal case, that is when $\mat{N} = \matZ$ and $\arr{f} = \arr{0}$, the ideal problem
\begin{align*}
  \partial_t \arr{U} = \mat{P} \nabla_{\arr{U}} \mathsf{H}
\end{align*}
can be split following two standard approaches \cite{McLachlan_Quispel_2002}:
\begin{itemize}
 \item Poisson splitting which consists in decomposing the Poisson matrix into skew-symmetric submatrices $\mat{P} = \mat{P}_1 + \mat{P}_2$, leading to
 \begin{equation*}
  \begin{aligned}
   \partial_t \arr{U} &= \mat{P}_1 \nabla_{\arr{U}} \mathsf{H}
   \\
   \partial_t \arr{U} &= \mat{P}_2 \nabla_{\arr{U}} \mathsf{H}.
  \end{aligned}
  \end{equation*}
 If these two subsystems can be solved exactly
 the resulting integrator preserves the Hamiltonian.
 \item Hamiltonian splitting which consists in decomposing the discrete Hamiltonian into $\mathsf{H} = \mathsf{H}_{1} + \mathsf{H}_{2}$, leading to
 \begin{equation*}
  \begin{aligned}
   \partial_t \arr{U} &= \mat{P} \nabla_{\arr{U}} \mathsf{H}_1
   \\
   \partial_t \arr{U} &= \mat{P} \nabla_{\arr{U}} \mathsf{H}_2.
   \end{aligned}
  \end{equation*}
If these two subsystems can be solved exactly
the resulting integrator preserves the Poisson bracket, and in particular its Casimir invariants.
\end{itemize}

In practice these subsystems are solved alternatively, leading to numerical flows which are composed to yield potentially high order time integration schemes. In this article we will consider Strang splitting, which is second-order.

\subsection{Splitting of the non-Hamiltonian terms}
\label{sec:split_non_ham}

For the Hamiltonian part, we can use the splitting methods described above.
We then need a criterion to decide in which subsystem to place the non-Hamiltonian terms
from \eqref{coeffs_semidiscrete_form_separated_Hamiltonian_part}--\eqref{PN},
namely the two boundary terms (associated with $\matA_1$ and $\arr{f}$) and the collision term
(associated with $\matM_{1, \hat{\nu}_e}$).

For the latter we do not see a particular constraint, but for the boundary terms there seems to be only 
one correct configuration, given some splitting of the Hamiltonian terms.
Consider indeed one equation in a subsystem containing the (weak) ``$\curl \bB$'' term and some (possibly vanishing) 
contributions from other terms, that is an equation of the form
$$
\langle \bs{F}, \partial_t \bE \rangle_{\Omega} 
- \langle \curl \bs{F}, \bB \rangle_{\Omega} 
+ \alpha \langle \nun \times \bs{F}, \nun \times \bE \rangle_{\Gamma_A} 
  + \gamma \langle \bs{F}, \nomegap \bY \rangle_{\Omega} 
= 
\beta \langle \nun \times \bs{F}, \nun \times \bs{s}^\inc \rangle_{\Gamma_A},
$$
for some parameters $\alpha, \beta, \gamma \in \RR$.
Then the corresponding boundary condition reads
$$
\nun \times (-\bB \times \nun + \alpha \bE) = \beta \nun \times \bs{s}^\inc
$$
By comparison with \eqref{SM_BC_general}
it appears that the choice $\alpha=\beta=1$ is the only option that ensures that the field satisfies the physically correct boundary conditions.
We confirmed numerically that splitting these boundary terms with $\alpha \neq \beta$ resulted in a very large error originating from the boundary. Indeed these terms are very large and they compensate only when they remain together.

This choice is also physical because it guarantees that in the case of $\hat{\nu}_e =0$ the energy-balance \eqref{deriv_Hamiltonian_discrete} is preserved, i.e. there is no mismatch between energy input and dissipation.

Another justification is based on guaranteeing the physical behaviour of the scattered field.
As we can see in \eqref{weak_form_scat}, the source of the scattered field $\bs{Y}^\scat$ is basically the interaction of the incoming field with the plasma. Thus, initially the scattered field is zero near the incoming boundary $\Gamma_A^\inc$.
Splitting the terms mentioned above, i.e. $\alpha \neq \beta$, introduces a boundary source in \eqref{weak_form_scat:ampere_eq}. As a result, a very large unphysical error appears initially at the incoming boundary, where the scattered field is supposed to be zero.

\subsection{Poisson splitting time scheme}

Following subsection~\ref{sec:splittings}, in this approach we split the Poisson matrix into two 
skew-symmetric matrices. Here a natural choice is to separate the curl-curl operators from the 
magnetic and plasma terms:
$$
\mat{P}_1 := 
\begin{bmatrix}
\matZ & \matM_1^{-1} \matC^T & \matZ \\
- \matC \matM_1^{-1} & \matZ & \matZ \\
\matZ & \matZ & \matZ
\end{bmatrix},
\qquad
\mat{P}_2 := 
\begin{bmatrix}
  \matZ & \matZ & - \matM_1^{-1} \matM_{1, \hat{\omega}_p} \matM_1^{-1}\\
  \matZ & \matZ & \matZ \\
  \matM_1^{-1} \matM_{1, \hat{\omega}_p} \matM_1^{-1} & \matZ & - \matM_1^{-1} \matR_{1, \hat{\omega}_c}  \matM_1^{-1}
 \end{bmatrix}.
$$
According to our discussion above the boundary contributions are kept together with the first (Maxwell) step which contains the $\curl \bB$ term.
Moreover the collision term associated to $\matM_{1, \hat{\nu}_e}$ is kept together with $\mat{P}_2$,
to avoid having to solve for $\bY$ in the first subsystem.

The splitting consists of solving the following subproblems:
\begin{itemize}
 \item a Maxwell subsystem, associated to $\mat{P}_1$,
\begin{subequations}
\label{Poisson_splitting_operatorA}
\begin{align}
 \matM_1 \partial_t \arr{E} & = \matC^T \matM_2 \arr{B} + \cos(t) \arr{S}^\inc_R  + \sin(t) \arr{S}^\inc_I -  \matA_{1} \arr{E}, \\
 \partial_t \arr{B} &= -\matC \arr{E},\\
 \partial_t \arr{Y} &= \arr{0},
\end{align}
\end{subequations}
and
\item a plasma subsystem, associated to $\mat{P}_2$,
\begin{subequations}
\label{Poisson_splitting_operatorBC}
\begin{align}
  \matM_1 \partial_t \arr{E} & = - \matM_{1, \hat{\omega}_p} \arr{Y}, \\
 \partial_t \arr{B} &= \arr{0},\\
 \matM_1 \partial_t \arr{Y} &= \matM_{1, \hat{\omega}_p} \arr{E} - (\matR_{1, \hat{\omega}_c} + \matM_{1, \hat{\nu}_e}) \arr{Y}.
\end{align}
\end{subequations}
\end{itemize}

\subsubsection{Time integration}
Each problem is integrated numerically in the interval $(t^n, t^{n+1} = t^n + \Delta t)$ 
using a second-order trapezoidal rule, leading to the following two discrete flows.
\begin{itemize}
 \item[$\bullet$] A first flow $(\arr{E}^{n+1}, \arr{B}^{n+1}, \arr{Y}^{n+1}) = \Phi_{\mathrm{Maxwell}}^{\mathrm{P}}(\arr{E}^n, \arr{B}^n, \arr{Y}^n; \Delta t, t^n)$
 obtained by integrating the Maxwell subsystem \eqref{Poisson_splitting_operatorA}:
 \begin{subequations}
 \label{integration_Poisson_splitting_operatorA}
  \begin{align}
   \left[ \matM_1 + \cfrac{\Delta t^2}{4} \matC^T \matM_2 \matC + \cfrac{\Delta t}{2} \matA_{1} \right] \arr{E}^{n+1/2} &= \matM_1 \arr{E}^n + \cfrac{\Delta t}{2} \matC^T \matM_2 \arr{B}^n
   +  \cfrac{s^n}{2} \arr{S}_R^\inc - \cfrac{c^n}{2} \arr{S}_I^\inc , \label{integration_Poisson_splitting_operatorA:avg_E}\\
   \arr{E}^{n+1} &= 2 \arr{E}^{n+1/2} - \arr{E}^{n},\\
   \arr{B}^{n+1} &= \arr{B}^n - \Delta t \matC \arr{E}^{n+1/2},\\
   \arr{Y}^{n+1} &= \arr{Y}^n,
  \end{align}
 \end{subequations}
  where $s^n = \sin(t^{n+1}) - \sin(t^n)$ and $c^n = \cos(t^{n+1}) - \cos(t^n)$. The system  \eqref{integration_Poisson_splitting_operatorA:avg_E} is solved using PCG (Preconditioned Conjugate Gradient) \cite{iterative_methods}. For the preconditioner we use a Kronecker mass solver.

 \item[$\bullet$] A second flow $(\arr{E}^{n+1}, \arr{B}^{n+1}, \arr{Y}^{n+1}) = \Phi_{\mathrm{plasma}}^{\mathrm{P}}(\arr{E}^n, \arr{B}^n, \arr{Y}^n; \Delta t)$ obtained by integrating the plasma subsystem~\eqref{Poisson_splitting_operatorBC}:
 \begin{subequations}
  \label{integration_Poisson_splitting_operatorBC}
  \begin{align}
   \begin{bmatrix}
    \matM_1 & \cfrac{\Delta t}{2} \matM_{1, \hat{\omega}_p}\\
    -\cfrac{\Delta t}{2} \matM_{1, \hat{\omega}_p} & \matM_1 + \cfrac{\Delta t}{2} \left( \matR_{1, \hat{\omega}_c} + \matM_{1, \hat{\nu}_e} \right)
   \end{bmatrix}
   \begin{bmatrix}  \arr{E}^{n+1} \\[1em] \arr{Y}^{n+1} \end{bmatrix}
   &=
   \begin{bmatrix}  \matM_1 \arr{E}^{n} - \cfrac{\Delta t}{2} \matM_{1, \hat{\omega}_p} \arr{Y}^{n}
   \\
   \Big( \matM_1 - \cfrac{\Delta t}{2} \left( \matR_{1, \hat{\omega}_c} + \matM_{1, \hat{\nu}_e} \right) \Big) \arr{Y}^{n} + \cfrac{\Delta t}{2} \matM_{1, \hat{\omega}_p} \arr{E}^{n} \end{bmatrix},\\
   \arr{B}^{n+1} &= \arr{B}^n.
  \end{align}
  \end{subequations}
Since the matrix is not symmetric, we use PBiCGStab (Preconditioned Biconjugate Gradient stabilized method) \cite{bicgstab}  to solve the system. The method is preconditioned with a block diagonal matrix, where each block is a Kronecker mass solver.
\end{itemize}

\subsubsection{Poisson splitting scheme}
Lastly, we compose the above flows using a Strang splitting scheme. Given the initial condition $\arrE^0, \arrB^0, \arrY^0$ at $t^0=0$, for each time-step $n=0,...,N-1$, compute
\begin{subequations}
\label{scheme_Poisson}
\begin{align}
  (\arr{E}^{n+1/3}, \arr{B}^{n+1/2}, \arr{Y}^n) &= \Phi_{\mathrm{Maxwell}}^{\mathrm{P}}(\arr{E}^n, \arr{B}^n, \arr{Y}^n; \Delta t/2, t^n),  \\
  (\arr{E}^{n+2/3}, \arr{B}^{n+1/2}, \arr{Y}^{n+1}) &= \Phi_{\mathrm{plasma}}^{\mathrm{P}}(\arr{E}^{n+1/3}, \arr{B}^{n+1/2}, \arr{Y}^n; \Delta t), \\
  (\arr{E}^{n+1}, \arr{B}^{n+1}, \arr{Y}^{n+1}) &= \Phi_{\mathrm{Maxwell}}^{\mathrm{P}}(\arr{E}^{n+2/3}, \arr{B}^{n+1/2}, \arr{Y}^{n+1}; \Delta t/2, t^n + \Delta t/2), \\
  t^{n+1} &= t^n + \Delta t,
\end{align}
\end{subequations}
where we note that $\arr{E}^{n+1/3}, \arr{B}^{n+1/2}, \arr{E}^{n+2/3}$ are intermediate solutions and not
proper approximate solutions for intermediate times.

\subsection{Hamiltonian splitting time scheme}

We now split the discrete Hamiltonian \eqref{Hamiltonian_semidiscrete} as described in subsection~\ref{sec:splittings},
Here we choose to separate $\mathsf{H}$ into two parts:
\begin{align}
\mathsf{H}_{\arr{E}} &:=  \cfrac{1}{2} \arr{E}^T \matM_1 \arr{E}, 
& 
\mathsf{H}_{\arr{B},\arr{Y}} &:=  \cfrac{1}{2} \left( \arr{B}^T \matM_2 \arr{B} + \arr{Y}^T \matM_1 \arr{Y} \right).
 \label{split_hamiltonian}
\end{align}
The alternative splitting given by $\mathsf{H}_{\arr{B}}=(1/2)\arr{B}^T \matM_2 \arr{B}$ and $\mathsf{H}_{\arr{E},\arr{Y}}=(1/2)( \arr{E}^T \matM_1 \arr{E} + \arr{Y}^T \matM_1 \arr{Y})$ yields a large subsystem involving the three fields, which is expensive to solve. In contrast, the choice in \eqref{split_hamiltonian} yields small subsystems which are less costly to solve.
Again, following subsection~\ref{sec:split_non_ham} we place the boundary terms with the second subsystem
($\mathsf{H}_{\arr{B},\arr{Y}}$) since it will contain the $\curl \bB$ term. The splitting consists of solving the following subproblems:

\begin{itemize}
 \item the electric subsystem, corresponding to $\mathsf{H}_{\arr{E}}$, which reads
\begin{subequations}
\label{Hamiltonian_splitting_operatorA}
\begin{align}
 \partial_t \arr{E} & = \arr{0}, \\
 \partial_t \arr{B} &= -\matC \arr{E},\\
 \matM_1 \partial_t \arr{Y} &= \matM_{1,\hat{\omega}_p} \arr{E},
 \end{align}
\end{subequations}
and
  \item the magnetic-plasma subsystem, corresponding to $\mathsf{H}_{\arr{B},\arr{Y}}$, which reads
  \begin{subequations}
  \label{Hamiltonian_splitting_operatorBC}
   \begin{align}
 \matM_1 \partial_t \arr{E} & = \matC^T  \matM_2 \arr{B} + \cos (t) \arr{S}^\inc_R + \sin (t) \arr{S}^\inc_I -  \matA_{1} \arr{E}  - \matM_{1,\hat{\omega}_p} \arr{Y},  \\
 \partial_t \arr{B} &= \arr{0},\\
 \matM_1 \partial_t \arr{Y} &= - \matR_{1, \hat{\omega}_c} \arr{Y} - \matM_{1, \hat{\nu}_e} \arr{Y}.
   \end{align}
  \end{subequations}
\end{itemize}
Note that here 
the collision term associated to $\matM_{1, \hat{\nu}_e}$ is kept for simplicity with the second subsystem but
could alternatively be placed differently, such as in \eqref{Hamiltonian_splitting_operatorA} or in both subsystems with some weighting.

\subsubsection{Time integration}
Each problem is integrated numerically in the interval $(t^n, t^{n+1} = t^n + \Delta t)$ yielding the following discrete flows.
\begin{itemize}
 \item[$\bullet$] A first flow $(\arr{E}^{n+1}, \arr{B}^{n+1}, \arr{Y}^{n+1}) = \Phi_{\arr{E}}^{\mathrm{H}}(\arr{E}^n, \arr{B}^n, \arr{Y}^n; \Delta t)$ obtained by integrating the electric subsystem \eqref{Hamiltonian_splitting_operatorA}:
 \begin{subequations}
\label{integration_Hamiltonian_splitting_operatorA}
\begin{align}
 \arr{E}^{n+1} &= \arr{E}^n, \\
 \arr{B}^{n+1} &= \arr{B}^n - \Delta t \matC \arr{E}^n, \\
 \matM_1 \arr{Y}^{n+1} &= \matM_1 \arr{Y}^n + \Delta t \matM_{1,\hat{\omega}_p}\arr{E}^n, \label{integration_Hamiltonian_splitting_operatorA:Y}
\end{align}
\end{subequations}
which corresponds to an exact integration and therefore $\Phi_{\arr{E}}^{\mathrm{H}}$ preserves the Poisson structure. Equation \eqref{integration_Hamiltonian_splitting_operatorA:Y} is solved using PCG with a Kronecker mass solver as preconditioner.

 \item[$\bullet$] A second flow $(\arr{E}^{n+1}, \arr{B}^{n+1}, \arr{Y}^{n+1}) = \Phi_{\arr{B},\arr{Y}}^{\mathrm{H}}(\arr{E}^n, \arr{B}^n, \arr{Y}^n; \Delta t, t^n)$ obtained by integrating the magnetic-plasma subsystem \eqref{Hamiltonian_splitting_operatorBC}:
\begin{subequations}
\label{integration_Hamiltonian_splitting_operatorBC}
 \begin{align}
 &\begin{bmatrix}
  \matM_1 + \cfrac{\Delta t}{2} \matA_{1} & \cfrac{\Delta t}{2} \matM_{1, \nomegap}\\
  \matZ & \matM_1 + \cfrac{\Delta t}{2} \left( \matR_{1, \nomegac} + \matM_{1, \hat{\nu}_e} \right)
 \end{bmatrix}
   \begin{bmatrix}  \arr{E}^{n+1} \\[1em] \arr{Y}^{n+1} \end{bmatrix} \nonumber\\
   &=
 \begin{bmatrix}
  \matM_1 \arrE^n - \cfrac{\Delta t}{2} \matA_{1} \arrE^n + \Delta t \matC^T \matM_2 \arr{B}^n - \cfrac{\Delta t}{2} \matM_{1, \nomegap} \arr{Y}^n + s^n \arr{S}_R^\inc - c^n \arr{S}_I^\inc \\
   \matM_1 \arrY^n - \cfrac{\Delta t}{2} \left( \matR_{1, \nomegac} + \matM_{1, \hat{\nu}_e} \right) \arrY^n
 \end{bmatrix}, \label{integration_Hamiltonian_splitting_operatorBC:EY} \\
 \arr{B}^{n+1} &= \arr{B}^n,
\end{align}
\end{subequations}
where $s^n = \sin(t^{n+1}) - \sin(t^n)$ and $c^n = \cos(t^{n+1}) - \cos(t^n)$. Here we used a second-order trapezoidal rule for the integration. The system \eqref{integration_Hamiltonian_splitting_operatorBC:EY} is solved with PBiCGStab, since the matrix is not symmetric. Again, the preconditioner is a block diagonal matrix, where the blocks are Kronecker mass solvers.
\end{itemize}

\subsubsection{Hamiltonian splitting scheme}
Lastly, we compose the flows using Strang splitting. Given an initial condition $\arrE^0, \arrB^0, \arrY^0$ 
at $t^0=0$, for each time-step $n=0,...,N-1$, compute
\begin{subequations}
\label{scheme_Hamiltonian}
\begin{align}
  (\arr{E}^n, \arr{B}^{n+1/2}, \arr{Y}^{n+1/3}) &= \Phi_{\arr{E}}^{\mathrm{H}}(\arr{E}^n, \arr{B}^n, \arr{Y}^n; \Delta t/2), \\
  (\arr{E}^{n+1}, \arr{B}^{n+1/2}, \arr{Y}^{n+2/3}) &= \Phi_{\arr{B},\arr{Y}}^{\mathrm{H}}(\arr{E}^{n}, \arr{B}^{n+1/2}, \arr{Y}^{n+1/3}; \Delta t, t^n), \\
  (\arr{E}^{n+1}, \arr{B}^{n+1}, \arr{Y}^{n+1}) &= \Phi_{\arr{E}}^{\mathrm{H}}(\arr{E}^{n+1}, \arr{B}^{n+1/2}, \arr{Y}^{n+2/3}; \Delta t/2),  \\
  t^{n+1} &= t^n + \Delta t,
\end{align}
\end{subequations}
where again the fractional superscripts just indicate intermediate solutions and not proper approximate solutions for intermediate times.

\subsection{Crank-Nicolson time scheme}

For the purpose of comparison we finally consider an implicit numerical integration of the full system
using a Crank-Nicolson method, which consists in integrating \eqref{semidiscrete_form} with the second-order trapezoidal rule. 
This yields a flow $(\arr{E}^{n+1}, \arr{B}^{n+1}, \arr{Y}^{n+1}) = \Phi^{\mathrm{CN}}(\arr{E}^n, \arr{B}^n, \arr{Y}^n; \Delta t, t^n)$ defined by
\begin{equation}
\begin{aligned}
  &\begin{bmatrix}
  \matM_1 + \cfrac{\Delta t}{2} \matA_{1} & -\cfrac{\Delta t}{2} \matC^T \matM_2 & \cfrac{\Delta t}{2} \matM_{1, \hat{\omega}_p}\\
  \cfrac{\Delta t}{2} \matC & \matI & \matZ\\
  -\cfrac{\Delta t}{2} \matM_{1, \hat{\omega}_p} & \matZ & \matM_1 + \cfrac{\Delta t}{2} (\matR_{1, \nomegac} + \matM_{1, \hat{\nu}_e})
  \end{bmatrix}
  \begin{bmatrix}  \arr{E}^{n+1} \\ \arr{B}^{n+1} \\ \arr{Y}^{n+1} \end{bmatrix} \nonumber
  \\
  &=
  \begin{bmatrix}
   \matM_1 \arr{E}^{n} - \cfrac{\Delta t}{2} \matA_{1} \arr{E}^{n} + \cfrac{\Delta t}{2} \matC^T \matM_2 \arr{B}^{n} -\cfrac{\Delta t}{2} \matM_{1, \hat{\omega}_p} \arr{Y}^{n} + s^n \arr{S}^\inc_R - c^n \arr{S}^\inc_I \\
  -\cfrac{\Delta t}{2} \matC \arr{E}^{n} + \arr{B}^{n} \\
  \cfrac{\Delta t}{2} \matM_{1, \hat{\omega}_p} \arr{E}^{n} + \matM_1 \arr{Y}^{n} - \cfrac{\Delta t}{2} \left( \matR_{1, \nomegac} + \matM_{1, \hat{\nu}_e} \right) \arr{Y}^{n}
  \end{bmatrix},
  \label{integration_Crank_Nicolson}
\end{aligned}
\end{equation}
where $\matI$ denotes the identity matrix, $s^n = \sin(t^{n+1}) - \sin(t^n)$ and $c^n = \cos(t^{n+1}) - \cos(t^n)$ 
as above.
Note that in the ideal case ($\mat{N} = \matZ$, $\arr{f}=\arr{0}$), this flow also preserves the Hamiltonian 
$\mathsf{H}$. Since the matrix is not symmetric, the system is solved using PBiCGStab. For the preconditioner we construct a block diagonal matrix, where each block is a Kronecker mass solver.

\subsubsection{Crank-Nicolson scheme}
Given an initial condition $\arrE^0, \arrB^0, \arrY^0$ at $t^0=0$, for each time-step $n=0,...,N-1$, compute
\begin{equation}
\label{scheme_Crank_Nicolson}
\begin{aligned}
  (\arr{E}^{n+1}, \arr{B}^{n+1}, \arr{Y}^{n+1}) &=  \Phi^{\mathrm{CN}}(\arr{E}^n, \arr{B}^n, \arr{Y}^n; \Delta t, t^n),  \\
  t^{n+1} &= t^n + \Delta t.
\end{aligned}
\end{equation}

\section{Long-time stability}
In this section we show that the solution cannot grow infinitely large after long times. We study the long-time stability at three levels: the continuous case, the semi-discrete case and the fully discrete case.

\subsection{Continuous level}
Here we show the long-time stability of the solution of the time-domain problem \eqref{weak_form} in terms of the evolution of the distance to the time-harmonic solution \eqref{cont_timeharmonicsol}, which involves solving the frequency-domain problem \eqref{freq_problem}.

The corresponding weak form consists of finding $\hat{\bE} \in H_\mathrm{imp}(\curl, \Omega)$ such that
\begin{equation}
 \langle \curl \bs{F}, \curl \hat{\bE} \rangle_\Omega - \langle  \bs{F}, \dieltens \hat{\bE} \rangle_\Omega - i \langle \nun \times \bs{F}, \nun \times \hat{\bE} \rangle_{\Gamma_A} = -i \langle \nun \times \bs{F}, \nun \times \hat{\bs{s}}^\inc \rangle_{\Gamma_A} \quad \forall \bs{F} \in H_\mathrm{imp}(\curl, \Omega). \label{weak_form_freq_pb}
\end{equation}
In the frequency-domain the underlying field of the spaces is $\mathbb{C}$, but for simplicity we maintain the same notation.

Given $\hat{\bE} \in H_\mathrm{imp}(\curl, \Omega)$ the solution to \eqref{weak_form_freq_pb}, the triplet $\bs{U}^{\mathrm{th}} = (\bE^{\mathrm{th}}, \bB^{\mathrm{th}}, \bY^{\mathrm{th}})$ defined in \eqref{cont_timeharmonicsol} with the fields $\hat{\bB} = -i \curl \hat{\bE}$ and $\nomegap \hat{\bY} = i(\hat{\bE} - \dieltens \hat{\bE})$
is a time-harmonic solution of the time-domain problem \eqref{weak_form} with initial condition $ \bs{U}^{\mathrm{th}}(0)$. As a result, $(\bE-\bE^{\mathrm{th}}, \bB-\bB^{\mathrm{th}}, \bY-\bY^{\mathrm{th}})$ is a solution of \eqref{weak_form} with $\bs{s}^\inc=\bs{0}$ and thus, from the energy balance  \eqref{deriv_Hamiltonian_continuous} we can deduce that
$\partial_t \mathcal{H}(\bE-\bE^{\mathrm{th}}, \bB-\bB^{\mathrm{th}}, \bY-\bY^{\mathrm{th}}) \le 0$, i.e.
\begin{equation}
\begin{aligned}
 \| \bs{U}(t) -\bs{U}^{\mathrm{th}}(t) \|_{L^2(\Omega)}
  \le  \| \bs{U}(0)-\bs{U}^{\mathrm{th}}(0) \|_{L^2(\Omega)} \qquad \forall t>0,
 \end{aligned}
\end{equation}
where $\| \bs{U} \|^2_{L^2(\Omega)} = \| \bE \|^2_{L^2(\Omega)} + \| \bB \|^2_{L^2(\Omega)} + \| \bY \|^2_{L^2(\Omega)} $.

\subsection{Semi-discrete level}
The long-time stability also holds after applying the space discretization. In terms of the spline coefficients the weak form of the frequency-domain problem \eqref{weak_form_freq_pb} becomes
\begin{equation}
 ( \matC^T \matM_1 \matC - \matM_{1,\dieltens} - i \matA_1 ) \hat{\arrE} = - i \hat{\arr{S}}^\inc, \label{semidiscrete_freqdomain}
\end{equation}
where
$$
(\matM_{1,\dieltens})_{\bs{i}, \bs{j}} = \langle \bs{\Lambda}^1_{\bs{i}}, \dieltens \bs{\Lambda}^1_{\bs{j}} \rangle_\Omega, \qquad (\hat{\arr{S}}^\inc)_{\bs{i}} = \langle \nun \times  \bs{\Lambda}^1_{\bs{i}} , \nun \times \hat{\bs{s}}^\inc \rangle_{\Gamma_A}.
$$
Again, the coefficient array $\arr{U}^{\mathrm{th}} = [\arrE^{\mathrm{th}}, \arrB^{\mathrm{th}}, \arrY^{\mathrm{th}}]^T$ with
\begin{equation*}
\arrE^{\mathrm{th}}(t) = \Re \{ \hat{\arrE} e^{-it}\}, \quad \arrB^{\mathrm{th}}(t) = \Re \{ \hat{\arrB} e^{-it}\}, \quad \arrY^{\mathrm{th}}(t) = \Re \{ \hat{\arrY} e^{-it}\}
\end{equation*}
is a solution of \eqref{semidiscrete_form} and therefore $[\arrE - \arrE^{\mathrm{th}}, \arrB - \arrB^{\mathrm{th}}, \arrY - \arrY^{\mathrm{th}}]^T$ is also a solution with $\arr{S}^\inc = \arr{0}$. From \eqref{deriv_Hamiltonian_discrete} we deduce that $\partial_t \mathsf{H}(\arrE - \arrE^{\mathrm{th}}, \arrB - \arrB^{\mathrm{th}}, \arrY - \arrY^{\mathrm{th}}) \le 0$, i.e. for every $t>0$,
\begin{equation}
\begin{aligned}
 \| \bs{U}_h(t)-\bs{U}_h^{\mathrm{th}}(t) \|_{L^2(\Omega)} \le  \| \bs{U}_h(0)-\bs{U}_h^{\mathrm{th}}(0) \|_{L^2(\Omega)} \qquad \forall t>0,
 \end{aligned} \label{semidiscrete_longtime_stability}
\end{equation}
where $\bs{\sigma}^1 (\bE_h)  = \arrE,
\bs{\sigma}^2 (\bB_h) = \arrB, \bs{\sigma}^1 (\bY_h) = \arrY$ are the coefficients of the solution $\bs{U}_h = (\bE_h, \bB_h, \bY_h)$.

\subsection{Fully discrete level}

At the fully-discrete level we show the long-time stability only for the Crank-Nicolson \eqref{scheme_Crank_Nicolson} and the Poisson splitting \eqref{scheme_Poisson} schemes. For the Hamiltonian-splitting scheme \eqref{scheme_Hamiltonian}, a similar result could be obtained by performing a backward error analysis and exhibiting a conserved modified Hamiltonian that is preserved and controls the true Hamiltonian under a CFL condition as it is done in \cite{DASILVA201524}. This analysis will not be performed here.

We derive a similar result to \eqref{semidiscrete_longtime_stability}, but with $\bs{U}_h$ being a solution of the Crank-Nicolson or the Poisson splitting scheme. Given a solution $\bs{U}_h = (\bs{E}_h, \bs{B}_h, \bs{Y}_h)$ with spline coefficients $\arr{U} = [\arrE, \arrB, \arrY]^T = [\bs{\sigma}_1(\bs{E}_h), \bs{\sigma}_2(\bs{B}_h), \bs{\sigma}_1(\bs{Y}_h)]^T$, we define the norm
\begin{equation*}
 \| \arr{U}^n\|_\matM = \sqrt{(\arr{U}^n)^T \matM \arr{U}^n} = \|\bs{U}_h(t^n)\|_{L^2(\Omega)} , \qquad \text{with} \quad \mat{M} =
 \begin{bmatrix}
  \matM_1 & 0 & 0 \\
  0 & \matM_2 & 0 \\
  0 & 0 & \matM_1
 \end{bmatrix},
\end{equation*}
for every $\arr{U}^n = \arr{U}(t^n) \in C_h$ and $n \ge 0$.

For the following analysis, we assume that the time-domain problem \eqref{weak_form} and the frequency-domain problem \eqref{weak_form_freq_pb} are well-posed.

\subsubsection{Crank-Nicolson scheme}

We recall the full evolution equation for the Crank-Nicolson scheme given in \eqref{semidiscrete_form}, that is,
\begin{align}
 \partial_t \arr{U} &= (\mat{P} + \mat{N}) \mat{M} \arr{U} + \arr{f}(t).
\label{semidiscrete_evol_CN}
\end{align}
The Crank-Nicolson scheme consists of the update rule
\begin{align}
 \arr{U}^{n+1} - \arr{U}^n = \cfrac{\Delta t}{2} (\mat{P} + \mat{N}) \mat{M} (\arr{U}^{n+1} + \arr{U}^n) + \int_{t^n}^{t^n + \Delta t} \arr{f}(t) dt,
\end{align}
In the case $\arr{f}=\arr{0}$, we can define the discrete evolution operator $\mat{K}$ satisfying  $\arr{U}^{n+1} = \mat{K} \arr{U}^n$ for every $n$, that is,
\begin{equation}
 \mat{K} = \left(\matI - \cfrac{\Delta t}{2} (\mat{P} + \mat{N}) \mat{M} \right)^{-1}\left(\matI + \cfrac{\Delta t}{2} (\mat{P} + \mat{N}) \mat{M} \right).
 \label{discrete_evolution_operator_CN}
\end{equation}
The discrete evolution operator is well defined since the matrix $\matI - (\Delta t/2) (\mat{P} + \mat{N}) \mat{M} $ is invertible.
\begin{proof}
Let $\arr{W}_k$ be an eigenvector with eigenvalue $\lambda_k$, then
\begin{align*}
 &\lambda_k \arr{W}_k^T \mat{M} \arr{W}_k = \left[ \left( \matI - (\Delta t/2) (\mat{P} + \mat{N}) \mat{M} \right) \arr{W}_k \right]^T \mat{M} \arr{W}_k =  \arr{W}_k^T \mat{M} \arr{W}_k - (\Delta t/2) \left[  (\mat{P} + \mat{N}) \mat{M}  \arr{W}_k \right]^T \mat{M} \arr{W}_k\\
 &\Rightarrow  (1 - \lambda_k) \arr{W}_k^T \mat{M} \arr{W}_k = (\Delta t/2) \left(  \mat{N} \mat{M}  \arr{W}_k \right)^T \mat{M} \arr{W}_k \le 0 \quad \Rightarrow \lambda_k \ge 1.
\end{align*}
In particular, $\lambda_k \neq 0 $, thus $\matI - (\Delta t/2) (\mat{P} + \mat{N}) \mat{M} $ is invertible.
\end{proof}

Furthermore, $\mat{K}$ is nonexpansive, in the sense that $\|\mat{K}\arr{U}\|_\matM \le \|\arr{U}\|_\matM$.
\begin{proof}
 Assume $\arr{f}=\arr{0}$, then $\arr{U}^{n+1} = \mat{K} \arr{U}^n$ and it holds that
 \begin{align*}
 \|\mat{K}\arr{U}^n\|^2_\matM - \|\arr{U}^n\|^2_\matM &= (\arr{U}^{n+1} - \arr{U}^n)^T \matM (\arr{U}^{n+1} - \arr{U}^n) = \cfrac{\Delta t}{2} \left [ (\mat{P} + \mat{N}) \matM (\arr{U}^{n+1} + \arr{U}^n) \right]^T  \matM (\arr{U}^{n+1} - \arr{U}^n)  \\
 &= \cfrac{\Delta t}{2} \left [ \mat{N} \matM (\arr{U}^{n+1} + \arr{U}^n) \right]^T \matM (\arr{U}^{n+1} - \arr{U}^n)  \le 0,
\end{align*}
since $\mat{P}$ is skew-symmetric and $\mat{N}$ is negative semi-definite.
\end{proof}

Next we consider the source term $\arr{f}$ and show the long-time stability result.
We consider the coefficients corresponding to a time-harmonic solution given by $\arr{U}^{\mathrm{th},n} = \Re \{ \hat{\arr{U}}_{\Delta t} e^{-in\Delta t} \}$, where the complex $\hat{\arr{U}}_{\Delta t}$ depends on $\Delta t$ but not on $n$. The corresponding sequence in $n$ yields a discrete evolution of the Crank-Nicolson scheme if for any $n$ it holds
\begin{equation}
 \arr{U}^{\mathrm{th},n} = \mat{K} \arr{U}^{\mathrm{th},n-1} + \left(\matI - \cfrac{\Delta t}{2} (\mat{P} + \mat{N}) \matM \right)^{-1} \int_{(n-1)\Delta t}^{n\Delta t} \arr{f} dt. \label{discrete_evol_CN}
\end{equation}
Equation \eqref{discrete_evol_CN} forces the definition of $\hat{\arr{U}}_{\Delta t}$, which reads
\begin{equation}
 \hat{\arr{U}}_{\Delta t} = i(e^{-i\Delta t} - 1) \left( e^{-i\Delta t} \matI - \mat{K} \right)^{-1} \left(\matI - \cfrac{\Delta t}{2} (\mat{P} + \mat{N}) \matM \right)^{-1} \hat{\arr{f}},
\label{coeffs_timeharmonic_CN}
\end{equation}
where $\hat{\arr{f}} = [\matM_1^{-1} \langle \nun \times  \bs{\Lambda}^1_{\bs{i}} , \nun \times \hat{\bs{s}}^\inc \rangle_{\Gamma_A}, \bs{0}, \bs{0}]^T$.
Here we assume that the matrix $(e^{-i\Delta t} \matI - \mat{K})$ is invertible, which is a natural consideration following the assumption that the frequency-domain problem is well-posed.

With this construction, for any iterates of the Crank-Nicolson scheme $\arr{U}^n$, the difference $\arr{U}^n - \arr{U}^{\mathrm{th},n}$ solves a  discrete evolution with $\arr{f}=\arr{0}$, i.e.
\begin{equation}
 \arr{U}^n - \arr{U}^{\mathrm{th},n} = \mat{K}^n \left( \arr{U}^0 - \arr{U}^{\mathrm{th},0} \right).
\end{equation}
Since $\mat{K}$ is nonexpansive, we conclude with the long-time stability result
\begin{equation}
 \| \bs{U}_h(t^n) - \bs{U}_h^\mathrm{th}(t^n)\|_{L^2(\Omega)} \le  \| \bs{U}_h(0) - \bs{U}_h^\mathrm{th}(0)\|_{L^2(\Omega)}\qquad \text{for every } n \ge 0,
 \label{long_time_stability_CN}
\end{equation}
where $\arr{U}^{\mathrm{th},n} $
are the spline coefficients of $\bs{U}_h^\mathrm{th}(t^n)$.

\subsubsection{Poisson splitting scheme}
First we recall the subproblems of the Poisson splitting scheme defined in \eqref{Poisson_splitting_operatorA} and  \eqref{Poisson_splitting_operatorBC}, that is,
\begin{equation}
\begin{aligned}
 \label{semidiscrete_evol_Poisson}
\partial_t \arr{U} &= \mat{T}_1 \arr{U} + \arr{f}(t),\\
\partial_t \arr{U} &= \mat{T}_2 \arr{U},
\end{aligned}
\end{equation}
with $\mat{T}_1 = (\mat{P}_1 + \mat{N}_1) \matM$, $\mat{T}_2 = (\mat{P}_2 + \mat{N}_2) \matM$,
$$
\mat{N}_1 = \begin{bmatrix}
            - \matM_1^{-1} \matA_1 \matM_1^{-1} & \matZ &  \matZ\\
            \matZ & \matZ &  \matZ\\
            \matZ & \matZ &  \matZ
            \end{bmatrix} \qquad \text{and} \qquad
\mat{N}_2 = \begin{bmatrix}
            \matZ & \matZ &  \matZ\\
            \matZ & \matZ &  \matZ\\
            \matZ & \matZ &  - \matM_1^{-1} \matM_{1, \hat{\nu}_e} \matM_1^{-1}
            \end{bmatrix}.
$$
The Poisson splitting scheme consists of the update rule
\begin{equation}
\begin{aligned}
 \arr{U}^{n+1/3} - \arr{U}^n &= \cfrac{\Delta t}{4} \mat{T}_1 (\arr{U}^{n+1/3} + \arr{U}^n) + \int_{t^n}^{t^n+\Delta t/2} \arr{f} dt,\\
 \arr{U}^{n+2/3} - \arr{U}^{n+1/3} &= \cfrac{\Delta t}{2} \mat{T}_2 (\arr{U}^{n+2/3} + \arr{U}^{n+1/3}),\\
 \arr{U}^{n+1} - \arr{U}^{n+2/3} &= \cfrac{\Delta t}{4} \mat{T}_1 (\arr{U}^{n+1} + \arr{U}^{n+2/3}) + \int_{t^n+\Delta t/2}^{t^n+\Delta t} \arr{f} dt.\\
\end{aligned}
\end{equation}
Hence any discrete evolution of the Poisson splitting scheme \eqref{scheme_Poisson} must satisfy
\begin{equation}
 \arr{U}^{n+1} = \mat{K} \arr{U}^n + \left(\matI - \cfrac{\Delta t}{4} \mat{T}_1 \right)^{-1} \int_{t^n+\Delta t/2}^{t^{n}+ \Delta t} \arr{f} dt + \mat{K} \left(\matI + \cfrac{\Delta t}{4} \mat{T}_1 \right)^{-1} \int_{t^n}^{t^n+\Delta t/2} \arr{f} dt,
\end{equation}
where the discrete evolution operator is given by
\begin{equation}
\mat{K} = \left(\matI - \cfrac{\Delta t}{4} \mat{T}_1 \right)^{-1}\left(\matI + \cfrac{\Delta t}{4} \mat{T}_1 \right) \left(\matI - \cfrac{\Delta t}{2} \mat{T}_2 \right)^{-1}\left(\matI + \cfrac{\Delta t}{2} \mat{T}_2 \right) \left(\matI - \cfrac{\Delta t}{4} \mat{T}_1 \right)^{-1}\left(\matI + \cfrac{\Delta t}{4} \mat{T}_1 \right),
\end{equation}
where the matrices $\matI - (\Delta t / 4) \mat{T}_1$ and $\matI - (\Delta t / 2) \mat{T}_2$ are invertible. The proof is analogous to the proof for the inverse of $\matI - (\Delta t/2) (\mat{P} + \mat{N}) \mat{M} $ in \eqref{discrete_evolution_operator_CN}.

In addition, $\mat{K}$ is nonexpansive, in the sense that $\|\mat{K}\arr{U}\|_\matM \le \|\arr{U}\|_\matM$ for every $\arr{U} \in C_h$. This is due to the fact that the operator is the composition of three nonexpansive operators (for each step). The proof that each step is nonexpansive is analogous to the proof showing that \eqref{discrete_evolution_operator_CN} is nonexpansive.

Following we show the form of the solution considering the source. As before, we consider the coefficient sequence given by $\arr{U}^{\mathrm{th},n} = \Re \{ \hat{\arr{U}}_{\Delta t} e^{-in\Delta t} \}$, where the complex $\hat{\arr{U}}_{\Delta t}$ depends on $\Delta t$ but not on $n$. They yield a discrete evolution of the Poisson splitting scheme \eqref{scheme_Poisson} if for every $n$ it holds
\begin{equation}
  \hat{\arr{U}}_{\Delta t} = i \left( e^{-i\Delta t} \matI - \mat{K} \right)^{-1}
  \left [ (e^{-i\Delta t} - e^{-i\Delta t/2}) \left(\matI - \cfrac{\Delta t}{4} \mat{T}_1 \right)^{-1} +
  (e^{-i\Delta t/2} - 1) \mat{K} \left(\matI + \cfrac{\Delta t}{4} \mat{T}_1 \right)^{-1} \right]
 \hat{\arr{f}}.
\label{coeffs_timeharmonic_Poisson}
\end{equation}
As in \eqref{coeffs_timeharmonic_CN}, here we assume that $( e^{-i\Delta t} \matI - \mat{K})$ is invertible.
For any discrete evolution $\arr{U}^n$, the difference yields an evolution of the Poisson splitting scheme \eqref{scheme_Poisson} without source. Consequently, it holds that
\begin{equation}
 \arr{U}^n - \arr{U}^{\mathrm{th},n} = \mat{K}^n \left( \arr{U}^0 - \arr{U}^{\mathrm{th},0} \right),
 \label{long_time_stability_coeffs_Poisson}
\end{equation}
which shows the same long-time stability result as in \eqref{long_time_stability_CN}, since $\mat{K}$ is also nonexpansive.

\section{Benchmarking and performance study}
\label{section_1Dtests}
In this section we verify the accuracy of our scheme on exact univariate solutions. For that we consider the domain $\Omega = [0,3\pi] \times [0,2\pi] \times [0, 2\pi]$ and the physical parameters $\nomegac=0.5$, $\nb_0 = (0,0,1)$ and $\nomegap(\bx)= x/100$. In addition, periodic boundary conditions are imposed along the $y$ and $z$ directions.
We study the accuracy of the proposed methods in two test cases. The first one follows an O-mode polarization (Transverse Magnetic), i.e. $\bE \parallel \nb_0$, and the second one follows an X-mode polarization (Transverse Electric), i.e. $\bE \perp \nb_0$.
The simulation of an X-mode wave is in general more complex than an O-mode wave. Firstly, in \eqref{weak_form:current_density_eq} the term $\langle \bs{G} \times \bY, \nomegac \nb_0  \rangle_\Omega$ is zero in the O-mode case. Secondly, as we have shown in the dispersion relation analysis of section \ref{section_freq_domain}, the X-mode configuration has a resonance. In our tests we tune the parameters in order to avoid the resonance. Since the X-mode case is more involved, we include a study on the numerical stability, conservation of physical quantities and performance.
A solution with a certain polarization can be manufactured following the steps:
\begin{enumerate}
 \item Choose a current density $\bY^\ex$ either polarized along $z$ (O-mode) or in the $x-y$ plane (X-mode).
 \item Find electric field $\bE^\ex$ using the current density equation \eqref{model:current_density_eq}.
 \item Find magnetic field $\bB^\ex$ using Faraday equation \eqref{model:faraday_eq}.
 \item Find the trace $\nun \times \bs{s}^\inc = \nun \times (\bE^\ex - \bB^\ex \times \bs{\nu}) $ so that the boundary conditions \eqref{SM_BC_general} hold.
 \item Add a volume source $\bs{S} := \partial_t \bE^{\mathrm{ex}} - \curl \bB^{\mathrm{ex}} + \hat{\omega}_p \bY^{\mathrm{ex}}$ to the Ampère-Maxwell equation \eqref{model:ampere_eq} corresponding to the residual.
\end{enumerate}
For convenience, we construct the above solution as the real part of a time-harmonic solution, that is 
$\bY^\ex = \Re\{\hat{\bY}^\ex e^{-it}\}$, $\dots$, $\bs{s}^\inc = \Re\{\hat{\bs{s}}^\inc e^{-it}\}$.
We further decompose the volume source as $\bs{S} = \bs{S}_R \cos(t) + \bs{S}_I \sin(t)$ and redefine the arrays in \eqref{SincRI} by
\begin{align}
(\arr{S}^\inc_R)_{\bs{i}} := \langle \nun \times \bs{\Lambda}^1_{\bs{i}}, \nun \times \Re \{ \hat{\bs{s}}^{\mathrm{inc}} \} \rangle_{\Gamma_A} + \langle \bs{\Lambda}^1_{\bs{i}}, \bs{S}_R \rangle_{\Omega},
& &
(\arr{S}^\inc_I)_{\bs{i}} := \langle \nun \times \bs{\Lambda}^1_{\bs{i}}, \nun \times \Im \{ \hat{\bs{s}}^{\mathrm{inc}} \} \rangle_{\Gamma_A} + \langle \bs{\Lambda}^1_{\bs{i}}, \bs{S}_I \rangle_{\Omega},
\label{sources_test}
\end{align}
so that the semi-discrete form in \eqref{semidiscrete_form} holds again, and the same schemes can be used.

\subsection{Numerical errors and discrete numerical parameters}

We denote the Courant number $\CFL=\Delta t / \Delta x$, where $\Delta t$ is the normalized time-step and $\Delta x$ is the diameter of a cell in the $x$ direction.
Each simulation is run for three periods, i.e. $t \in [0, 3T]$, where $T=2\pi$. The grid size and time-step are controlled by the parameters
$$ \mathrm{PPW} = \cfrac{2\pi}{\Delta x}, \qquad \mathrm{PPP} = \cfrac{2\pi}{\Delta t},$$
which correspond to the number of points per wavelength and the number of time-points per period, respectively.
Since the parameters and exact solution only depend on $x$ we use a single cell in the $y$ and $z$ directions.
The B-spline degree of the space $V^0_h$ is fixed as $\bs{p} =(3,1,1)$, thus the degree of the basis functions for spaces $V^1_h$ and $V^2_h$ is set, as it is shown in \eqref{Bspline_bases}.
We distinguish the following three errors:
 \begin{itemize}
    \item Projection errors
     defined as
   \begin{align}
    \varepsilon^{\mathrm{proj}}_{E} &= \|\bE^{\mathrm{ex}} - P_1(\bE^{\mathrm{ex}}) \|_{L^2(\Omega)}, & \varepsilon^{\mathrm{proj}}_{B} &= \|\bB^{\mathrm{ex}} - \Pi_2(\bB^{\mathrm{ex}}) \|_{L^2(\Omega)}, & \varepsilon^{\mathrm{proj}}_{Y} &= \|\bY^{\mathrm{ex}} - P_1(\bY^{\mathrm{ex}}) \|_{L^2(\Omega)}, \label{def_proj_error}
   \end{align}
   where $P_1:(L^2(\Omega))^3 \rightarrow V_h^1$ denotes the $L^2$-projection.  For consistency reasons with discrete Gauss laws, it is natural to define the reference discrete fields using the $L^2$-projection for $\bE^{\mathrm{ex}}$, $\bY^{\mathrm{ex}}$ and a commuting projection $\Pi_2$ for $\bB^{\mathrm{ex}}$.
   These errors only depend on the space discretization, i.e. $\mathrm{PPW}$ and $\bs{p}$.

  \item Total errors
  defined as
  \begin{align}
  \varepsilon^{\mathrm{total}}_{E} &= \|\bE^{\mathrm{ex}} - \bE_h \|_{L^2(\Omega)}, & \varepsilon^{\mathrm{total}}_{B} &= \|\bB^{\mathrm{ex}} - \bB_h \|_{L^2(\Omega)}, & \varepsilon^{\mathrm{total}}_{Y} &= \|\bY^{\mathrm{ex}} - \bY_h \|_{L^2(\Omega)}, \label{def_solution_error}
  \end{align}
which depend on the space and time discretizations, thus $\mathrm{PPW}$, $\bs{p}$ and $\mathrm{PPP}$.

  \item Errors of the time-evolution solver  defined as
  \begin{align}
  \varepsilon^{\mathrm{solver}}_{E} &= \|P_1(\bE^{\mathrm{ex}}) - \bE_h \|_{L^2(\Omega)}, & \varepsilon^{\mathrm{solver}}_{B} &= \|\Pi_2(\bB^{\mathrm{ex}}) - \bB_h \|_{L^2(\Omega)}, & \varepsilon^{\mathrm{solver}}_{Y} &= \|P_1(\bY^{\mathrm{ex}}) - \bY_h \|_{L^2(\Omega)}.
  \label{def_solver_error}
  \end{align}
Since the discrete operators are approximations of continuous operators, the solver errors also depend on the space and time discretizations, thus $\mathrm{PPW}$, $\bs{p}$ and $\mathrm{PPP}$.
 \end{itemize}
These errors are defined pointwise in time. In practice, we study their maximum value in time.

%%%%%%%%%%%%%%%%%%%%%%%%%%%%%%%%%%%%%%%%%%%%%%%%%%%%%%%%%%%%%%%%%%%%%%%%%%%%%%

\subsection{O-mode test}
Following the steps mentioned previously, we construct a solution with O-mode configuration, which reads
\begin{equation}
 \label{exact_solution_Omode}
 \begin{aligned}
\bE^{\mathrm{ex}} &=
 \begin{bmatrix}
  0\\
  0\\
  \cos (x) \cos (t) + \sin(x) \sin(t)
 \end{bmatrix}, &
\bB^{\mathrm{ex}} &= -
 \begin{bmatrix}
  0\\
  \cos (x) \cos (t) + \sin(x) \sin(t)\\
  0
 \end{bmatrix},\\
\bY^{\mathrm{ex}} &= \nomegap(x)
 \begin{bmatrix}
   0\\
   0\\
   \cos(x) \sin(t) - \sin (x) \cos (t)
  \end{bmatrix},
  \end{aligned}
 \end{equation}
for which the sources in \eqref{sources_test} are given by
\begin{equation}
 \bs{S}_R = \begin{bmatrix}
           0 \\
           0 \\
            - \nomegap^2 \sin(x)
          \end{bmatrix}, \quad
\bs{S}_I = \begin{bmatrix}
           0 \\
           0 \\
           \nomegap^2 \cos(x)
          \end{bmatrix}, \quad
 \bs{S}^\inc_R = \begin{bmatrix}
           0 \\
           0 \\
           (1 - \nu_1) \cos(x)
          \end{bmatrix}, \quad
 \bs{S}^\inc_I = \begin{bmatrix}
           0 \\
           0 \\
           (1 - \nu_1) \sin(x)
          \end{bmatrix}.
 \label{exact_sources_Omode}
\end{equation}

\subsubsection{Verification of the high order space approximation}

\begin{figure}[H]
         \centering %[trim=left bottom right top, clip]
         \includegraphics[width=0.45\textwidth, trim=0 0.9cm 0.5cm 3.15cm, clip]{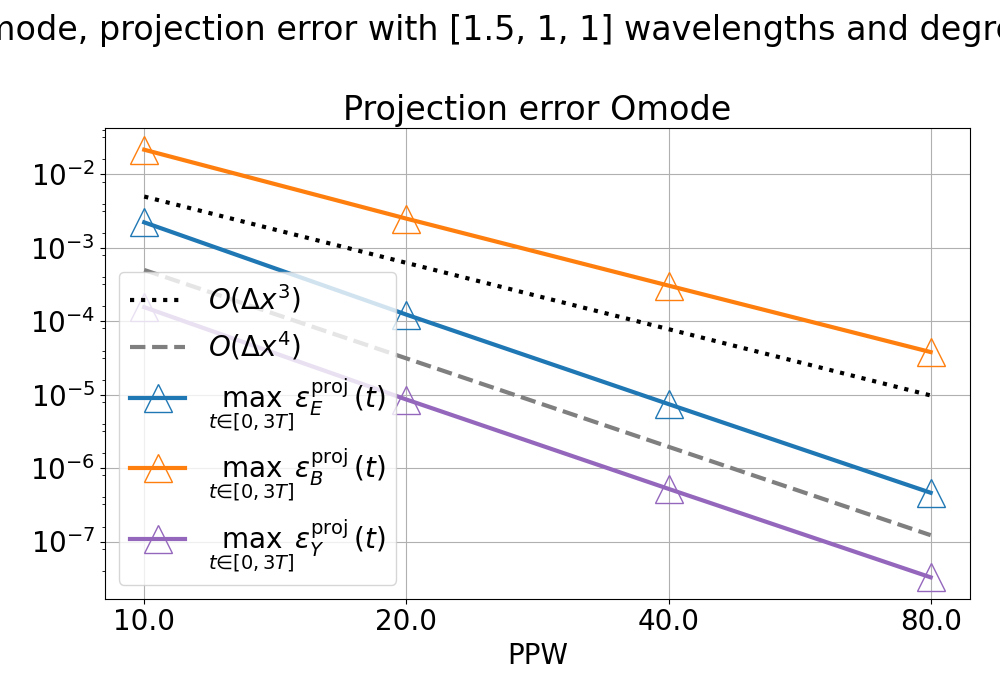}
         \includegraphics[width=0.45\textwidth, trim=0 0.9cm 0.5cm 3.15cm, clip]{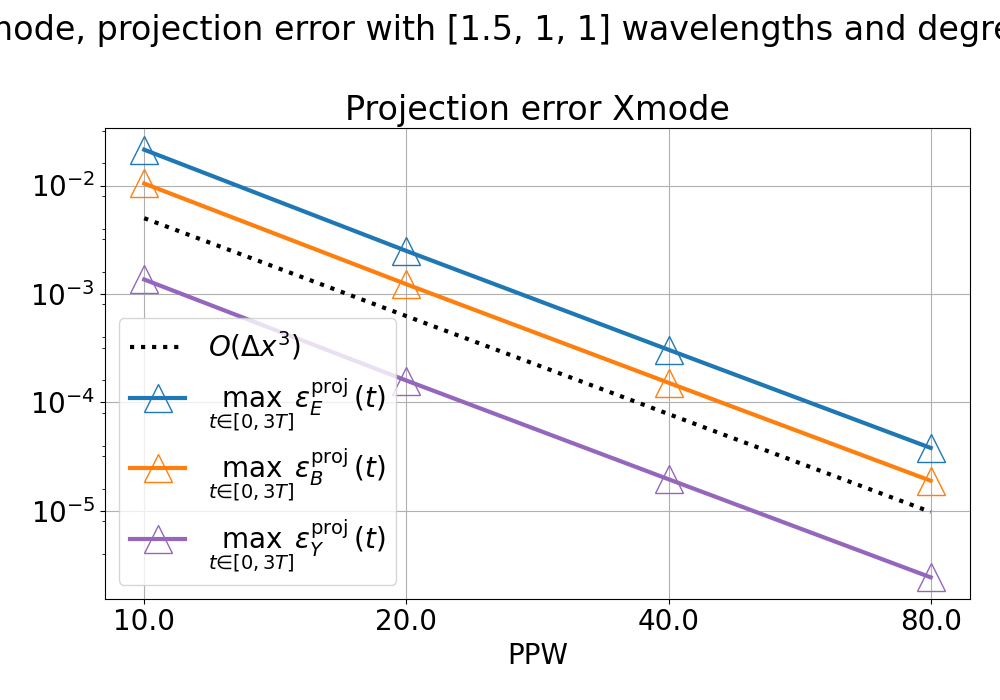}
         \caption{Projection errors of the O-mode solution \eqref{exact_solution_Omode} (left) and X-mode solution \eqref{exact_solution_Xmode} (right) for a decreasing grid size with B-spline degree $\bs{p}=(3,1,1)$ in $V_h^0$.
         }
         \label{pics_manufactured_proj_error}
 \end{figure}

Before studying the convergence of the solver, we verify the order of the projection errors for the manufactured solutions \eqref{exact_solution_Omode} and \eqref{exact_solution_Xmode}. In figure \ref{pics_manufactured_proj_error}, we see that the projection errors decrease at least as $O(\Delta x^3)$. This is expected given the definition of the discretized spaces in \eqref{Bspline_bases} and the fixed anisotropic B-spline degrees.

\subsubsection{Convergence}
\begin{figure}[H]
         \centering %[trim=left bottom right top, clip]
        \includegraphics[width=\textwidth, trim=0.8cm 0.8cm 1.5cm 2.1cm, clip]{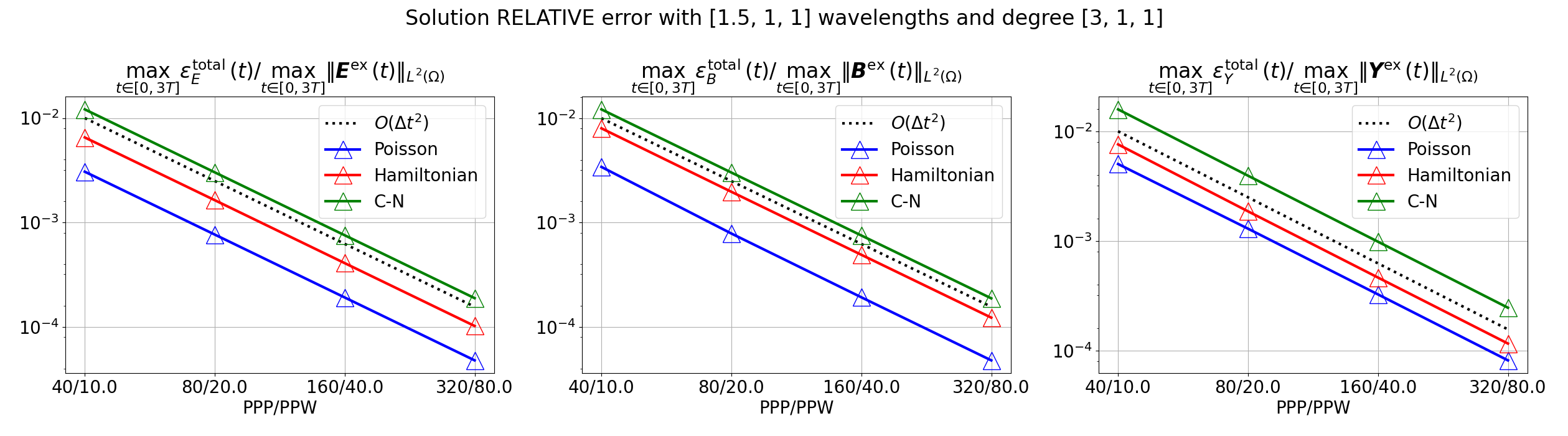}
         \caption{(O-mode solution) Decay of the relative total errors as the grid is refined with fixed $\CFL=0.25$.}
         \label{pics_manufactured_Omode_rel_errors}
 \end{figure}

Here we study the decay of the total errors as the grid is refined with constant $\CFL=0.25$. The error plots are gathered in figure \ref{pics_manufactured_Omode_rel_errors}.
 Figure \ref{pics_manufactured_Omode_rel_errors} shows that the three schemes are second order, since the decay of the total errors for consecutive refinements under constant $\CFL=0.25$ is quadratic, i.e. the total errors behave as $O(\Delta t^2)$. As shown in figure \ref{pics_manufactured_proj_error}, the projection errors decay faster than the total errors and therefore they are negligible.

 Finally, we observe that the Poisson splitting method produces the smallest total error with a remarkable difference with respect to the Crank-Nicolson method, which produces the largest error.

\subsection{X-mode test}
\label{sec:Xmodetest}
Here we construct an exact solution with X-mode configuration given by
\begin{equation}
 \bE^{\mathrm{ex}} =
 \begin{bmatrix}
  - \cos(x) \sin(t) \\
  - \hat{\omega}_c \cos(x) \cos(t) \\
  0
 \end{bmatrix},
\qquad
 \bB^{\mathrm{ex}} =
 \begin{bmatrix}
  0\\
  0\\
  - \hat{\omega}_c \sin(x) \sin(t)
 \end{bmatrix},
 \qquad
 \bY^{\mathrm{ex}} =
 \begin{bmatrix}
  \hat{\omega}_p(x) \cos(x) \cos(t)\\
  0\\
  0
 \end{bmatrix},
 \label{exact_solution_Xmode}
\end{equation}
for which the sources in \eqref{sources_test} are given by
\begin{equation}
 \bs{S}_R = \begin{bmatrix}
           (\nomegap^2(x) - 1) \cos(x)\\
           0 \\
           0
          \end{bmatrix}, \qquad
\bs{S}_I = \bs{0}, \qquad
 \bs{S}^\inc_R = \begin{bmatrix}
           0 \\
           - \nomegac \cos(x) \\
           0
          \end{bmatrix}, \qquad
 \bs{S}^\inc_I = \begin{bmatrix}
           - \cos(x) \\
           \nomegac \nu_1(x) \sin(x)\\
           0
          \end{bmatrix}.
 \label{exact_sources_Xmode}
\end{equation}

\subsubsection{Convergence}
\label{subsection_Xmode_convergence}
\begin{figure}[H]
         \centering %[trim=left bottom right top, clip]
         \includegraphics[width=\textwidth, trim=0.8cm 0.8cm 1.5cm 2.1cm, clip]{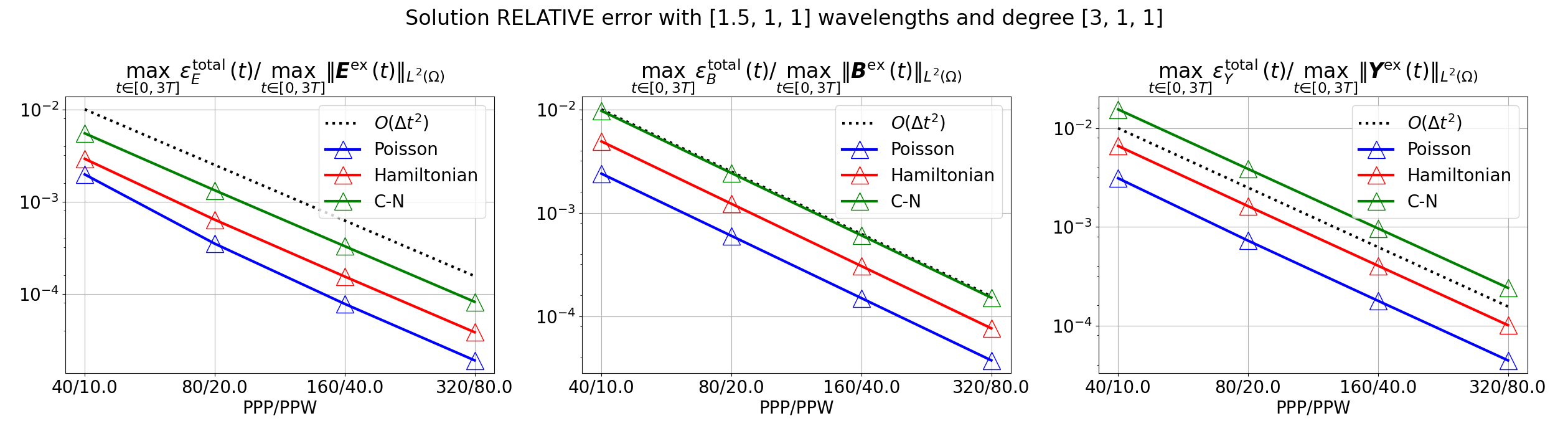}
         \includegraphics[width=\textwidth, trim=0.8cm 0.8cm 1.5cm 2.1cm, clip]{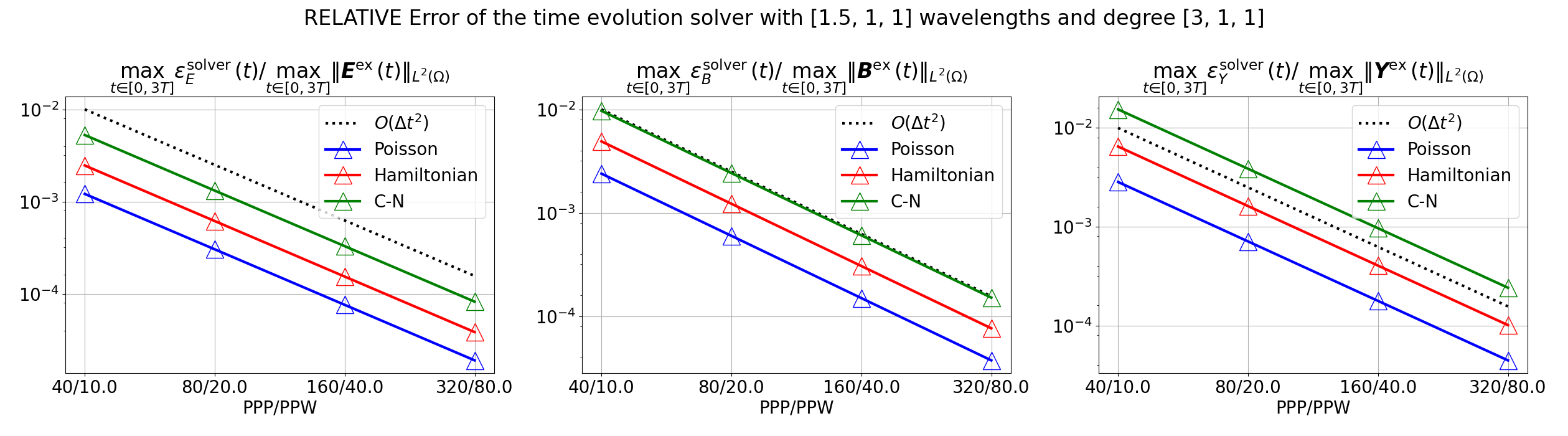}
         \caption{(X-mode solution) Decay of the relative total (top) and solver (bottom) errors as the grid is refined with fixed $\CFL=0.25$.}
         \label{pics_manufactured_rel_errors}
 \end{figure}

As before, we study the error decay as the grid is refined with constant $\CFL=0.25$. In particular, figure \ref{pics_manufactured_rel_errors} shows the decay of the relative total and solver errors.

Figure \ref{pics_manufactured_rel_errors} shows that the decay of the total and solver errors is quadratic as the grid is refined with fixed $\CFL=0.25$, which implies that the three methods are second-order. Furthermore, the solver error is slightly lower than the total error and the difference is only remarkable when the grid is coarse, thus when the projection errors are large.
Finally, we observe that the Poisson splitting scheme produces the smallest error with a large margin to the Crank-Nicolson method.

\subsubsection{Stability}
\begin{figure}[H]
         \centering %[trim=left bottom right top, clip]
         \includegraphics[width=\textwidth, trim=0.8cm 0.8cm 0.8cm 2.1cm, clip]{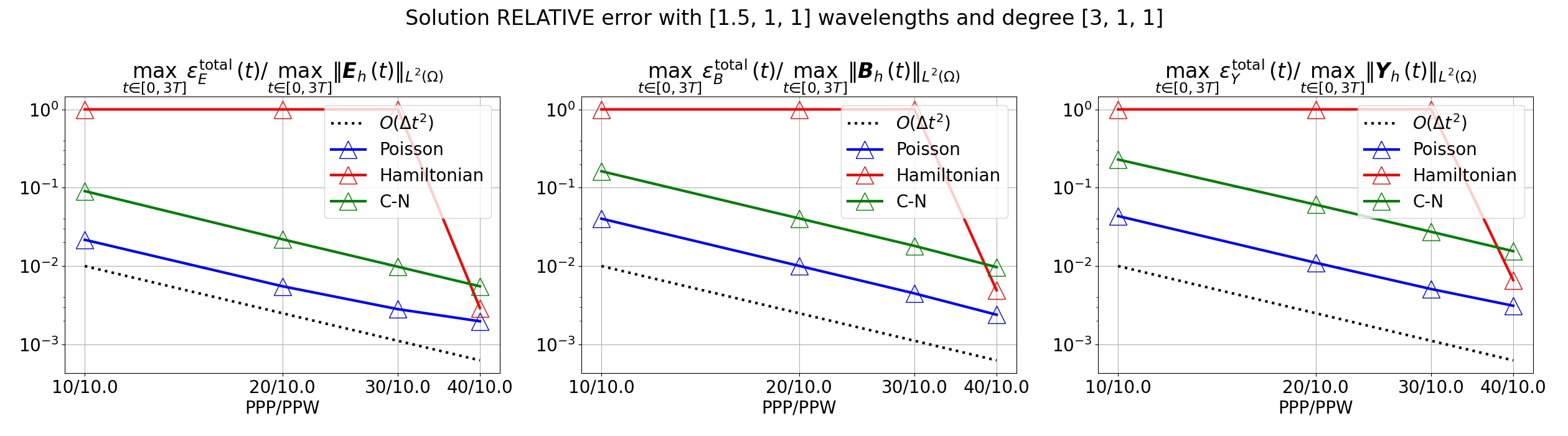}
         \caption{(X-mode solution) Relative total errors with fixed $\mathrm{PPW}=10$ and decreasing $\CFL$ (from left to right in the horizontal axis).}
         \label{pics_manufactured_solerror_fixedncells}
 \end{figure}

In this part we study the numerical stability in the sense of measuring the total errors as $\Delta t$ increases and $\Delta x$ is kept constant, thus as $\CFL$ becomes larger.
We recall that a method is conditionally stable, if it requires a $\CFL$-condition to converge. If it converges without a restriction on the ratio between $\Delta x$ and $\Delta x$, the method is unconditionally stable. In Figure \ref{pics_manufactured_solerror_fixedncells} we show the relative total errors as $\Delta t$ decreases and $\mathrm{PPW}=10$ is fixed.
Figure \ref{pics_manufactured_solerror_fixedncells} shows that the Hamiltonian splitting scheme is conditionally stable (with $\CFL \simeq 0.25$) and the Poisson splitting and Crank-Nicolson schemes are unconditionally stable, since they converge with second order under any $\CFL$ values. This numerical experiment confirms the long-time stability for Crank-Nicolson and Poisson splitting schemes (see \ref{long_time_stability_CN}).
For instance, with $\CFL = 0.33$, the solution of the Hamiltonian splitting method reaches the magnitude of $10^{30}$. For this reason, in the figure we plot the total error relative to the numerical solution.

\subsubsection{Conservation of physical quantities}
\begin{figure}[H]
         \centering %[trim=left bottom right top, clip]
         \includegraphics[width=\textwidth, trim=0.8cm 0.8cm 1.5cm 2.1cm, clip]{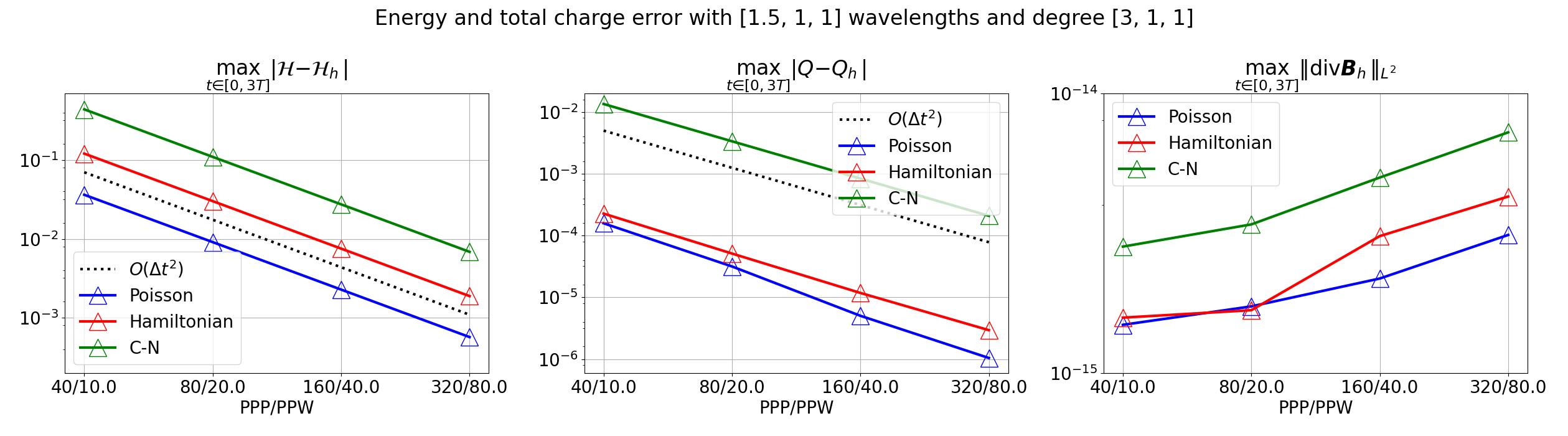}
         \caption{(X-mode solution) Error in the energy (left), total charge (center) and divergence of the magnetic field (right) as the grid is refined with fixed $\CFL=0.25$. In the left plot $\mathcal{H}$ denotes $\mathcal{H}(\bE^\mathrm{ex}, \bB^\mathrm{ex}, \bY^\mathrm{ex})$ and $\mathcal{H}_h$ refers to $\mathcal{H}(\bE_h, \bB_h, \bY_h)$.
         }
         \label{pics_manufactured_energy_charge_conv}
 \end{figure}

Here we study the error of the schemes at approximating the following quantities of interest:
\begin{itemize}
 \item The discrete Hamiltonian $\mathcal{H}(\bE_h, \bB_h, \bY_h)$, which is compared to the exact Hamiltonian given by
\begin{equation*}
  \mathcal{H}(\bE^\mathrm{ex}, \bB^\mathrm{ex}, \bY^\mathrm{ex}) =
   \pi^2 (\hat{\omega}_c^2 + 1)L \sin^2(t)
  + \pi^2 \cos^2(t) \left (  \hat{\omega}_c^2 L  + \cfrac{4L^3 + (6L^2-3) }{12 \times 10^4} \right ), \qquad \text{ where }L=3\pi.
\end{equation*}
 \item The total charge $Q_h(t) = \int_{\Omega} \widetilde{\mathrm{div}}_h \bE_h d\bx$, which is compared to the exact total charge $Q(t) =  \int_{\Omega} \mathrm{div} \bE^\mathrm{ex} d\bx = 2\pi \sin(t)$. The weak divergence operator $\widetilde{\mathrm{div}}_h: V_h^1 \rightarrow V_h^0$ is defined as
\begin{align}
 \langle \Lambda^0_{\bs{i}} , \widetilde{\mathrm{div}}_h \bs{\Lambda}^1_{\bs{j}} \rangle_{\Omega} &= - \langle  \mathrm{grad} \Lambda^0_{\bs{i}} , \bs{\Lambda}^1_{\bs{j}} \rangle_{\Omega} +  \langle \Lambda^0_{\bs{i}}, \bs{\Lambda}^1_{\bs{j}} \cdot \bs{\nu} \rangle_{\partial \Omega} = (- \matG^T\matM_1 + \matB_1)_{\bs{i}, \bs{j}},
 \label{weak_div}
 \end{align}
 where the matrix $\matB_1$ is defined by the boundary term. Then we can compute $Q_h(t) = \sum_{\bs{i}}  (- \matG^T\matM_1 + \matB_1)_{\bs{i}} \arrE$.

 \item The divergence of the magnetic field $\mathrm{div} \bB_h$, which should remain zero as it is an invariant \eqref{invariant_divB}. Following the de Rham diagram \eqref{deRham_diagram}, it can be computed exactly as $\bs{\sigma}^3(\mathrm{div} \bB_h) = \matD \arrB$, since $\bB_h \in V^2_h \subset H(\mathrm{div}, \Omega)$.
\end{itemize}
The evolution of the respective errors for grid refinements under constant $\CFL=0.25$ is shown in figure \ref{pics_manufactured_energy_charge_conv}.

Figure \ref{pics_manufactured_energy_charge_conv} shows that the physical quantities of interest are well preserved by the presented schemes under $\CFL=0.25$. In particular, the energy and the total charge errors convergence with second order and the divergence of the magnetic field remains below $10^{-14}$, which is below the tolerance of $10^{-12}$ used for PCG and PBiCGStab. Finally we highlight that the Poisson splitting scheme is remarkably more accurate than the Crank-Nicolson scheme. In particular, the energy error produced by the Poisson splitting scheme is approximately one order of magnitude and in the case of the total charge, the difference is approximately two orders of magnitude.

\subsection{Performance study}

Here we study and compare the performance of the presented methods for the X-mode test case described in section \ref{sec:Xmodetest}.

\subsubsection{Operation metrics}

Since the operators have different sizes and they are applied a different number of times per time-step, the number of iterations is not enough to compare the performance of the schemes. Consequently we consider the metric
\begin{equation*}
 \mathrm{MVBP} = \text{number of Matrix-Vector Block Products (per time step)},
\end{equation*}
where a block correponds to a field:
$\arrE, \arrB$ or $\arrY$. The blocks do not have the same size, since the fields belong to different spaces, however the size difference is very small and has a negligible effect on the performance.

The $\mathrm{MVBP}$ metric does not depend on the number of coefficients within each block and therefore it cannot be used to compare different costs as the grid is refined. The number of time steps, or points per period 
($\mathrm{PPP}$), may also vary between two simulations. 
For these reasons we introduce a second metric,
namely 
\begin{equation*}
 \mathrm{LFOps} = \text{number of Local Field Operations} := \mathrm{PPP} \times \mathrm{MVBP} \times \mathrm{dim}
\end{equation*}
where $\mathrm{dim}$ is the dimension of the space $V^1_h$ (close to that of
the space $V^2_h$).
Our performance study is based on these two metrics: $\mathrm{LFOps}$ which accounts for the number of operations in one period of a given simulation, and $\mathrm{MVBP}$ per time step which indicates how well the linear systems 
are conditionned.

Next we show how to compute the average $\mathrm{MVBP}$ per time-step. Firstly, we need to count the number of matrix-vector products computed by the preconditioned conjugate gradient (PCG) and the preconditioned biconjugate gradient stabilized method (PBICGSTAB) in order to reach a fixed tolerance. This quantity is completely general, it holds for a general preconditioned linear system $P A \bs{u} = P \bs{b}$. Certainly, it only depends on the number of iterations needed to reach convergence. In Table \ref{table_performance_PCG_PBICGSTAB}, we show the general formulae.

\begin{table}[H]
\begin{center}
\begin{tabular}{ |c|c|c|c| }
\hline
\multicolumn{4}{|c|}{Number of matrix-vector products to solve $P A \arr{u} = P \arr{b}$} \\ \hline
 Solver & products with $A$ & products with $P$ & Total matrix-vector products\\ \hline
 PCG & 1 + $n$ & 1 + $n$ & 2 + 2$n$ \\ \hline
 PBiCGStab & 1 + 2$n$ & 1 + 2$n$ & 2 + 4$n$ \\ \hline
\end{tabular}
\end{center}
\caption{Number of matrix-vector products for each solver with respect to the number of iterations $n$ required to reach a fixed tolerance. Here $A$ denotes a general matrix and $P$ is the preconditioner. }
\label{table_performance_PCG_PBICGSTAB}
\end{table}
Then, we count the number of times each iterative method is called within a time-step, which is referred to as the number of \textit{solves}, as well as the number of blocks involved in each solve. These two quantities enable the computation of the $\mathrm{MVBP}$ per time-step involved in the matrix-inversion, which is shown in Table \ref{table_performance_MVBP_perniter}. In addition, we add the $\mathrm{MVBP}$ involved in the computation of the right-hand side, which can be directly infered from the schemes written previously. The total $\mathrm{MVBP}$ per time-step in terms of the average number of iterations needed for convergence of each solve is collected in the last column of Table \ref{table_performance_MVBP_perniter}.

\begin{table}[H]
\begin{center}
\begin{tabular}{ |p{18mm}|p{25mm}|p{32mm}|p{38mm}|p{39mm}| }
\hline
\multicolumn{5}{|c|}{$\mathrm{MVBP}$ per time-step} \\
 \hline
 Scheme & \#solves PCG (\#blocks) & \#solves PBiCGStab (\#blocks) & $\mathrm{MVBP}$ matrix-inversion & $\mathrm{MVBP}$ per time-step \\ \hline
 C-N & \hfil 0 & \hfil 1 (3) & \hfil 6 + 12$n$ & \hfil 15 + 12$n$ \\ \hline
 Poisson & \hfil 2 (1) & \hfil 1 (2) & 8 + 4$n_{\mathrm{Maxwell}}$ + 8$n_{\mathrm{plasma}}$ & 17 + 4$n_{\mathrm{Maxwell}}$ + 8$n_{\mathrm{plasma}}$ \\ \hline
 Hamiltonian & \hfil 2 (1) & \hfil 1 (2) & \hfil 8 + 4$n_{\arrE}$ + 8$n_{\arrB, \arrY}$ & \hfil 18 + 4$n_{\arrE}$ + 8$n_{\arrB, \arrY}$ \\ \hline
\end{tabular}
\end{center}
\caption{Number of $\mathrm{MVBP}$ per time-step for each scheme. Here $n, n_{\mathrm{Maxwell}}, n_{\mathrm{plasma}}, n_{\arrE}$ and $n_{\arrB, \arrY}$ denote the number of iterations taken to invert the matrices of $\Phi^{\mathrm{CN}}, \Phi^P_{\mathrm{Maxwell}}, \Phi^P_{\mathrm{plasma}}, \Phi^H_{\arrE}$ and $\Phi^H_{\arrB, \arrY}$ respectively.
}
\label{table_performance_MVBP_perniter}
\end{table}
All the matrices are stored with the \textit{StencilMatrix} format of Psydac, which only stores the coefficients of the corresponding stencil and thus it is very memory-efficient.

\subsubsection{Benchmarking the X-mode test case: number of MVBP per time step}

In order to compute the $\mathrm{MVBP}$ as in Table \ref{table_performance_MVBP_perniter}, it is necessary to collect the average number of iterations per time-step over three periods.
Tables \ref{table_fixedCFL} and \ref{table_fixedPPW} gather the average number of iterations required for each matrix inversion and the average $\mathrm{MVBP}$ for each method. In particular, table \ref{table_fixedCFL} exhibits the performance under a constant $\CFL$, while table \ref{table_fixedPPW} does for a varying $\CFL$ with constant space resolution $\mathrm{PPW}=10$.

\begin{table}[H]
\begin{center}
\begin{tabular}{ |c|c|c|c|c|c|c|c| }
\cline{3-8}
 \multicolumn{2}{c|}{} & \multicolumn{3}{c|}{Average number of iterations} & \multicolumn{3}{c|}{Average $\mathrm{MVBP}$} \\ \cline{3-8}
 \hline
 $\mathrm{PPW}$ & $\mathrm{PPP}$ & C-N & Poisson & Hamiltonian & C-N & Poisson & Hamiltonian \\ \hline
 10 & 40 & $n=11.8$ & $n_{\mathrm{Maxwell}}= 8.7, n_{\mathrm{plasma}}= 4$ & $n_{\arr{E}}=2, n_{\arr{B}, \arr{Y}}=4$ & 147.6 & 74.8 & 48 \\ \hline
 20 & 80 & $n=11.1$ & $n_{\mathrm{Maxwell}}= 7.8, n_{\mathrm{plasma}}= 3.4$ & $n_{\arr{E}}=2, n_{\arr{B}, \arr{Y}}=4$ & 147.6 & 74.8 & 48 \\ \hline
 40 & 160 & $n=10.7$ & $n_{\mathrm{Maxwell}}= 7.6, n_{\mathrm{plasma}}= 3$ & $n_{\arr{E}}=2, n_{\arr{B}, \arr{Y}}=4$ & 134.4 & 62.4 & 48\\ \hline
 80 & 320 & $n=9.95$ & $n_{\mathrm{Maxwell}}= 7.1, n_{\mathrm{plasma}}= 3$ & $n_{\arr{E}}=2, n_{\arr{B}, \arr{Y}}=4$ & 125.4 & 60.4 & 48 \\ \hline
\end{tabular}
\end{center}
\caption{Number of iterations and $\mathrm{MVBP}$ per time-step with fixed $\CFL=0.25$, averaged over three periods.
}
\label{table_fixedCFL}
\end{table}

\begin{figure}[H]
         \centering %[trim=left bottom right top, clip]
         \begin{subfigure}[t]{0.32\textwidth}
         \includegraphics[width=\textwidth, trim=0.8cm 0.8cm 0.5cm 2.4cm, clip]{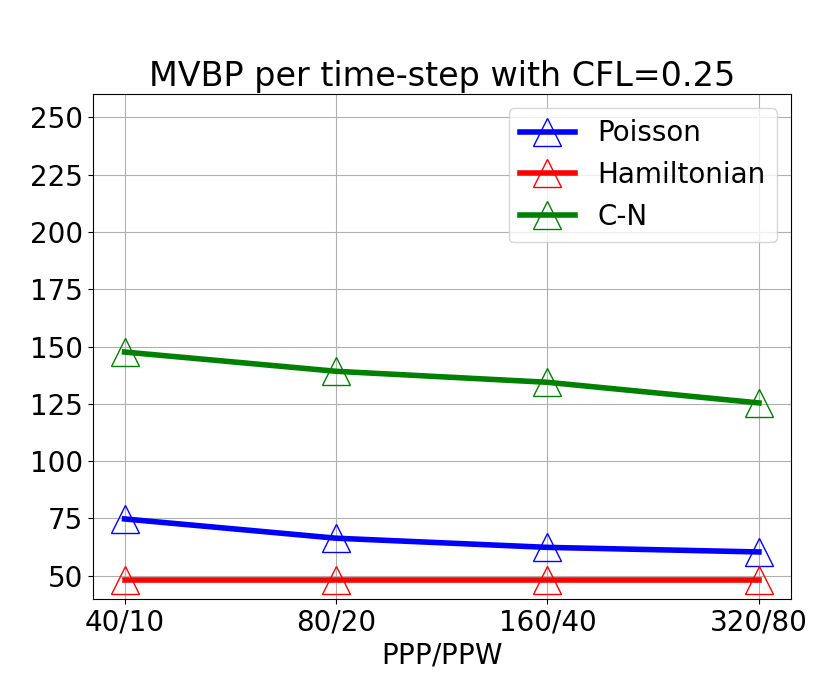}
         \caption{Average $\mathrm{MVBP}$ per time-step as grid is refined with $\CFL=0.25$.}
         \label{pics_manufactured_MVBP:a}
         \end{subfigure}
         \hfill
         \begin{subfigure}[t]{0.32\textwidth}
         \includegraphics[width=\textwidth, trim=0.8cm 0.8cm 0.5cm 2.4cm, clip]{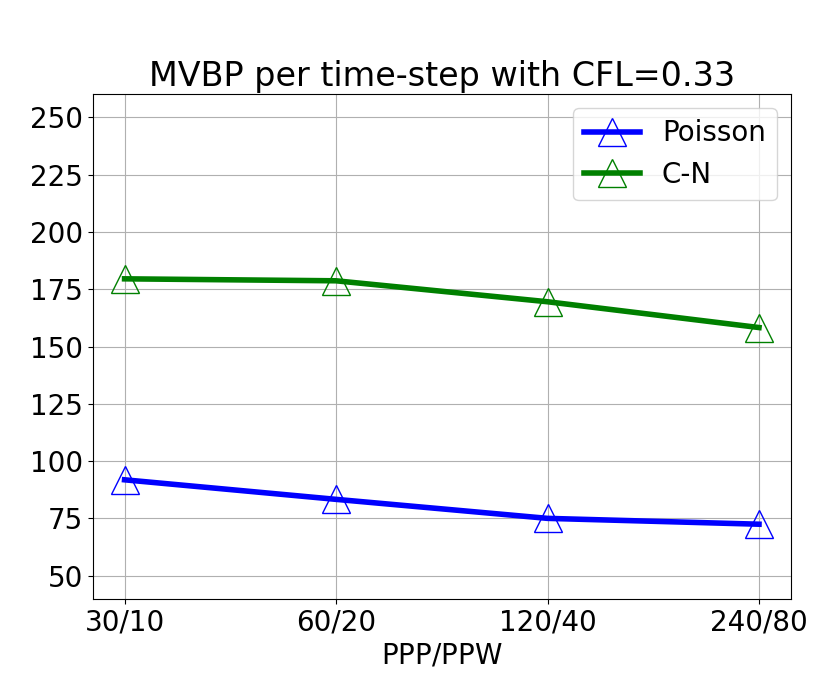}
         \caption{Average $\mathrm{MVBP}$ per time-step as grid is refined with $\CFL=0.33$.}
         \label{pics_manufactured_MVBP:b}
         \end{subfigure}
         \hfill
         \begin{subfigure}[t]{0.32\textwidth}
         \includegraphics[width=\textwidth, trim=0.8cm 0.8cm 0.5cm 2.4cm, clip]{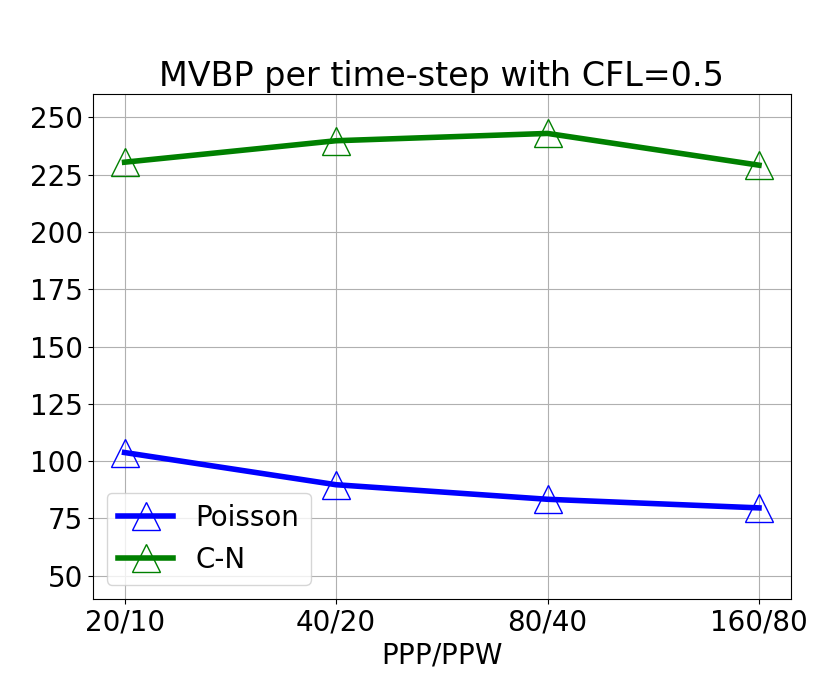}
         \caption{Average $\mathrm{MVBP}$ per time-step as grid is refined with $\CFL=0.5$.}
         \label{pics_manufactured_MVBP:c}
         \end{subfigure}
          \caption{Evolution of the average $\mathrm{MVBP}$ per time-step as the grid is refined under a constant $\CFL$.
          }
         \label{pics_manufactured_MVBP}
 \end{figure}

\begin{table}[H]
\begin{center}
\begin{tabular}{ |c|c|c|c|c|c|c|c| }
\cline{3-8}
 \multicolumn{2}{c|}{} & \multicolumn{3}{c|}{Average number of iterations} & \multicolumn{3}{c|}{Average $\mathrm{MVBP}$} \\ \cline{3-8}
\hline
 $\CFL$ &$\mathrm{PPP}$ & C-N & Poisson & Hamiltonian & C-N & Poisson & Hamiltonian \\ \hline
% 44 (CFL=0.23) & &  & \\ \hline
 0.25 & 40 & $n=11.8$ & $n_{\mathrm{Maxwell}}= 8.7, n_{\mathrm{plasma}}= 4$ & $n_{\arr{E}}=2, n_{\arr{B}, \arr{Y}}=4$ & 147.6 & 74.8 & 48 \\ \hline
 0.33 & 30 & $n=13.7$ & $n_{\mathrm{Maxwell}}= 9.6, n_{\mathrm{plasma}}= 4.5$ & n.c. & 170.4 & 82.4 & n.c. \\ \hline
 0.5 & 20 & $n=18$ & $n_{\mathrm{Maxwell}}= 10.9, n_{\mathrm{plasma}}= 5.4$ & n.c. & 222 & 94.8 & n.c. \\ \hline
 1 & 10 & $n=34.2$ & $n_{\mathrm{Maxwell}}= 13.9, n_{\mathrm{plasma}}= 6.3$ & n.c. & 416.4 & 114 & n.c. \\ \hline
\end{tabular}
\end{center}
\caption{Number of iterations and $\mathrm{MVBP}$ per time-step with fixed $\mathrm{PPW}=10$, averaged over three periods. For the cases where the solution did not converge \ref{pics_manufactured_solerror_fixedncells},  we use the abbreviation n.c. denoting ``not converged''. }
\label{table_fixedPPW}
\end{table}
Figure \ref{pics_manufactured_MVBP} shows the evolution of the average $\mathrm{MVBP}$ per time-step as the grid is refined with a fixed $\CFL$.
Figure \ref{pics_manufactured_MVBP} shows that the time-splitting schemes produce less average $\mathrm{MVBP}$ per time-step than the Crank-Nicolson scheme. This is expected, as the Crank-Nicolson scheme solves the full system (three blocks) in every time-step, while the splitting schemes involve solving smaller subsystems. Besides, we observe that the average $\mathrm{MVBP}$ per time-step produced by the Poisson splitting scheme is less sensitive to the value of $\CFL$. In other words, as $\CFL$ increases the conditioning of the problem deteriorates faster for the Crank-Nicolson scheme than for the Poisson splitting scheme.

Another observation is that in general the average $\mathrm{MVBP}$ per time-step decreases as the grid is refined with constant $\CFL$. This shows the effect of the block Kronecker mass preconditioner, which is more effective as the time-step reduces. More concretely, the diagonal of the matrices to invert is dominated by a block mass matrix and the perturbation is at least of the order of $\Delta t$. In other words, in all the schemes, the matrices to invert in the limit as $\Delta t \rightarrow 0$ become the preconditioner.

In \ref{pics_manufactured_MVBP:a} we observe that for $\CFL=0.25$ the Hamiltonian splitting scheme produces the lowest average $\mathrm{MVBP}$ per time-step . In particular, in table \ref{table_fixedCFL} we see that the number of iterations required by the Hamiltonian splitting scheme remains constant and small as the grid is refined. On the contrary, while the number of iterations required by $\Phi^P_\mathrm{plasma}$ remains more or less constant, the number of iterations of $\Phi^P_\mathrm{Maxwell}$ increases and is overall larger.
Indeed, operator $\Phi^H_{\arrE}$ \eqref{integration_Hamiltonian_splitting_operatorA} involves inverting a mass matrix, which coincides exactly with our preconditioner, and operator $\Phi^H_{\arrB, \arrY}$ involves inverting a matrix, the diagonal of which is a block mass matrix only perturbed by boundary terms of the order of $\Delta t$. On the contrary, the matrix to invert in \eqref{integration_Poisson_splitting_operatorA:avg_E} contains the term $(\Delta t^2/4) \matC^T \matM_2 \matC$, which introduces a perturbation from the preconditioner everywhere.

\subsubsection{Benchmarking the X-mode test case: total cost}
\begin{figure}[H]
         \begin{subfigure}[t]{0.32\textwidth}
         \centering %[trim=left bottom right top, clip]
         \includegraphics[width=\textwidth, trim=0.8cm 0.8cm 0.5cm 2.4cm, clip]{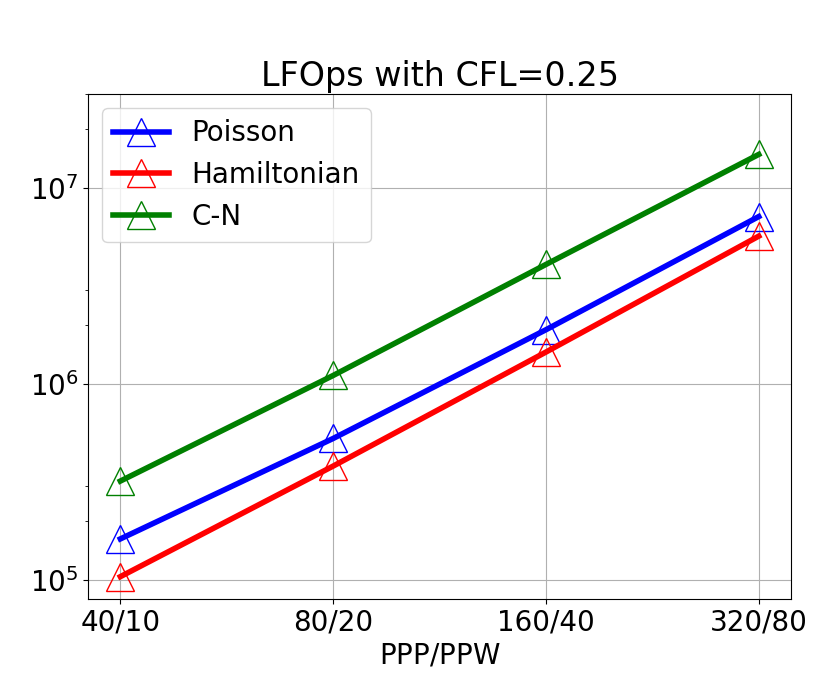}
         \caption{$\mathrm{LFOps}$ as grid is refined with $\CFL=0.25$.}
         \label{pics_manufactured_LFOps:a}
         \end{subfigure}
         \hfill
         \begin{subfigure}[t]{0.32\textwidth}
         \includegraphics[width=\textwidth, trim=0.8cm 0.8cm 0.5cm 2.4cm, clip]{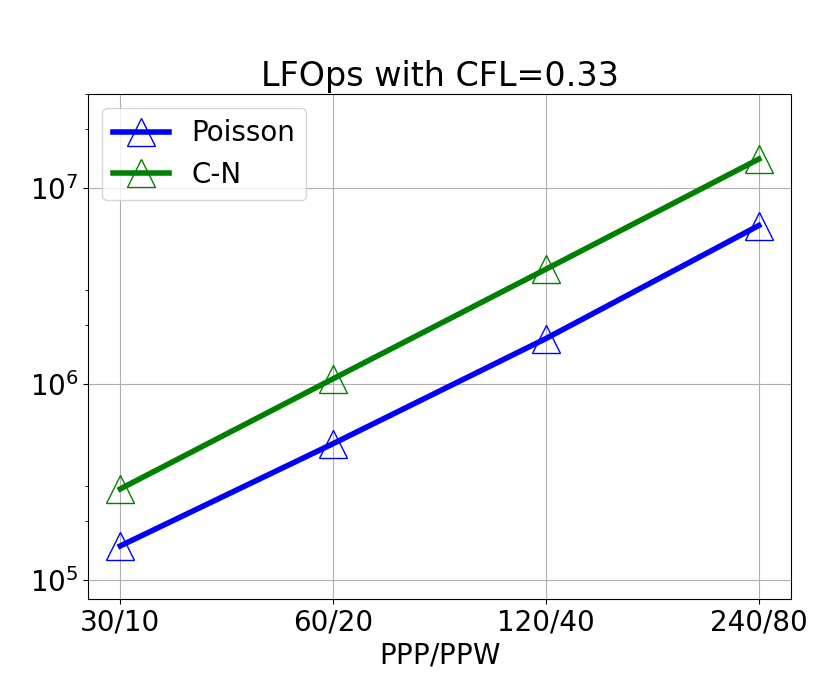}
         \caption{$\mathrm{LFOps}$ as grid is refined with $\CFL=0.33$.}
         \label{pics_manufactured_LFOps:b}
         \end{subfigure}
         \hfill
         \begin{subfigure}[t]{0.32\textwidth}
         \includegraphics[width=\textwidth, trim=0.8cm 0.8cm 0.5cm 2.4cm, clip]{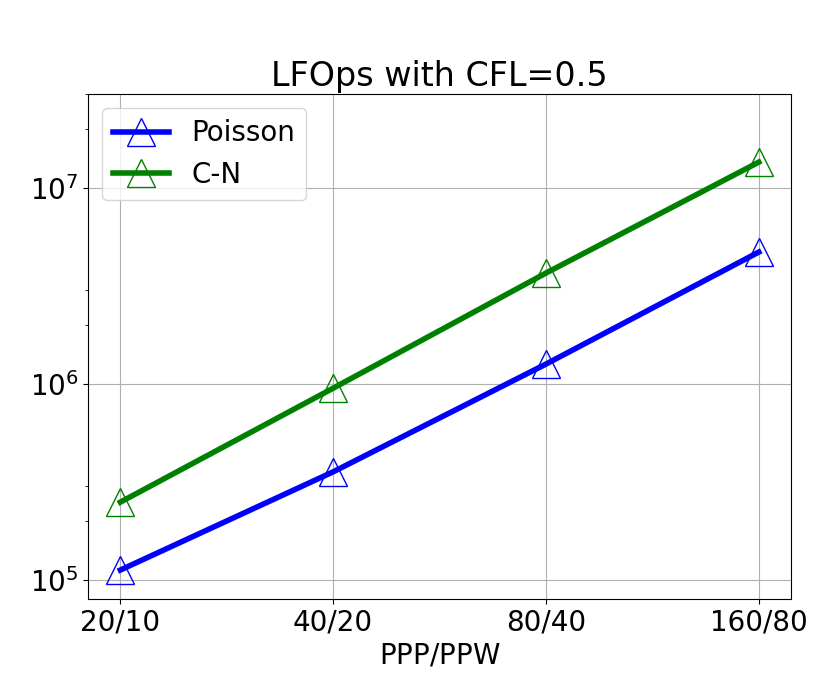}
         \caption{$\mathrm{LFOps}$ as grid is refined with $\CFL=0.5$.}
         \label{pics_manufactured_LFOps:c}
         \end{subfigure}
         \\[0.3em]
         \begin{subfigure}[t]{0.32\textwidth}
         \includegraphics[width=\textwidth, trim=0.8cm 0.8cm 0.5cm 2.4cm, clip]{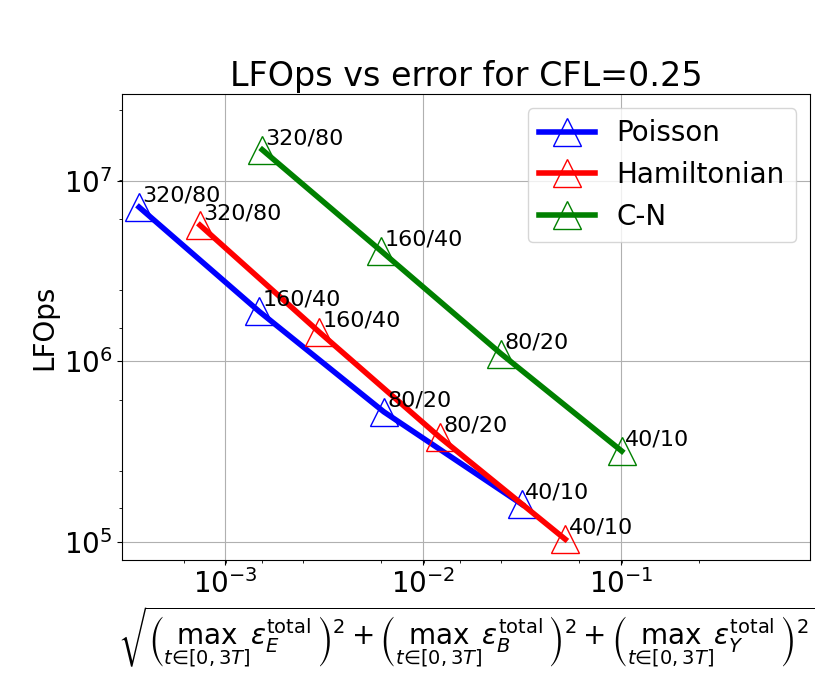}
         \caption{$\mathrm{LFOps}$ and error as grid is refined with $\CFL=0.25$.}
         \label{pics_manufactured_LFOps:d}
         \end{subfigure}
         \hfill
         \begin{subfigure}[t]{0.32\textwidth}
         \includegraphics[width=\textwidth, trim=0.8cm 0.8cm 0.5cm 2.4cm, clip]{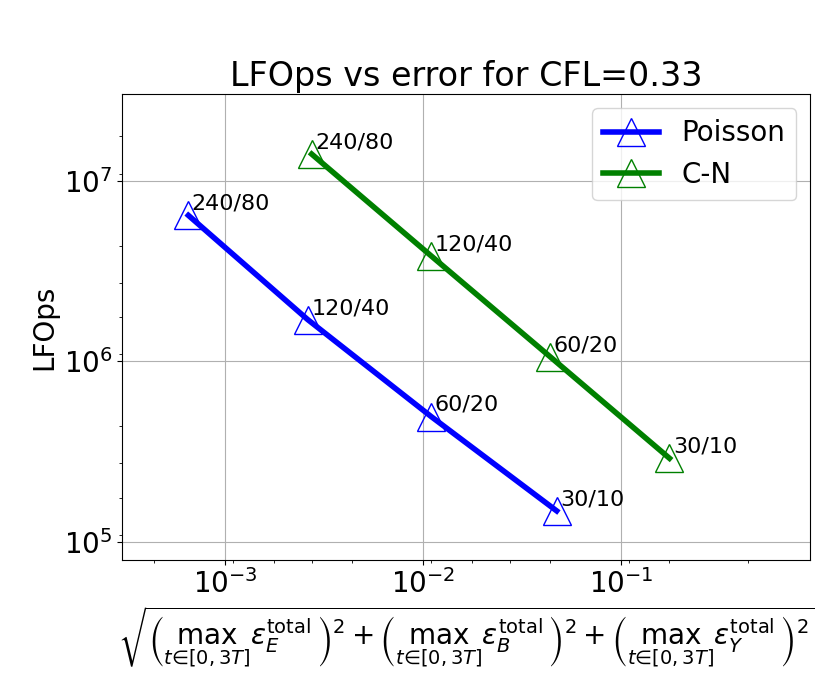}
         \caption{$\mathrm{LFOps}$ and error as grid is refined with $\CFL=0.33$.}
         \label{pics_manufactured_LFOps:e}
         \end{subfigure}
         \hfill
         \begin{subfigure}[t]{0.32\textwidth}
         \includegraphics[width=\textwidth, trim=0.8cm 0.8cm 0.5cm 2.4cm, clip]{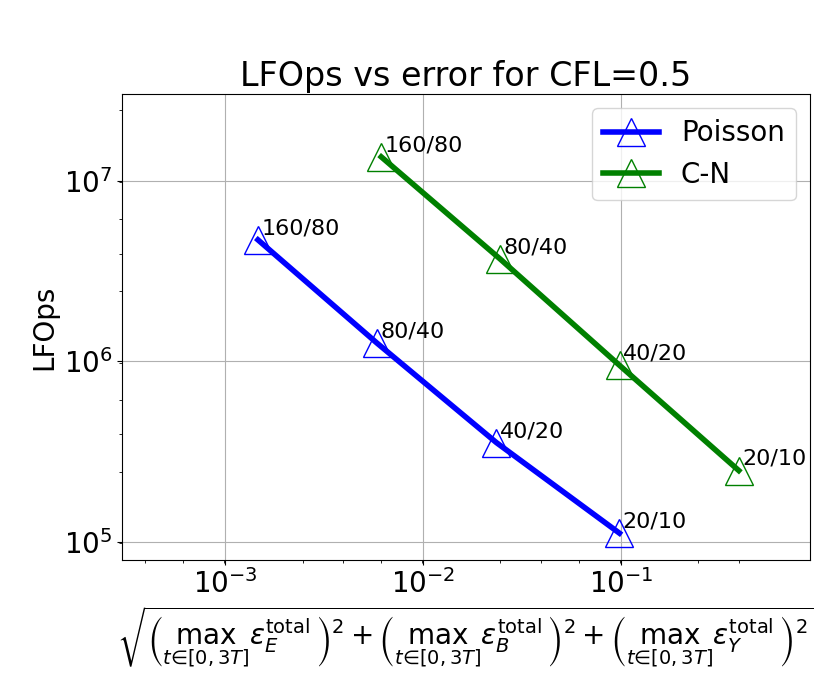}
         \caption{$\mathrm{LFOps}$ and error as grid is refined with $\CFL=0.5$.}
         \label{pics_manufactured_LFOps:f}
         \end{subfigure}
         \caption{Top: Evolution of $\mathrm{LFOps}$ as the grid is refined with constant $\CFL$. Bottom: Evolution of $\mathrm{LFOps}$ with respect to the total error as the grid is refined with fixed $\CFL$. }
         \label{pics_manufactured_LFOps}
 \end{figure}

Figures \ref{pics_manufactured_LFOps:a}, \ref{pics_manufactured_LFOps:b} and \ref{pics_manufactured_LFOps:c} show the evolution of the cost of each method in terms of $\mathrm{LFOps}$ as the grid is refined with constant $\CFL$. The plots  \ref{pics_manufactured_LFOps:d}, \ref{pics_manufactured_LFOps:e} and \ref{pics_manufactured_LFOps:f} show the evolution of the cost with respect to the accuracy.
From \ref{pics_manufactured_LFOps:a}, \ref{pics_manufactured_LFOps:b} and \ref{pics_manufactured_LFOps:c} we observe that the cost is dominated by the increase in $\mathrm{PPW}$ and $\mathrm{PPP}$. More concretely, it behaves as $\mathrm{LFOps} \sim (1/\Delta x)^k$ for some order $k > 0$ (corresponding to the slope), which is similar for all the schemes. Moreover, for every $\mathrm{PPW}$ the cost decreases as $\CFL$ increases (smaller $\mathrm{PPP}$).
The second row of plots, i.e. \ref{pics_manufactured_LFOps:d}, \ref{pics_manufactured_LFOps:e} and \ref{pics_manufactured_LFOps:f}, shows clearly that the cost increases as the error decreases. Besides, to achieve a certain error, the cost increases as $\CFL$ increases and the Poisson splitting method has the lowest cost in almost all cases. The figure also shows that to reach a certain error tolerance, the Crank-Nicolson method is remarkably more expensive than the methods based on time-splitting.

We conclude that the Poisson splitting method is the most performant. Furthermore, it can be used with larger $\CFL$ values without a significant sacrifice in performance.

\section{Application to 2D turbulent plasma densities}
In the performance study we have seen that the Poisson splitting scheme \eqref{scheme_Poisson} has the best cost-accuracy relation. In this last section we use it to run two-dimensional experiments. The domain is fixed to $\Omega = [0, 24\pi] \times [0, 24\pi] \times [0, 2\pi]$, which corresponds to $12$ wavelengths in the $x$ and $y$ directions. The solution is forced to be constant in $z$ by setting periodic boundary conditions and imposing a single cell along this direction. For the $x,y$ directions,
we set the number of points per wavelength in the $x$ and $y$ directions as $\mathrm{PPW}_x = \mathrm{PPW}_y =7$ and impose Silver-Müller boundary conditions. The incoming wave is the Gaussian beam defined in \eqref{gaussian_beam_3d} with $\Gamma_s^\inc=\{x=0\}\cap \partial \Omega$.
The beam is launched normal to the boundary and centered in the $y$-axis. Furthermore, we consider an envelope of the form
\begin{equation}
\chi(t) = 2/\pi \arctan (t/ (20 \Delta t)),
\label{gaussian_beam_envelope}
\end{equation}
which is used to avoid introducing high frequencies in the spectrum during the first time-steps.
As in section \ref{section_1Dtests}, the background magnetic field is set constant $\nb_0 = (0,0,1)$ with $\nomegac=0.5$ and the electron-collision coefficient is $\hat{\nu}_e=0$. The initial condition is set to zero and $\mathrm{PPP}=32$, which corresponds to $\CFL \simeq 0.22$.

For the experiments we manufacture two density patterns. The first pattern is made of the superposition of Gaussian functions centered at different points. The superposition is multiplied by an envelope of the form $\xi(\bx)=x$, to ensure that the incoming wave is launched from vacuum. Besides, the pattern is scaled so that it contains the cutoff contour but not the resonances. The second pattern consists on a single blob of high density, which is placed slightly below the beam axis. In total we run three experiments. The first involves an O-mode polarization and the complicated density pattern. The second consists of an X-mode polarization with the single blob density. Finally, the third one involves an X-mode polarization with the complicated density pattern. Our code allows more general polarizations, nevertheless here we focus only on pure X-mode and O-mode configurations for simplicity of the analysis.

In addition, the electric field of the solution is compared with the time-harmonic field $\bE^\mathrm{th}_h$, which corresponds to the electric field of the computed time-harmonic solution with spline coefficients $\arrE^\mathrm{th} = \Re \{ \hat{\arrE} e^{-it}\}$. The array $\hat{\arrE}$ is the solution of the discrete frequency-domain problem \eqref{semidiscrete_freqdomain}, which we solve using the parallel direct solver MUMPS \footnote{\href{https://mumps-solver.org/index.php}{https://mumps-solver.org/index.php}}, which is accessible from Psydac through PETSc \cite{petsc}, \cite{petsc4py}. To study the convergence to the time-harmonic regime, for simplicity we track in time the quantity
\begin{equation}
\bs{R}(\bx, t) = \cfrac{\bE_h(\bx, t) - \bE^\mathrm{th}_h(\bx, t) }{\max_t \{ \| \bE^\mathrm{th}_h \|_{L^2(\Omega)}, \| \bE_h \|_{L^2(\Omega)} \}}.
 \label{rel_diff_conv_timeharmonic}
\end{equation}
We denote by $\bE^\mathrm{freq}$ the complex solution of the frequency-domain problem \eqref{weak_form_freq_pb}, such that $\bs{\sigma}^1(\bE^\mathrm{freq}) = \hat{\arrE}$.

 To our present knowledge there is not a proved limited amplitude principle that ensures the convergence $\lim_{t\rightarrow \infty} \|\bs{R}\|_{L^2(\Omega)} = 0$ in the case of a general cold plasma. Physically, the wave goes through a transient phase until it reaches an equilibrium, where the behaviour is similar to $\bE^\mathrm{th}_h$. Therefore, we use $\bs{R}$ as an indicator of the qualitative behaviour of the solution after long simulation times.

\begin{remark}[Relation to the long-time stability]
 The long-time stability result shown for the Poisson splitting scheme in \eqref{long_time_stability_coeffs_Poisson} does not imply that $\bs{R}$ decreases in time.
 The analysis shows that under certain assumptions, we can construct a time-harmonic solution $\bs{U}_h^\mathrm{th}$ such that $ \| \bs{U}_h(t^n) - \bs{U}_h^\mathrm{th}(t^n)\|_{L^2(\Omega)}$ cannot increase as $n$ increases. Here $\bs{U}_h^\mathrm{th}$ is a solution of the time-domain problem \eqref{weak_form} and not the frequency-domain problem \eqref{weak_form_freq_pb}. In fact, we note that the coefficients of $\bs{U}_h^\mathrm{th}$ given in \eqref{coeffs_timeharmonic_Poisson} depend on $\Delta t$. In other words, the electric field of $\bs{U}_h^\mathrm{th}$ is an approximation of $\bs{E}_h^\mathrm{th}$ using the time-evolution problem and therefore it involves the errors related to the time discretization. In contrast, the frequency-domain problem \eqref{semidiscrete_freqdomain} does not involve the time discretization and therefore it provides a more reliable reference solution.

 In fact, there are certain reasons to suspect that in practice $\lim_{t\rightarrow \infty} \|\bs{R}\|_{L^2(\Omega)} \neq 0$. Firstly, the source spectrum in the time-domain and frequency-domain problem is not the same. The time-domain problem has an initialization, which introduces high frequencies in the spectrum. Besides, the envelope for the Gaussian beam given in \eqref{gaussian_beam_envelope} introduces certain artificial frequencies in the first time-steps. Unless the corresponding modes are part of the transient regime and propagate outwards through the boundary $\Gamma_A$, they inevitably contribute to $\bs{R}$.
 A secondary factor that may influence the convergence of $\bs{R}$ is that, contrary to the numerical solution, $\bs{E}_h^\mathrm{th}$ is not affected by the time-splitting error.
 \label{remark_conv_timeharmonic}

\end{remark}

\subsection{O-mode with high density fluctuations}
\begin{figure}[H]
        \centering %[trim=left bottom right top, clip]
        \includegraphics[height=0.25\textheight, trim=0 2.5cm 0cm 5cm, clip]{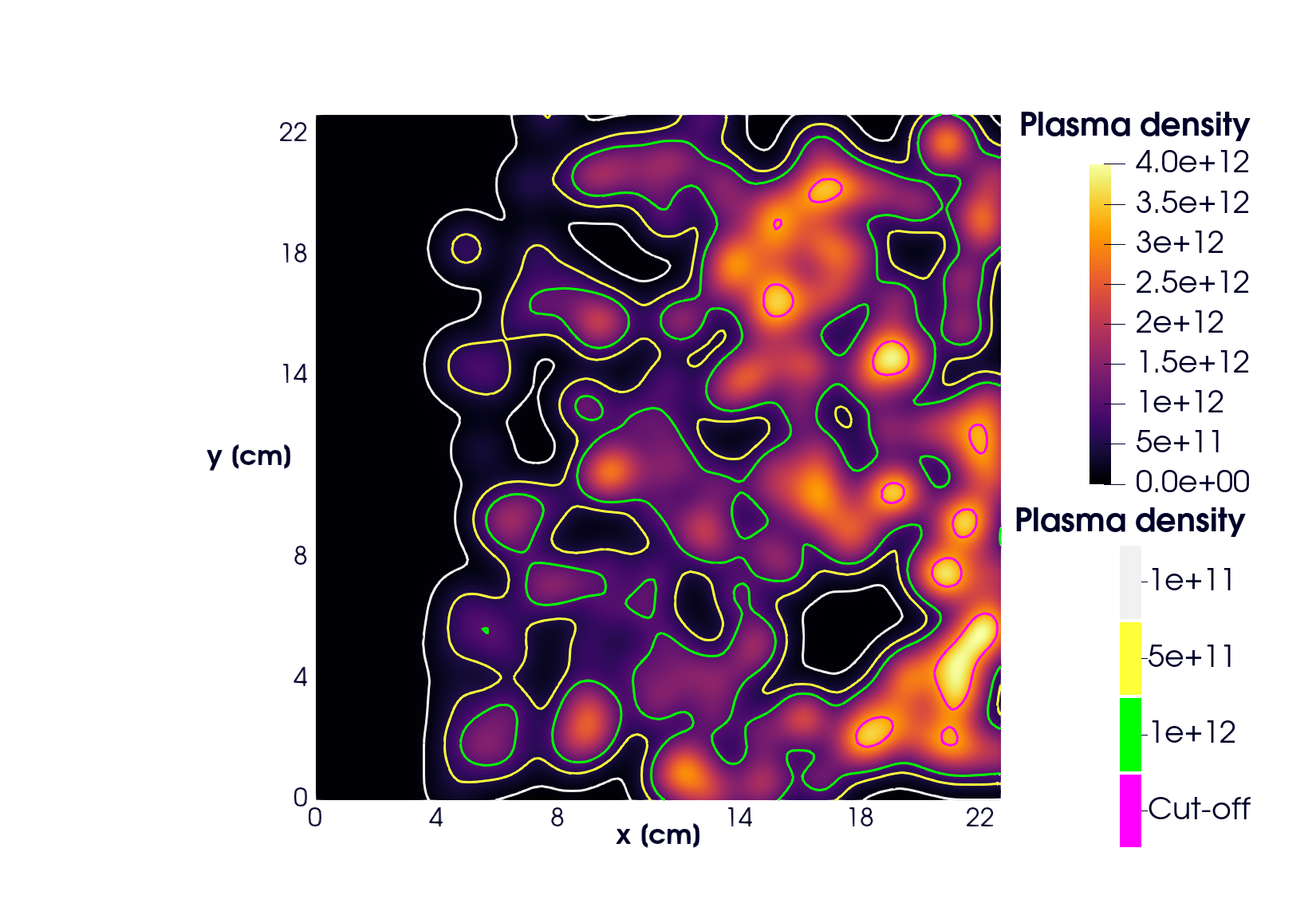}
        \caption{Two-dimensional electron plasma density with high fluctuations. The O-mode cutoff contour (magenta contour) is given by the relation $\nomegap^2(\bx)=1$ \cite{chen2013introduction}, which corresponds approximately to $n_e(\bx) =3 \times 10^{12} $ electrons per cubic centimeter.}
        \label{pics_example_density}
\end{figure}

\begin{figure}[H]
         \centering %[trim=left bottom right top, clip]
          \includegraphics[width=0.49\textwidth, trim=9.5cm 2.5cm 4cm 4.5cm clip]{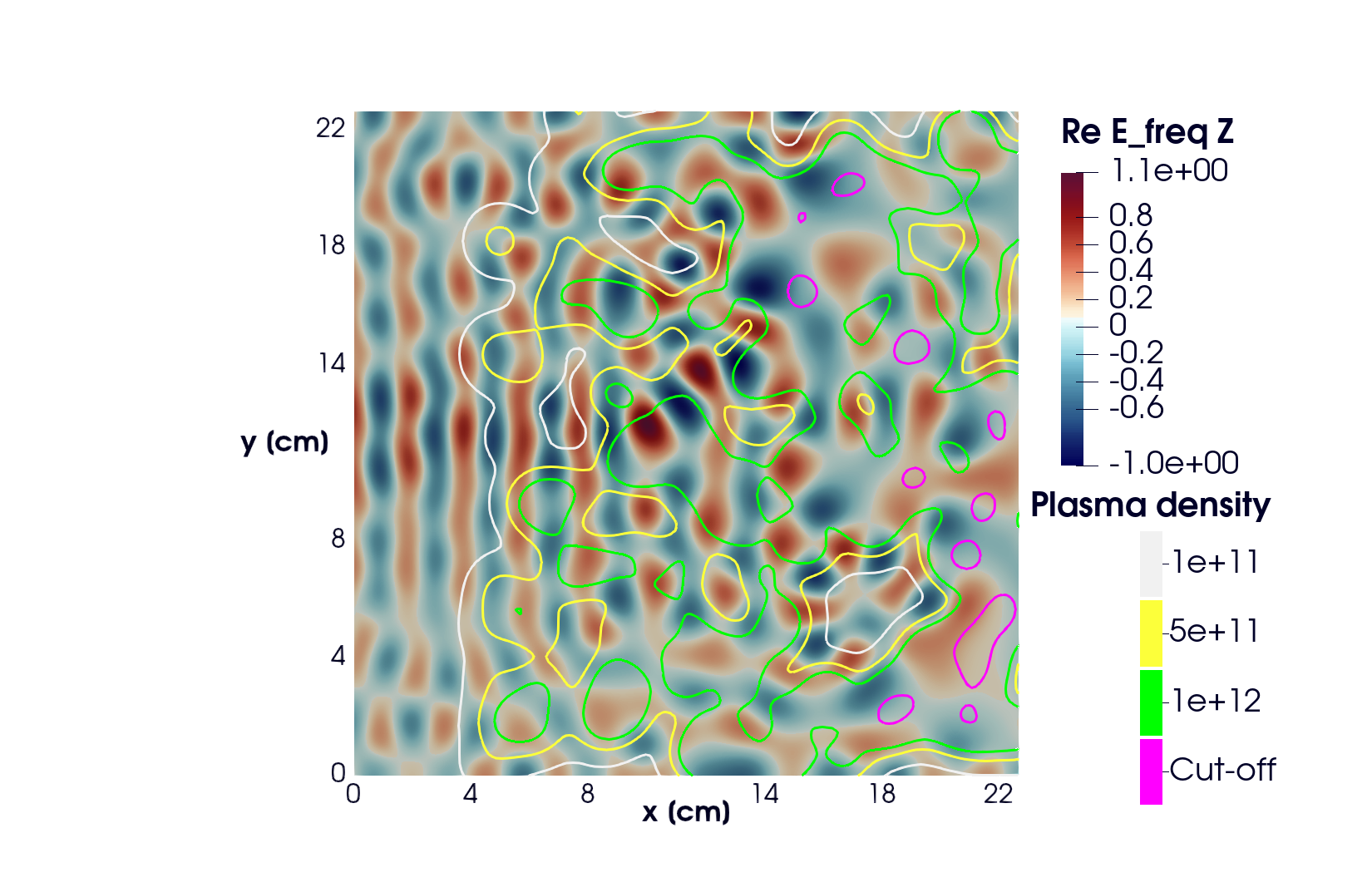}
          \includegraphics[width=0.49\textwidth, trim=9.5cm 2.5cm 4cm 4.5cm, clip]{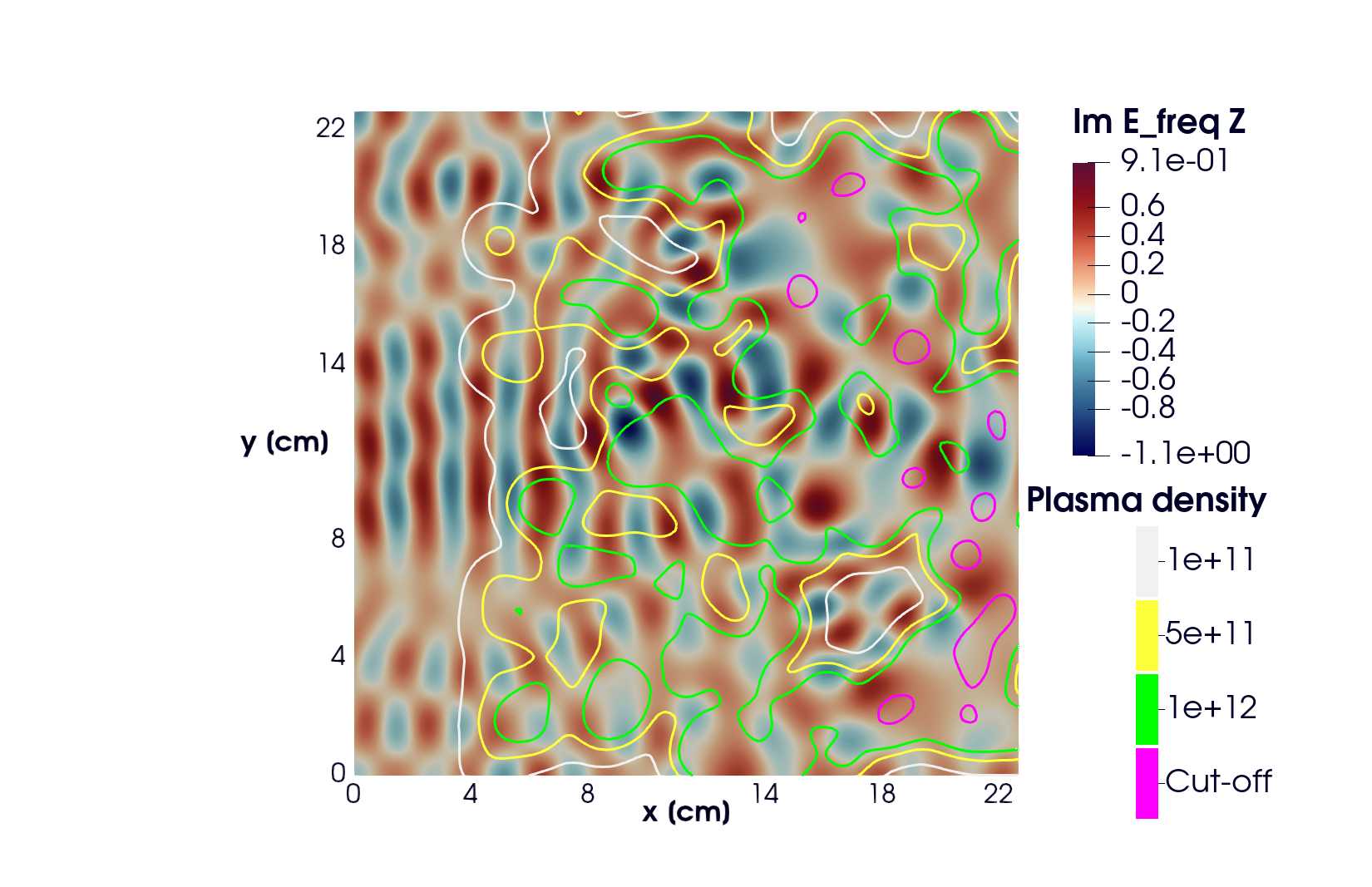}
          \includegraphics[width=0.24\textwidth, trim=14.8cm 6cm 14.8cm 4.5cm, clip]{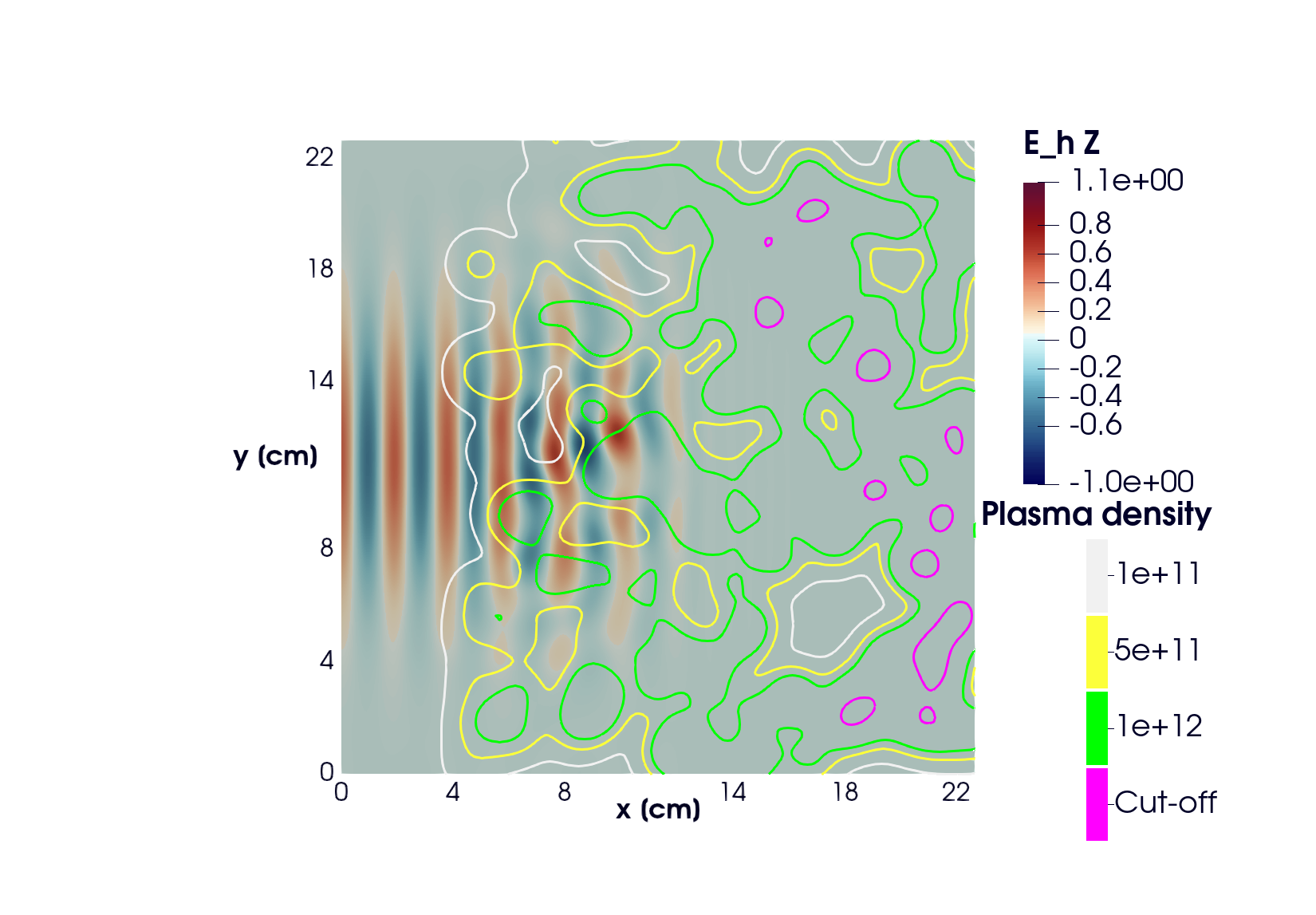}
          \includegraphics[width=0.24\textwidth, trim=14.8cm 6cm 14.8cm 4.5cm, clip]{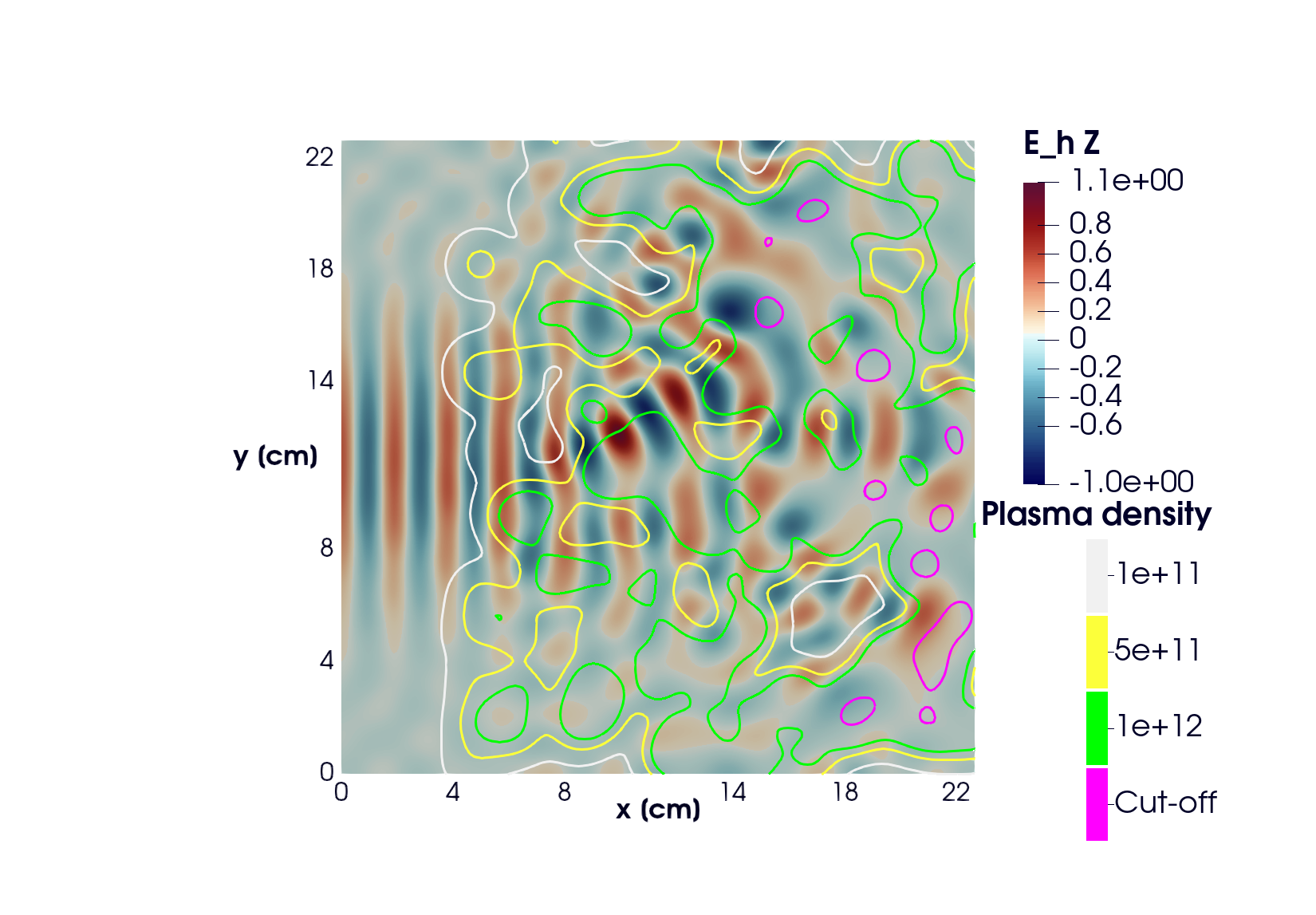}
          \includegraphics[width=0.24\textwidth, trim=14.8cm 6cm 14.8cm 4.5cm, clip]{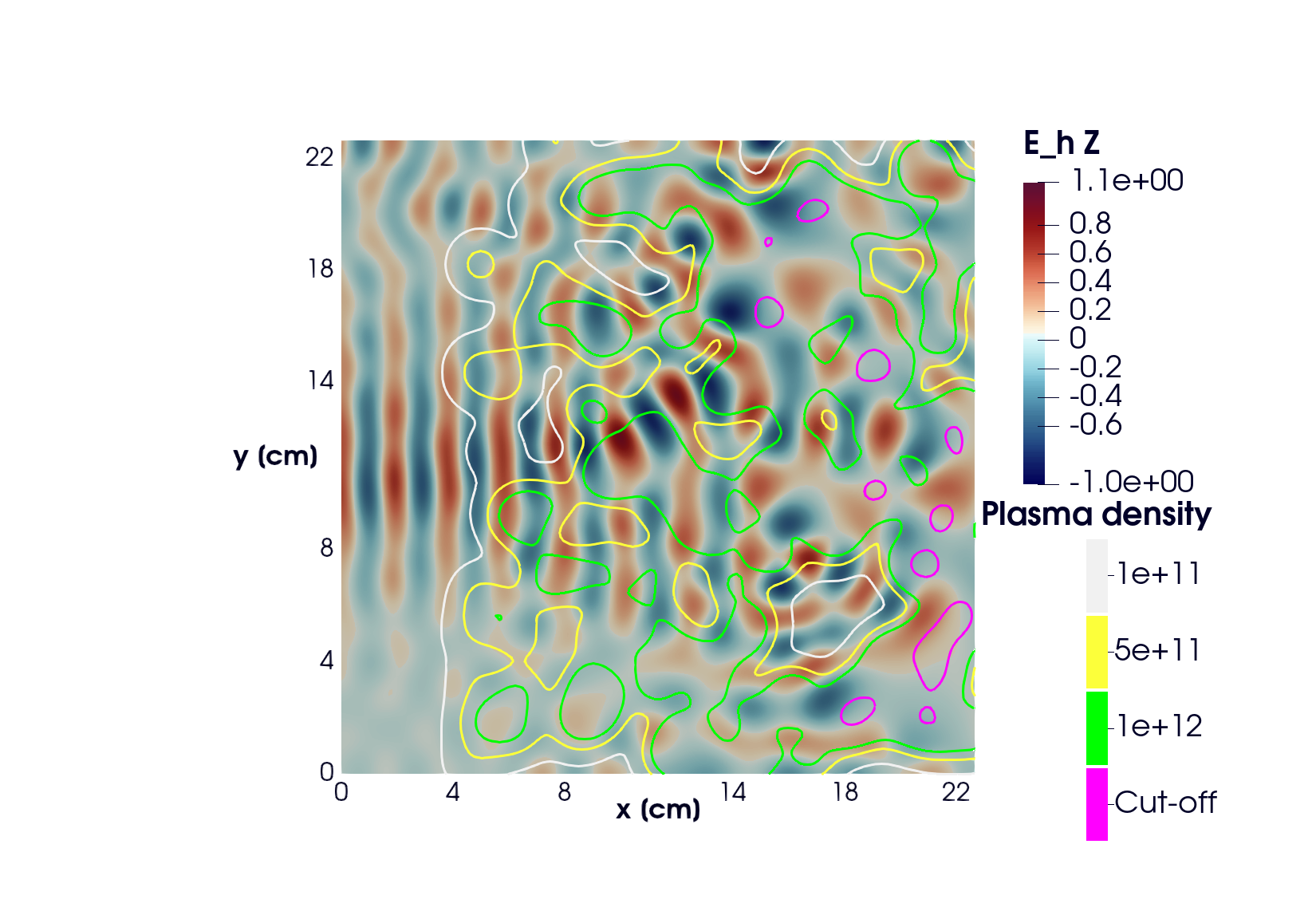}
          \includegraphics[width=0.24\textwidth, trim=14.8cm 6cm 14.8cm 4.5cm, clip]{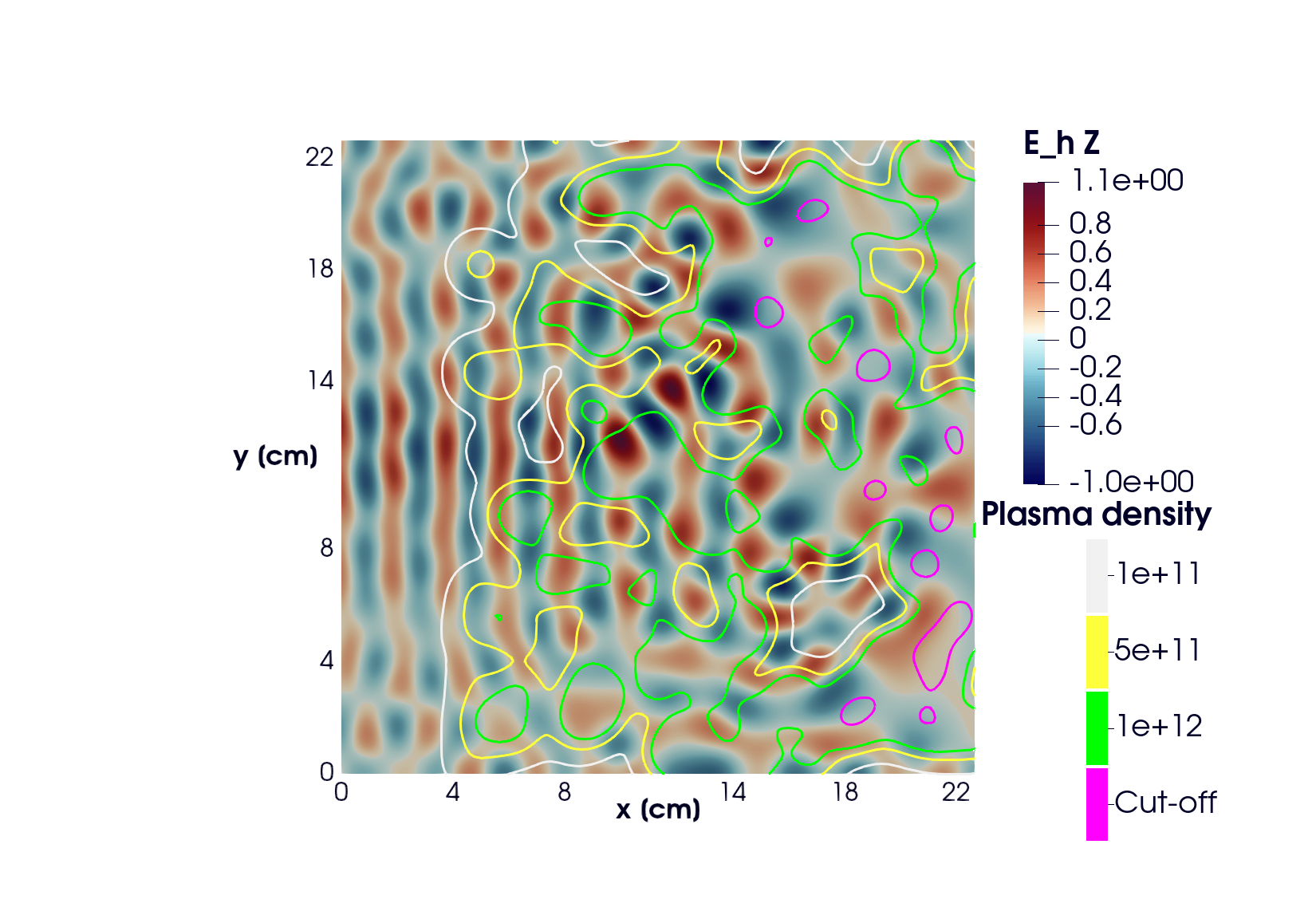}
         \caption{Top: real and imaginary part of $( \bE^{\mathrm{freq}} )_z$. Bottom: evolution of $( \bE_h )_z$ at $t = 7T, 16T, 23T,50T$ (increasing from left to right). The color map is fixed to the amplitude of $( \Re \{\bE^{\mathrm{freq}}\} )_z$, which is $[-1, 1.1]$.}
         \label{pics_example_Omode_solution}
 \end{figure}

\begin{figure}[H]
          \centering %[trim=left bottom right top, clip]
          \includegraphics[width=0.4\textwidth, trim=0.8cm 1cm 0.8cm 4.6cm, clip]{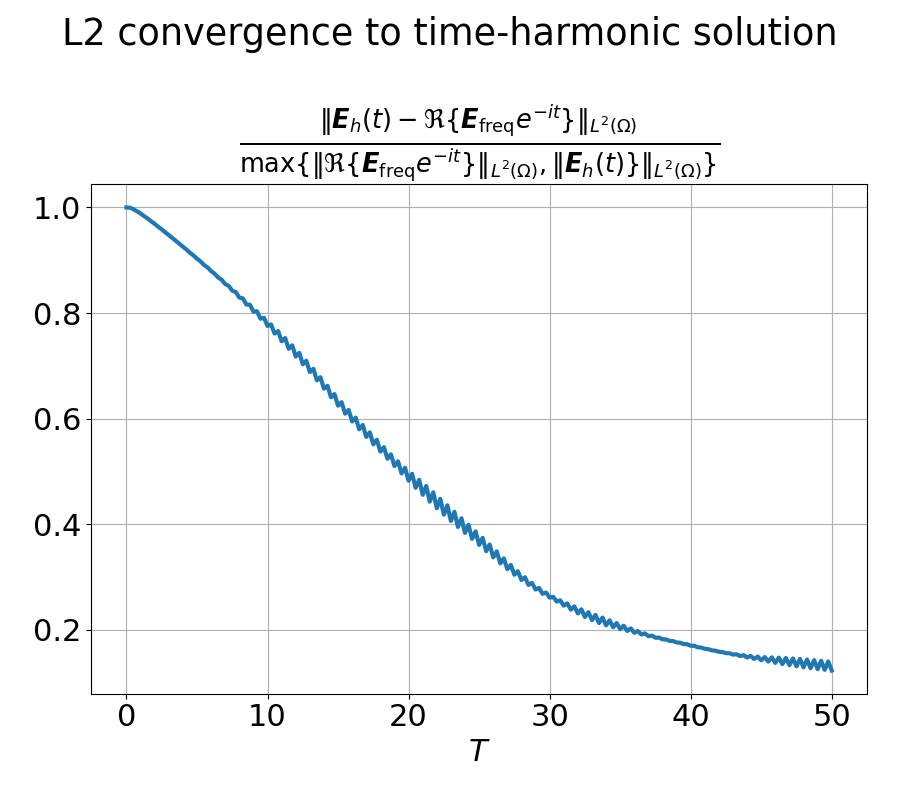}
          \includegraphics[width=0.59\textwidth, trim=9.5cm 2.5cm 3.5cm 4.5cm, clip]{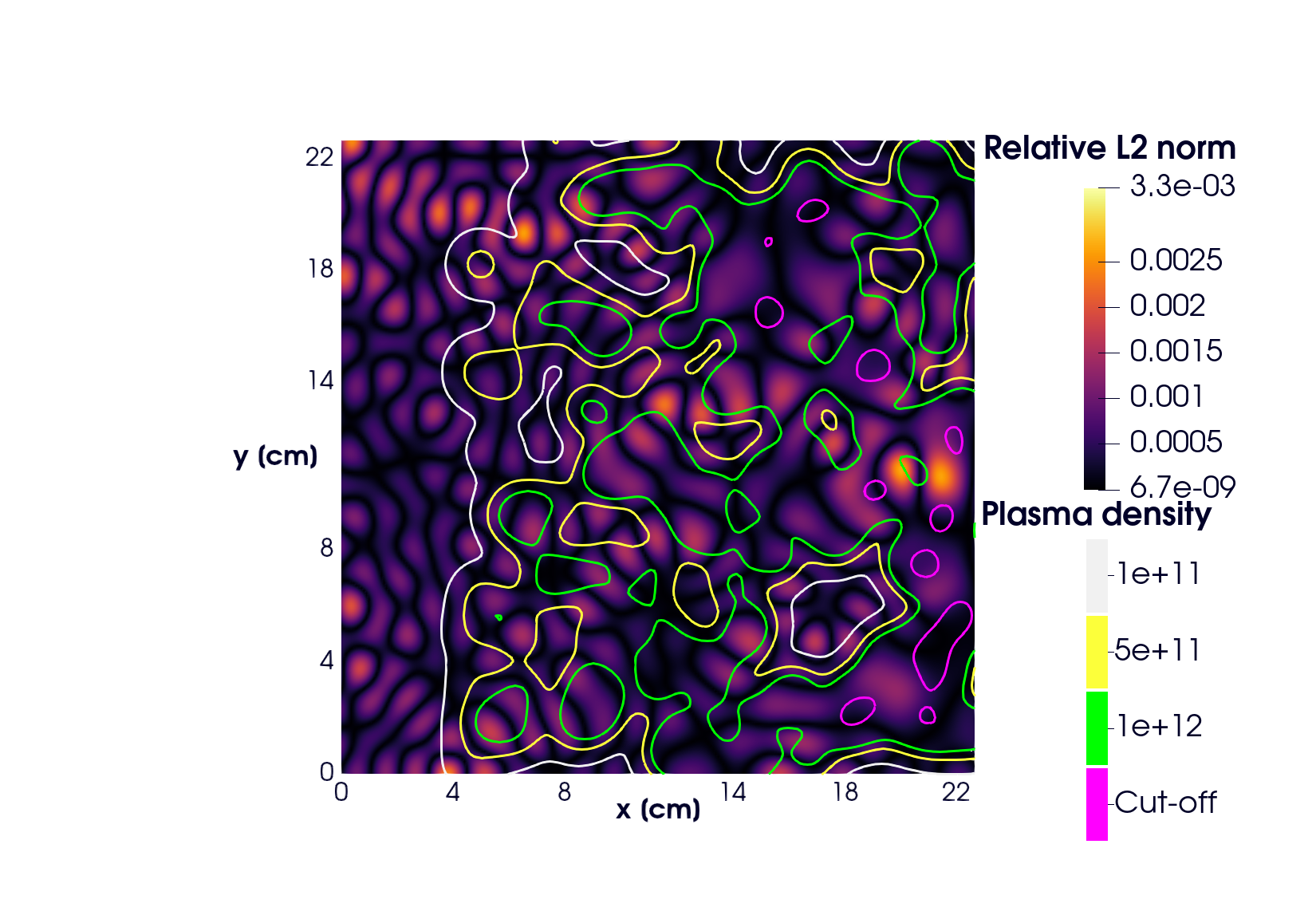}
         \caption{Time evolution of $\|\bs{R}\|_{L^2(\Omega)}$ (left) and pointwise absolute value of $\bs{R}(\bx,50T)$ (right).   }
         \label{pics_example_Omode_conv_timeharmonic}
 \end{figure}

In this experiment we consider an O-mode configuration, i.e. $\bs{e} = \bs{e}_z$ in \eqref{gaussian_beam_3d} and thus, $(\bE_h)_x = (\bE_h)_y = \bs{0}$. Figure \ref{pics_example_density} shows the electron plasma density and the small regions delimited by the cutoff contour (magenta contour), across which the wave cannot propagate.
The top row of figure \ref{pics_example_Omode_solution} shows the real (left) and imaginary (right) parts of the $z$-component of the frequency-domain solution $\bE^\mathrm{freq}$. The $x$ and $y$ components are zero. The bottom row shows the time evolution of $(\bE_h)_z$, the electric field of the time-domain solution (time increases from left to right). As expected, the beam propagates steadily until it reaches the printed contours, where the electron plasma density is non-negligible. Henceforth, the trajectory of the field is refracted by the high density regions, leading to a complicated solution pattern, which ressembles $\Re \{ \bE^\mathrm{freq} \}$.
 The left-hand side plot in figure \ref{pics_example_Omode_conv_timeharmonic} shows the time evolution of $\|\bs{R}\|_{L^2(\Omega)}$ defined in \eqref{rel_diff_conv_timeharmonic}. The right-hand side plot shows the pointwise value at the last time point, which corresponds to $|\bs{R}(\bx, 50T)|$.
The left plot in figure \ref{pics_example_Omode_conv_timeharmonic} shows that $\|\bs{R}\|_{L^2(\Omega)}$ decreases until it reaches a threshold of approximately $15\%$. In other words, the convergence of the electric field to the time-harmonic regime is clear in the beginning but stagnates after a large time. Moreover, the plot of $|\bs{R}(\bx, 50T)|$ shows that the difference is small and it is not concentrated in any particular region. The observations are coherent with the explanations given in remark \ref{remark_conv_timeharmonic}.

 \subsection{X-mode with a single blob of density}
\begin{figure}[H]
        \centering %[trim=left bottom right top, clip]
        \includegraphics[height=0.25\textheight, trim=0 2.5cm 0cm 4.8cm, clip]{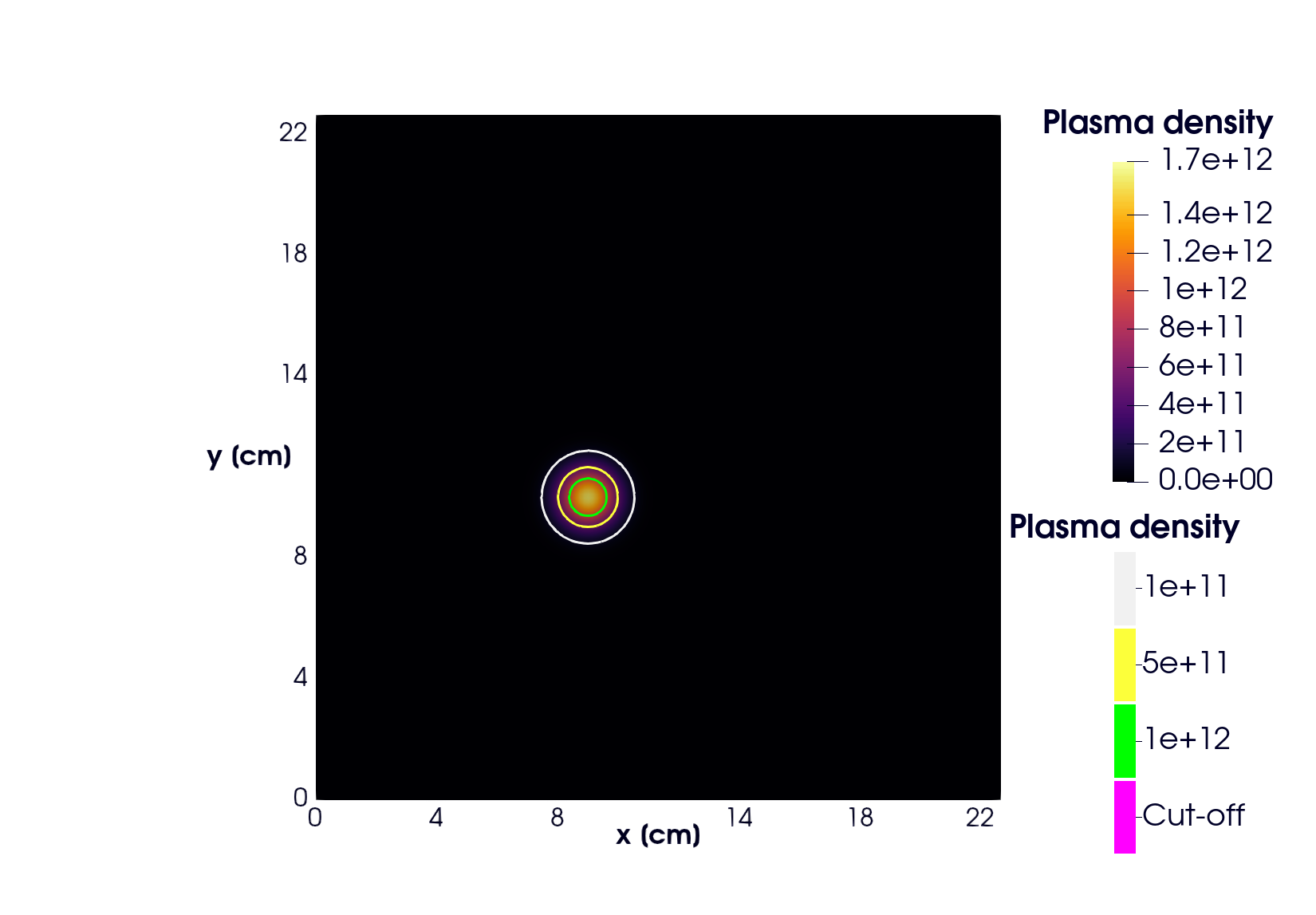}
        \caption{Electron-plasma density consisting of a single blob generated with a Gaussian function. The position is set to be slightly below the propagation axis of the beam. The maximum is below the cutoff value.}
        \label{pics_example_1lobXmode_density}
\end{figure}

\begin{figure}[H]
         \centering %[trim=left bottom right top, clip]
         \includegraphics[width=0.49\textwidth, trim=9cm 2.5cm 3cm 4.5cm, clip]{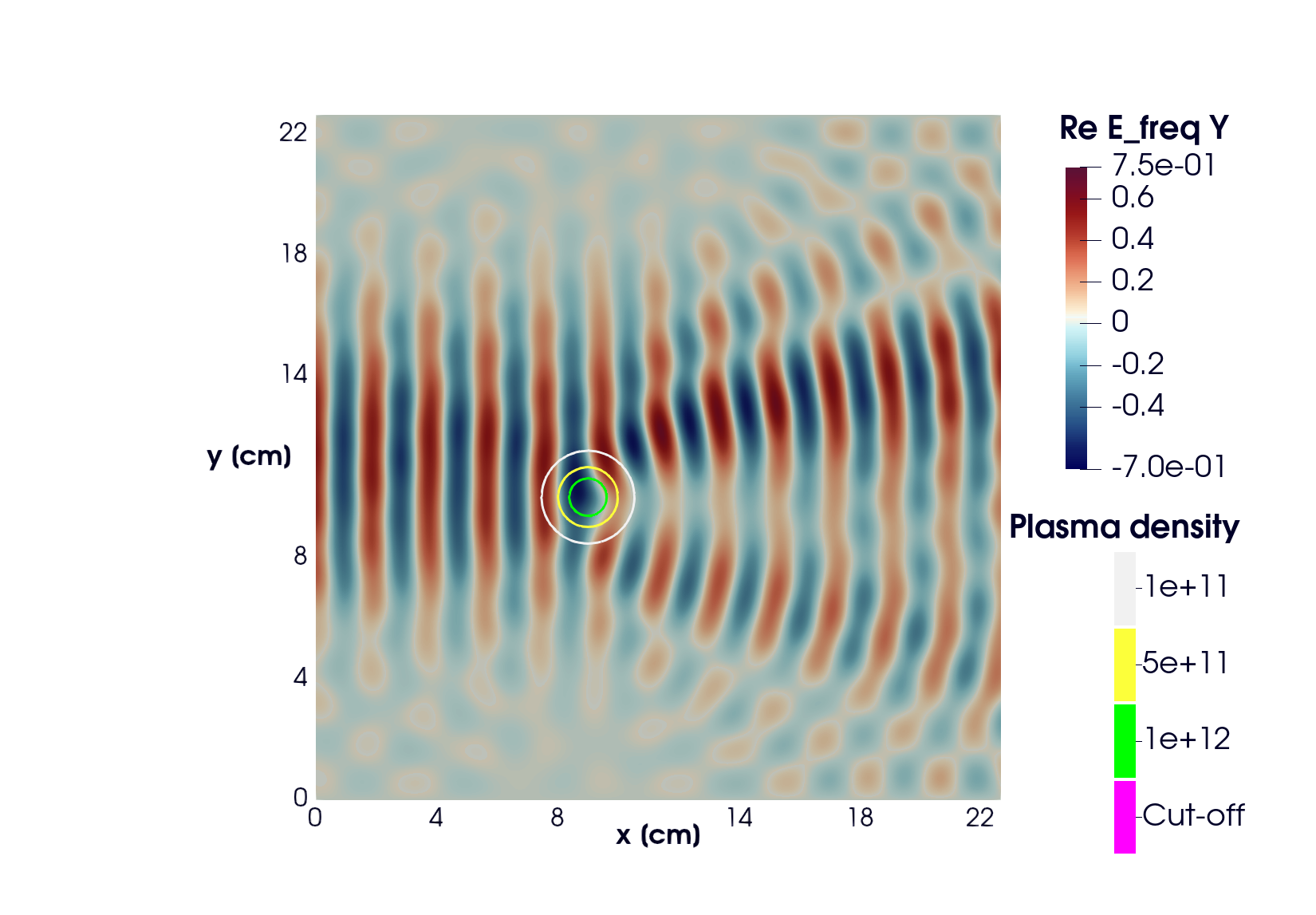}
         \includegraphics[width=0.49\textwidth, trim=9cm 2.5cm 3cm 4.5cm, clip]{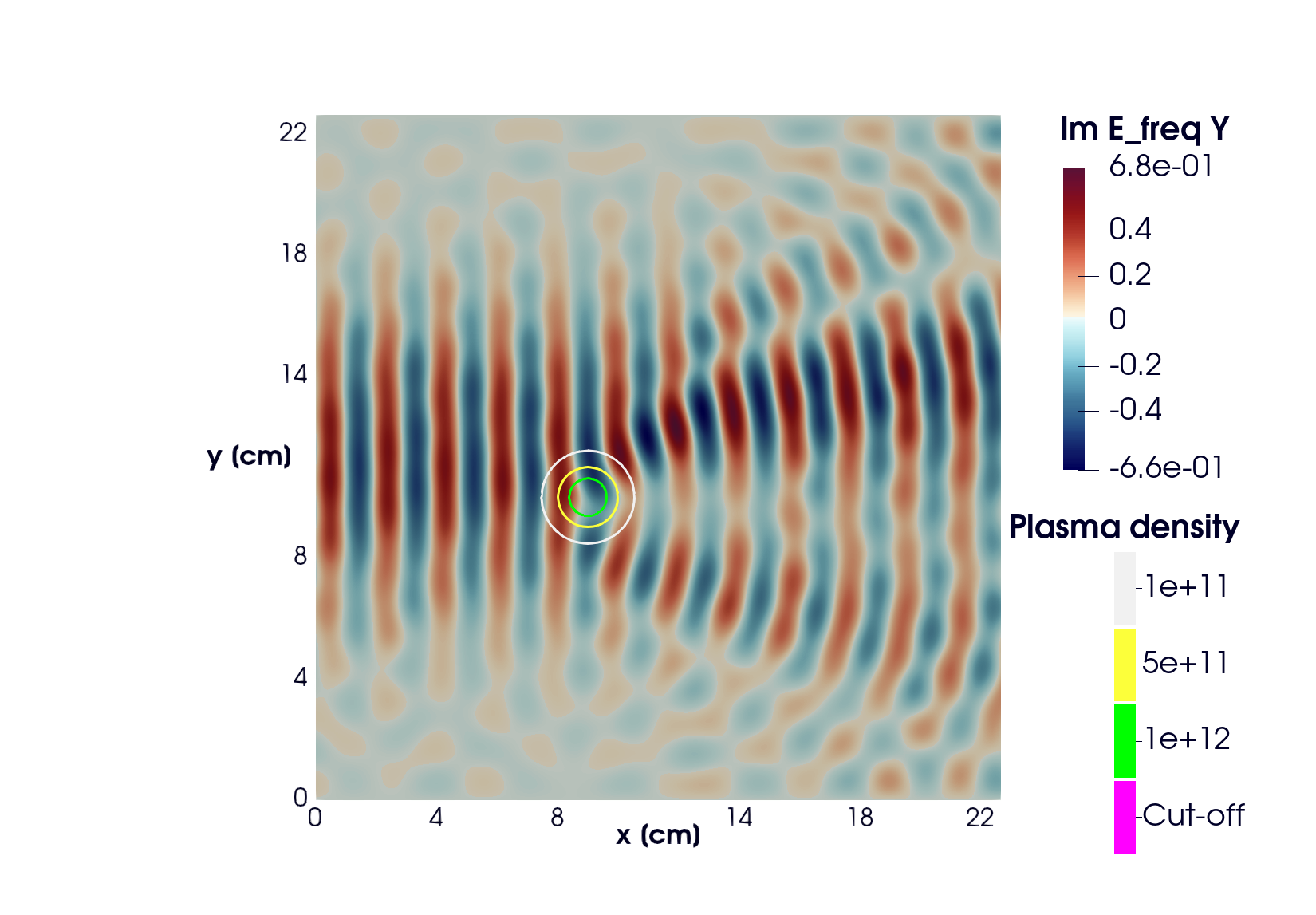}
         \includegraphics[height=0.25\textheight, trim=13.6cm 4.8cm 13.5cm 3.5cm, clip]{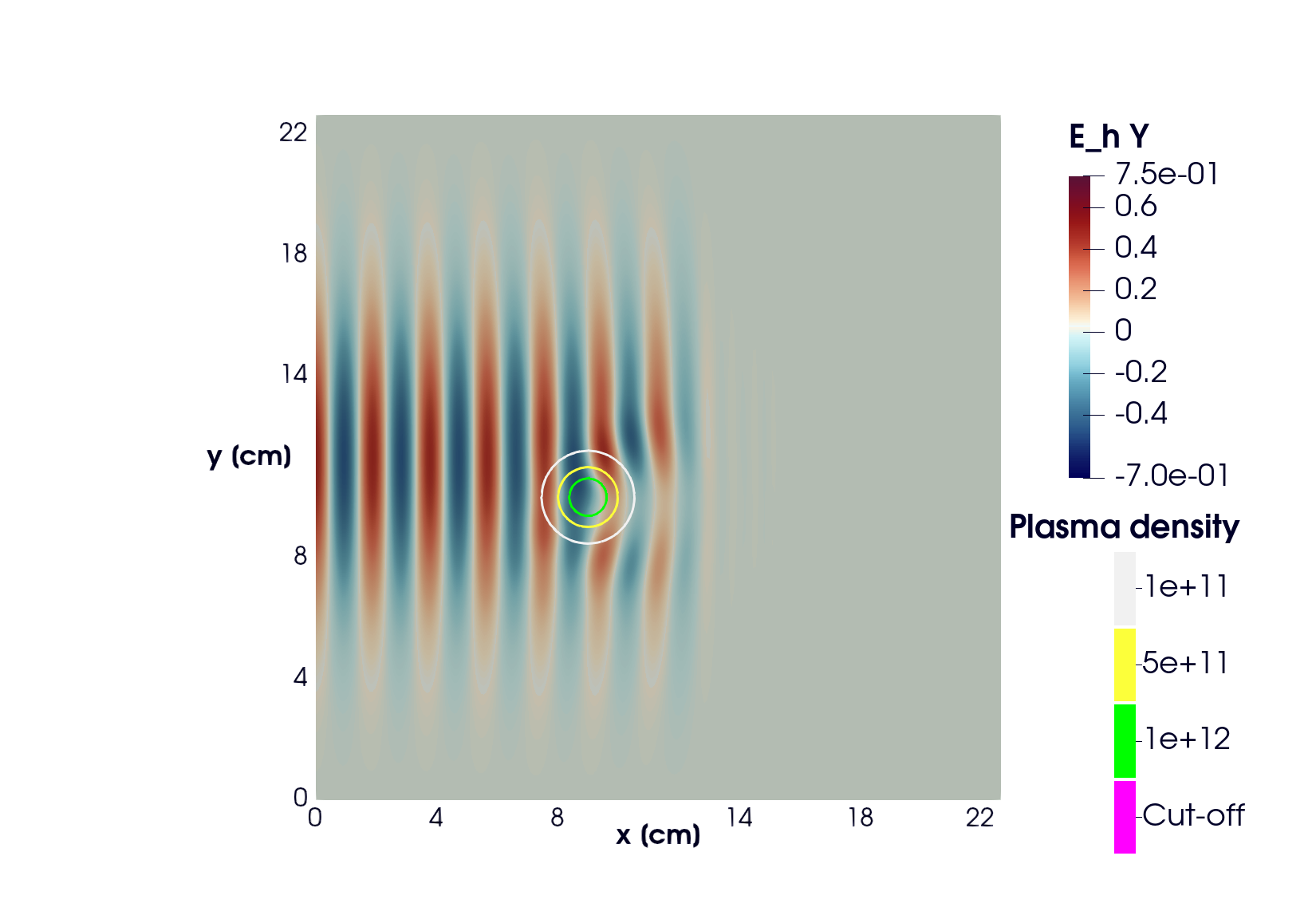}
         \includegraphics[height=0.25\textheight, trim=13.6cm 4.8cm 13.5cm 3.5cm, clip]{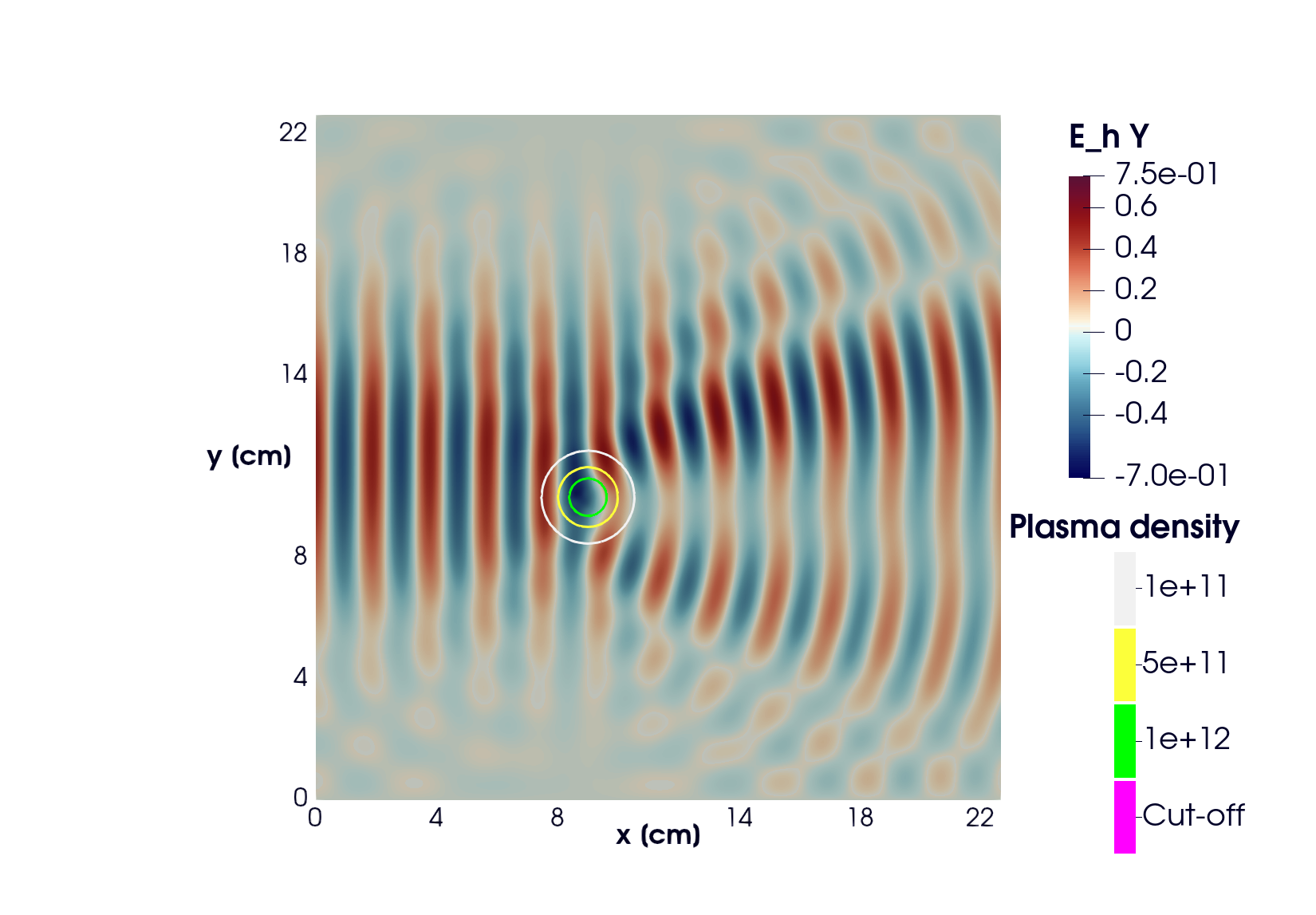}
         \caption{Top: real and imaginary part of $( \bE^{\mathrm{freq}} )_y$. Bottom: evolution of $( \bE_h )_y$ at $t = 7T, 15T$ (increasing from left to right). The color map is fixed to the amplitude of $( \Re \{\bE^{\mathrm{freq}}\} )_y$, which is $[-0.7, 0.75]$.}
         \label{pics_example_1blobXmode_solution_Y}
\end{figure}

\begin{figure}[H]
        \centering
         \includegraphics[width=0.49\textwidth, trim=9.1cm 2.5cm 3.1cm 4.5cm, clip]{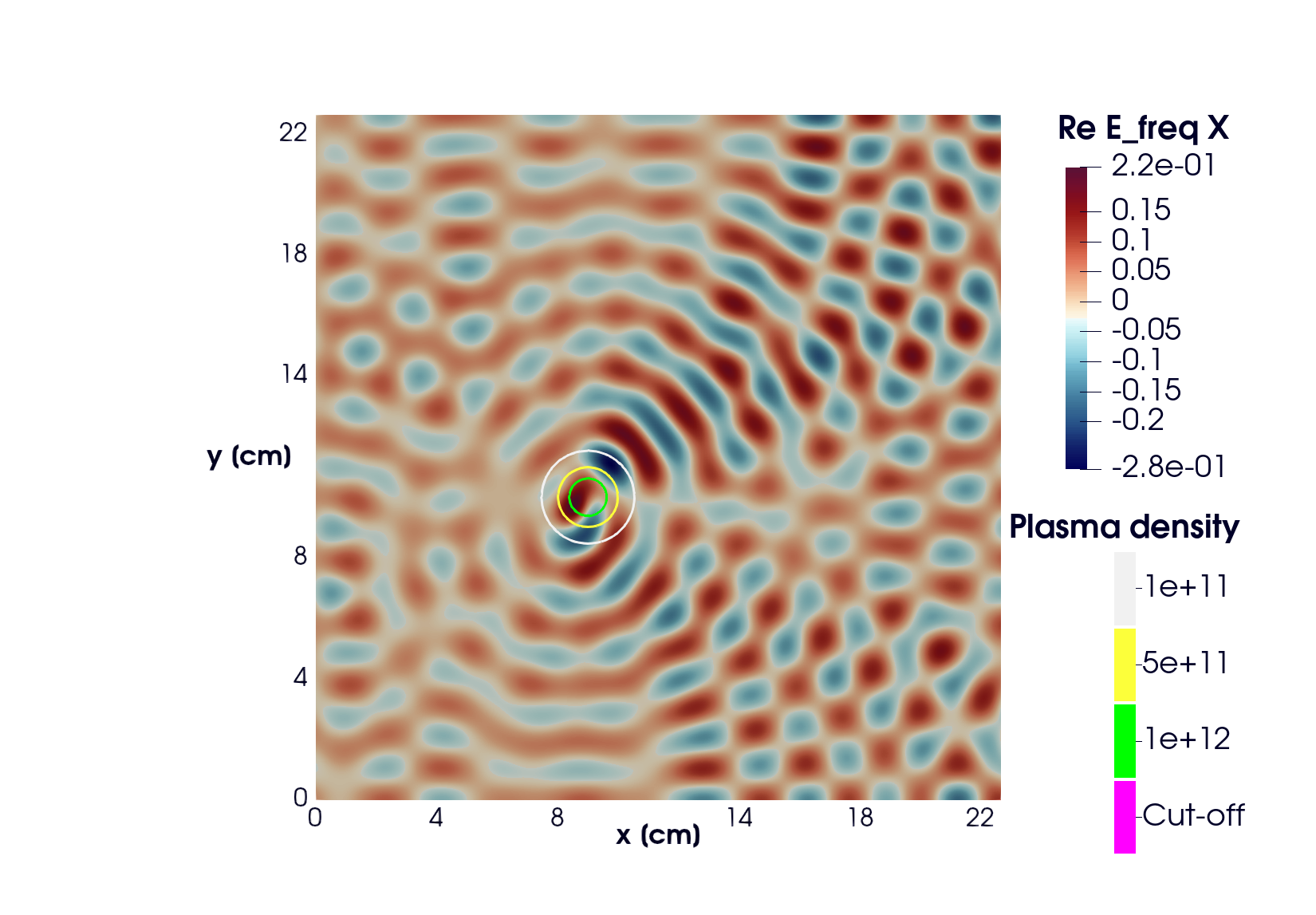}
        \includegraphics[width=0.49\textwidth, trim=9.1cm 2.5cm 3.1cm 4.5cm, clip]{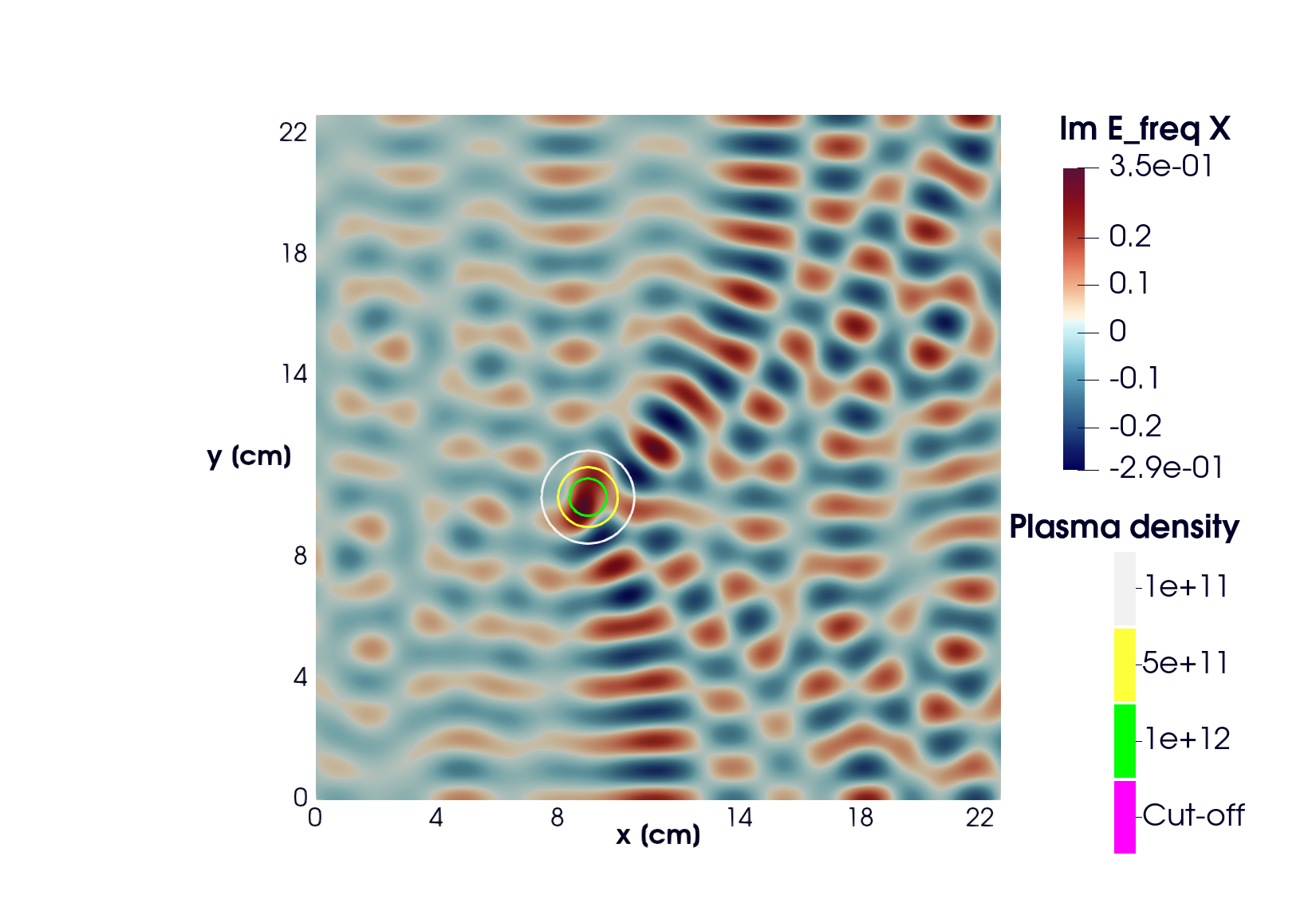}
        \includegraphics[height=0.25\textheight, trim=13.6cm 4.8cm 13.5cm 3.5cm, clip]{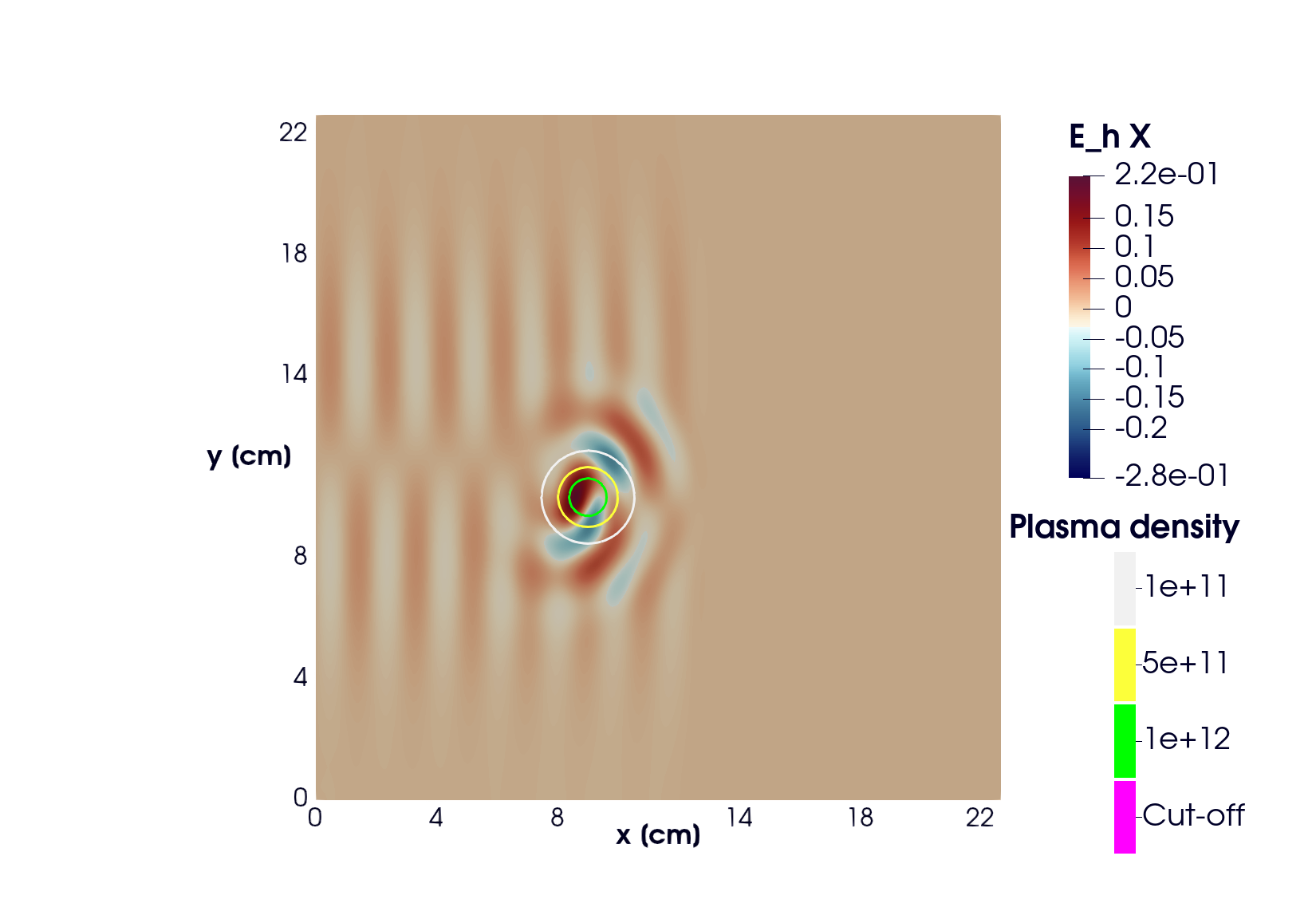}
        \includegraphics[height=0.25\textheight, trim=13.6cm 4.8cm 13.5cm 3.5cm, clip]{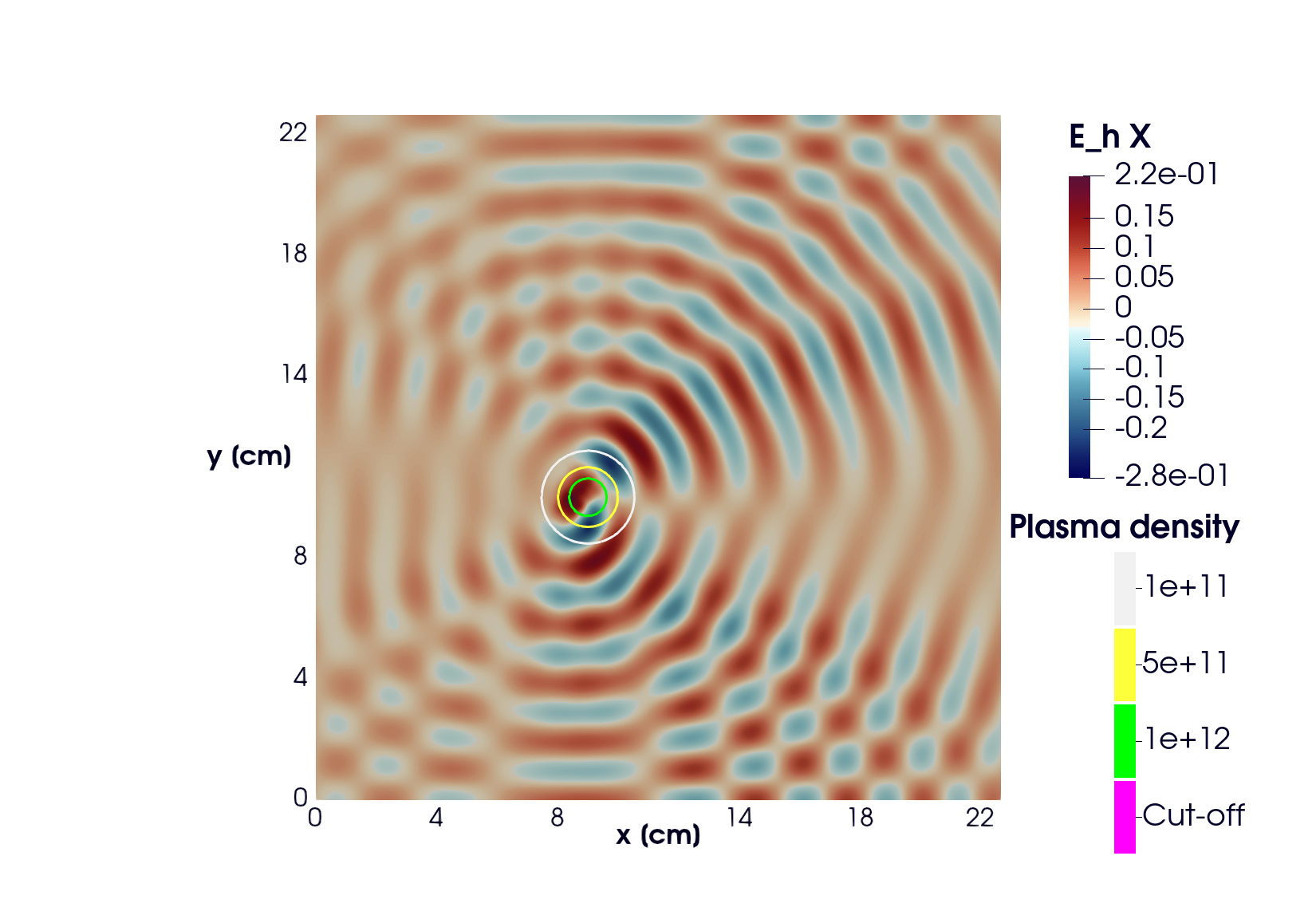}
         \caption{Top: real and imaginary part of  $( \bE^{\mathrm{freq}} )_x$. Bottom: evolution of $( \bE_h )_x$ at $t = 7T, 15T$ (increasing from left to right). The color map is fixed to the amplitude of $( \Re \{\bE^{\mathrm{freq}}\} )_x$, which is $[-0.28, 0.22]$.
         }
         \label{pics_example_1blobXmode_solution_X}
 \end{figure}

 In this experiment we consider an X-mode configuration, i.e. $\bs{e} = \bs{e}_y$ in \eqref{gaussian_beam_3d} and thus, $(\bE_h)_z = \bs{0}$. The goal of this test is to check that the method does not produce unphysical asymmetries. For that, we construct the electron-plasma density profile shown in figure \ref{pics_example_1lobXmode_density}, which consists of a single plasma blob centered slightly below the beam propagation axis. The maximum is lightly below the cutoff to allow the wave to propagate through.
Figures \ref{pics_example_1blobXmode_solution_Y} ($y$-component) and \ref{pics_example_1blobXmode_solution_X} ($x$-component) show the real and imaginary parts of $\bE^\mathrm{freq}$ (top) and the time-evolution of $\bE_h$ (bottom).

In figure \ref{pics_example_1blobXmode_solution_Y} we can see that the $y$-component of the field is split into two branches at the high-density blob. Since the blob is beneath the beam axis, the amplitude of the superior branch is higher, thus the solution reflects the natural asymmetry in the problem and the code does not produce unphysical asymmetries.
Since the launched beam is polarized in $y$, the $x$-component shown in figure \ref{pics_example_1blobXmode_solution_X} only gains intensity at the high-density blob and shows the gyrating effect caused by the interaction with the plasma.
We observe that the frequency-domain field exhibits a visible interference pattern. This can be caused by the reflections at the boundary produced by the Silver-Müller boundary conditions when the wavefront is not parallel to the boundary, which is the case almost everywhere in the boundary in this experiment. We observe that the time-domain field also develops this intereference pattern, but to a lesser extent.

\begin{figure}[H]
          \centering %[trim=left bottom right top, clip]
          \includegraphics[width=0.4\textwidth, trim=0.8cm 1cm 0.8cm 4.6cm, clip]{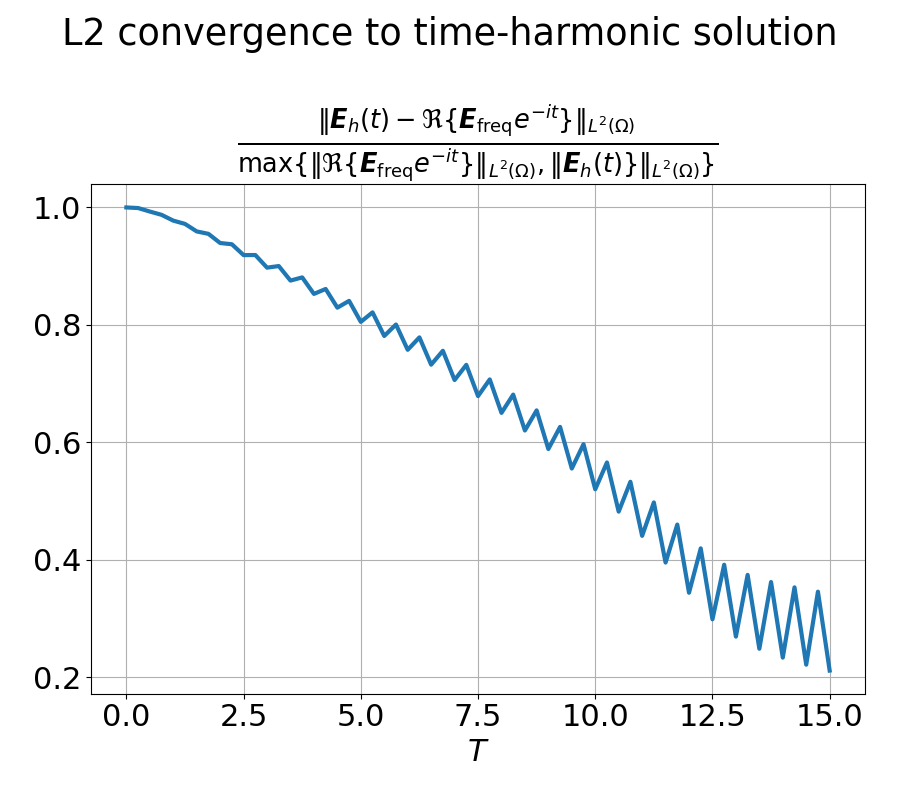}
          \includegraphics[width=0.59\textwidth, trim=9cm 2.5cm 2.3cm 4.5cm, clip]{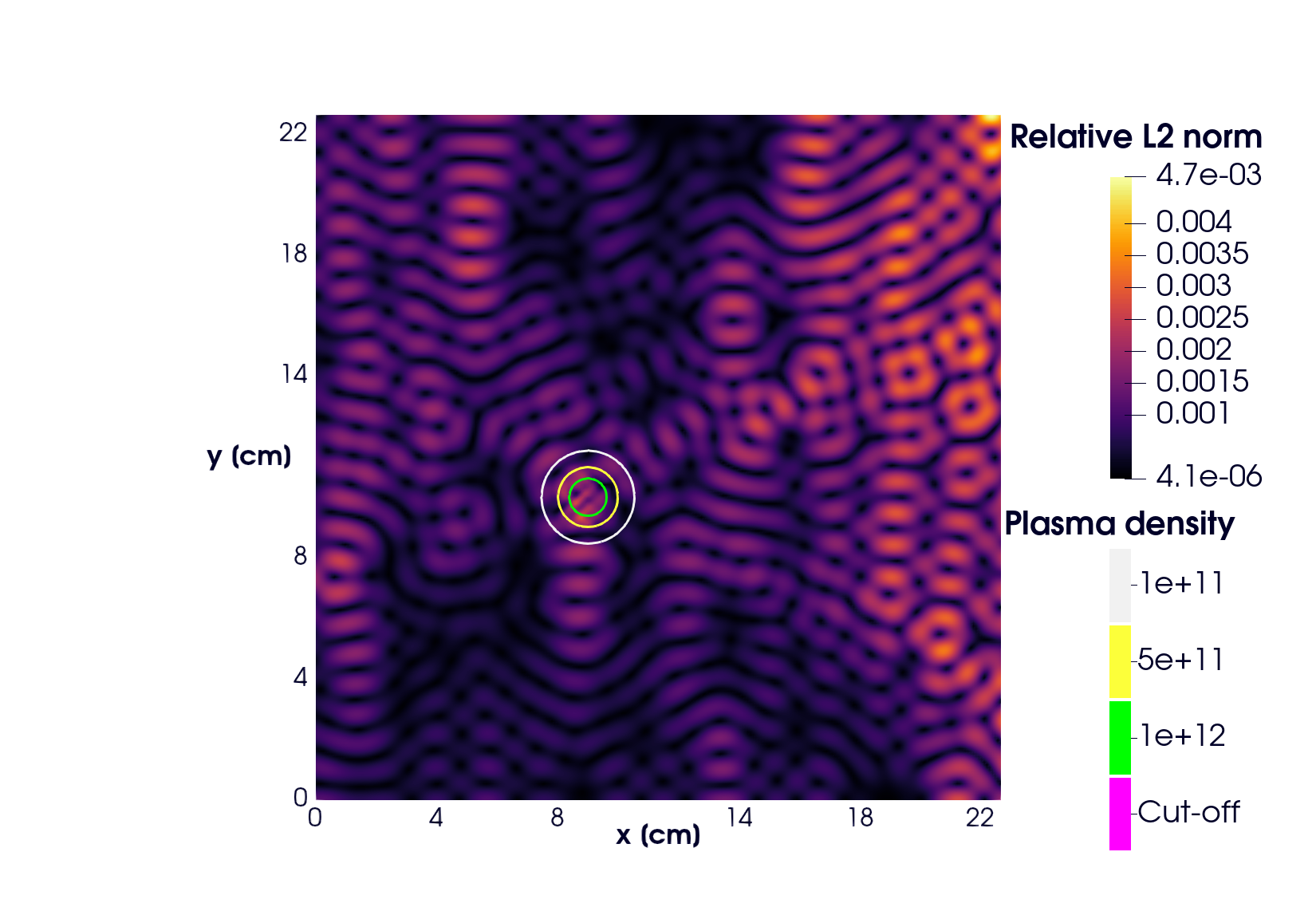}
         \caption{Time evolution of $\|\bs{R}\|_{L^2(\Omega)}$ (left) and pointwise Euclidian norm of $\bs{R}(\bx,15T)$ (right).
         }
         \label{pics_example_1lobXmode_conv_timeharmonic}
 \end{figure}

The left-hand side plot in figure \ref{pics_example_1lobXmode_conv_timeharmonic} shows the time evolution of $\|\bs{R}\|_{L^2(\Omega)}$. The right-hand side panel shows the pointwise Euclidian norm of the difference at the last time point, which is $t=15T$.
Figure \ref{pics_example_1lobXmode_conv_timeharmonic} shows that, as before, the convergence of $\|\bs{R}\|_{L^2(\Omega)}$ stagnates after a certain time. Specifically in this experiment we observe oscillations that grow large as time increases. This could be caused by a frequency mismatch between the two solution components, since the $x$-component is only excited after a certain number of periods, where the phase is distorted by the accumulated splitting error.  The right plot in figure \ref{pics_example_1lobXmode_conv_timeharmonic} shows that the difference at the last simulation point is slightly larger close to the boundary. This can indicate that the solution requires more time to converge or that the modes with different frequencies are not able to propagate outwards through the boundary $\Gamma_A$.

  \subsection{X-mode with high density fluctuations}
\begin{figure}[H]
        \centering %[trim=left bottom right top, clip]
        \includegraphics[height=0.3\textheight, trim=0 2.5cm 0cm 4.5cm, clip]{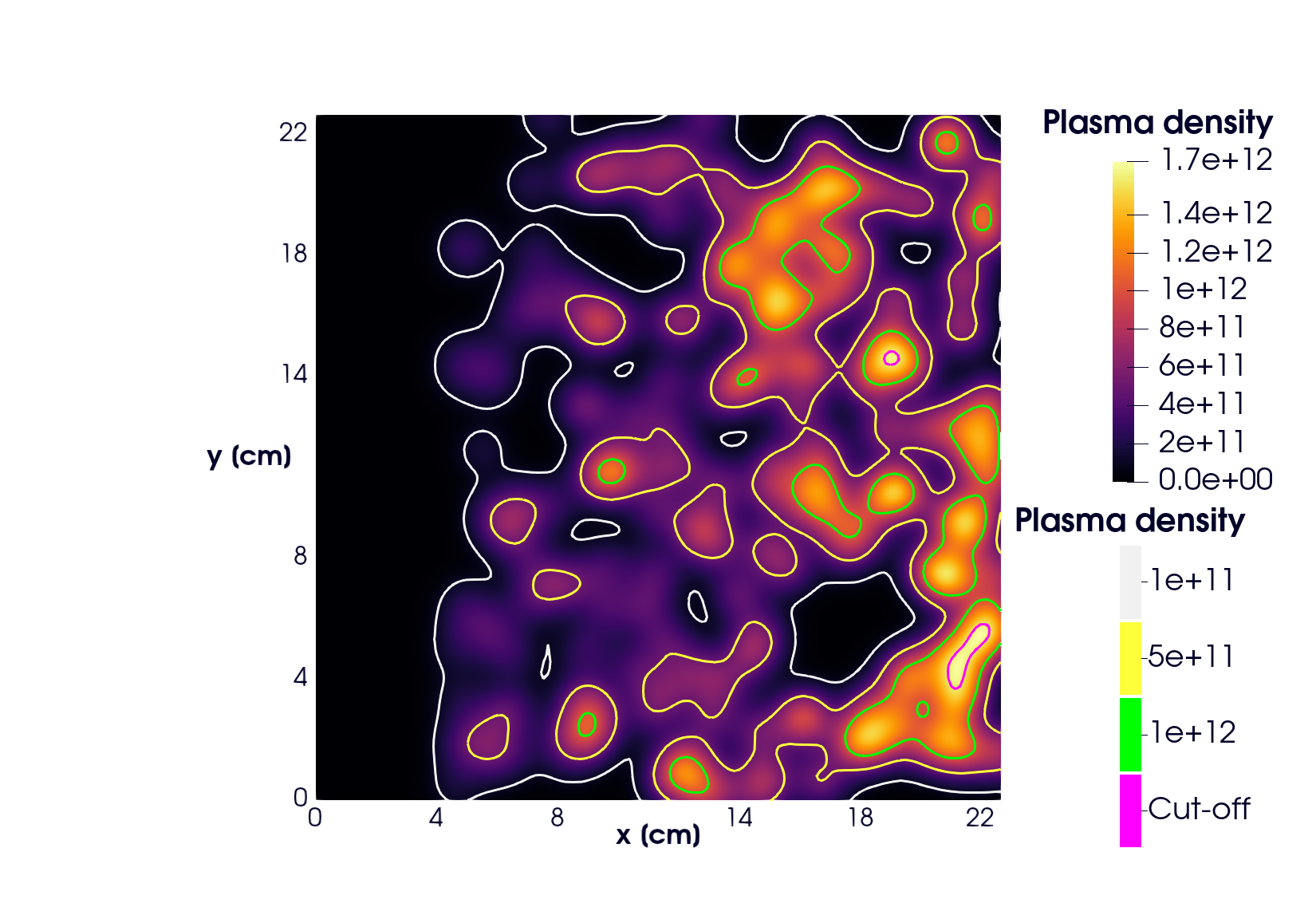}
        \caption{Electron plasma density with high fluctuations. The profile does not reach the upper-hybrid resonance. The X-mode cutoff contour (magenta contour) is given by the relation $\nomegap^2=1 - \nomegac$ \cite{chen2013introduction}, which corresponds approximately to $n_e(\bx) = 1.6 \times 10^{12}$ electrons per cubic centimeter.}
        \label{pics_example_high_fluctuations_Xmode_density}
\end{figure}

\begin{figure}[H]
         \centering %[trim=left bottom right top, clip]
         \includegraphics[width=0.49\textwidth, trim=9cm 2.5cm 3cm 4.5cm, clip]{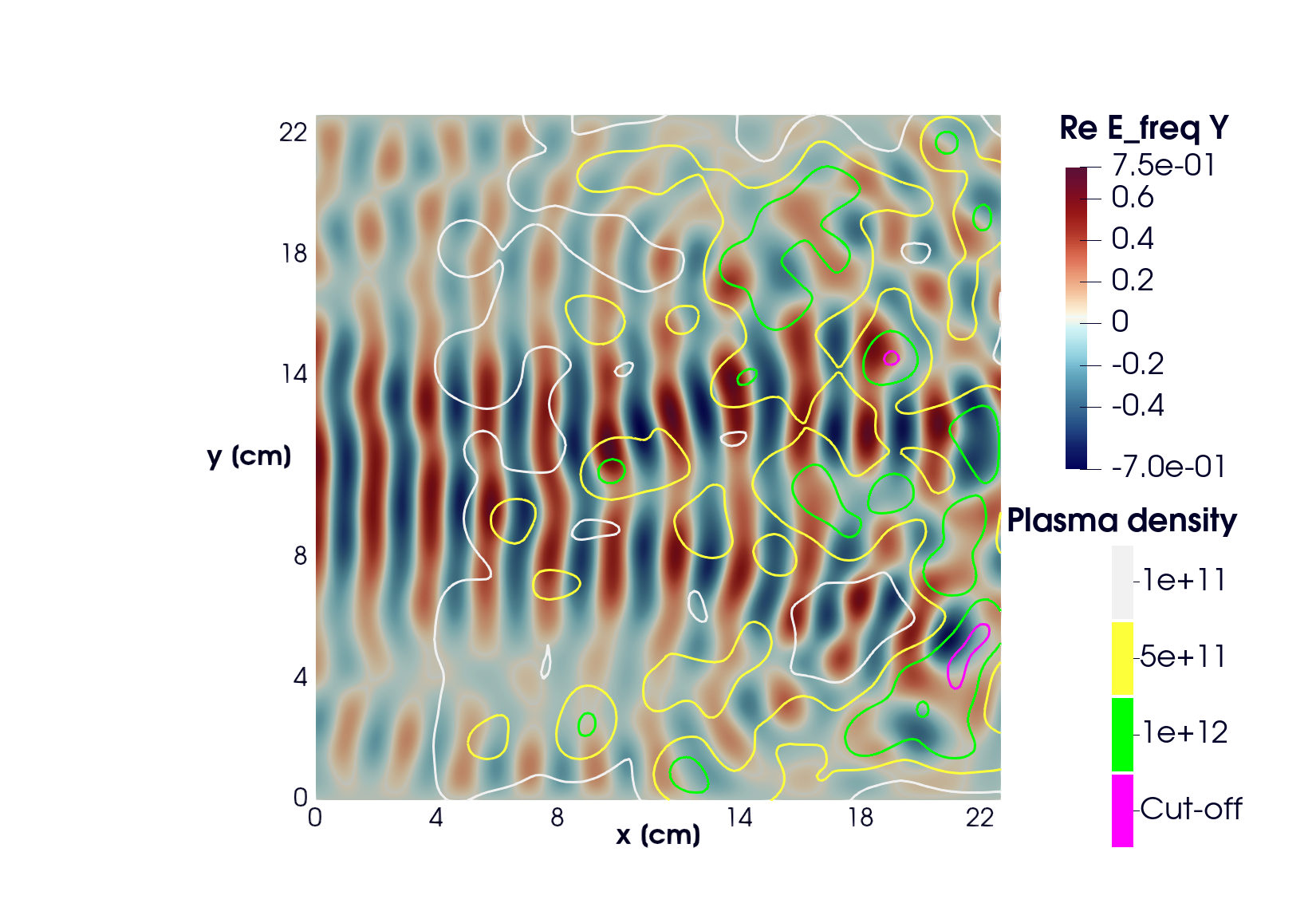}
        \includegraphics[width=0.49\textwidth, trim=9cm 2.5cm 3cm 4.5cm, clip]{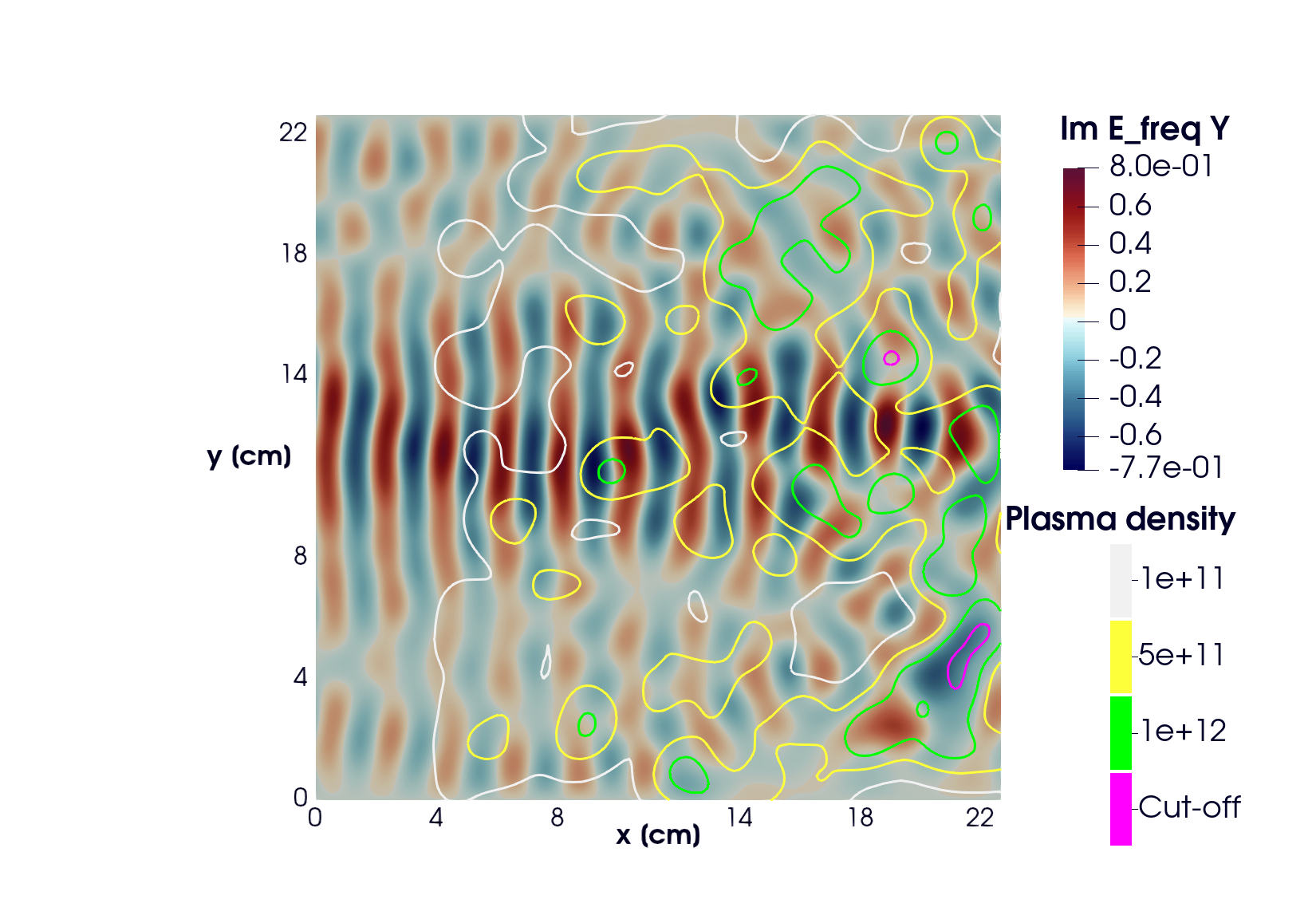}
          \includegraphics[width=0.194\textwidth, trim=14cm 4.8cm 13.9cm 3.5cm, clip]{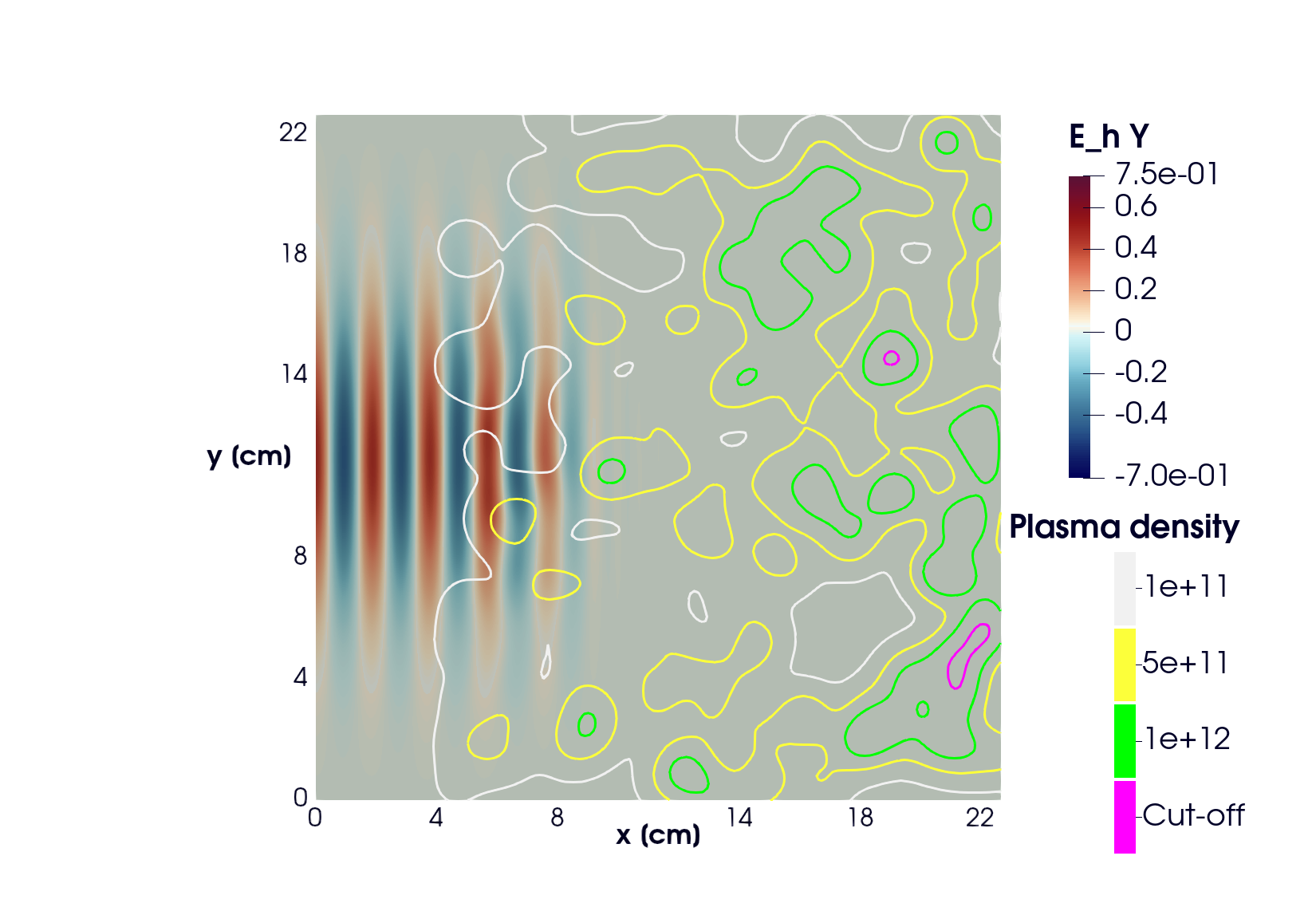}
          \includegraphics[width=0.194\textwidth, trim=14cm 4.8cm 13.9cm 3.5cm, clip]{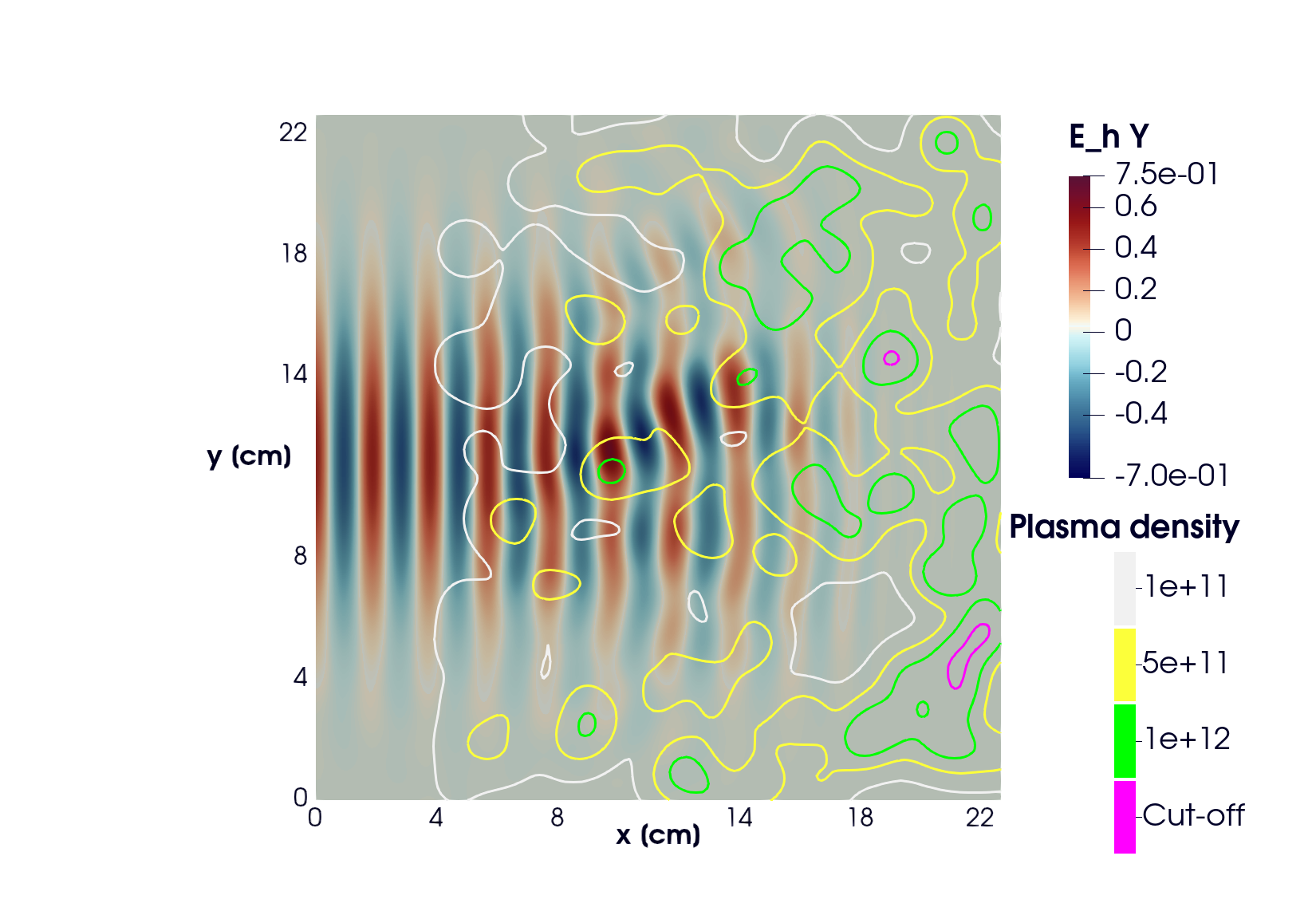}
          \includegraphics[width=0.194\textwidth, trim=14cm 4.8cm 13.9cm 3.5cm, clip]{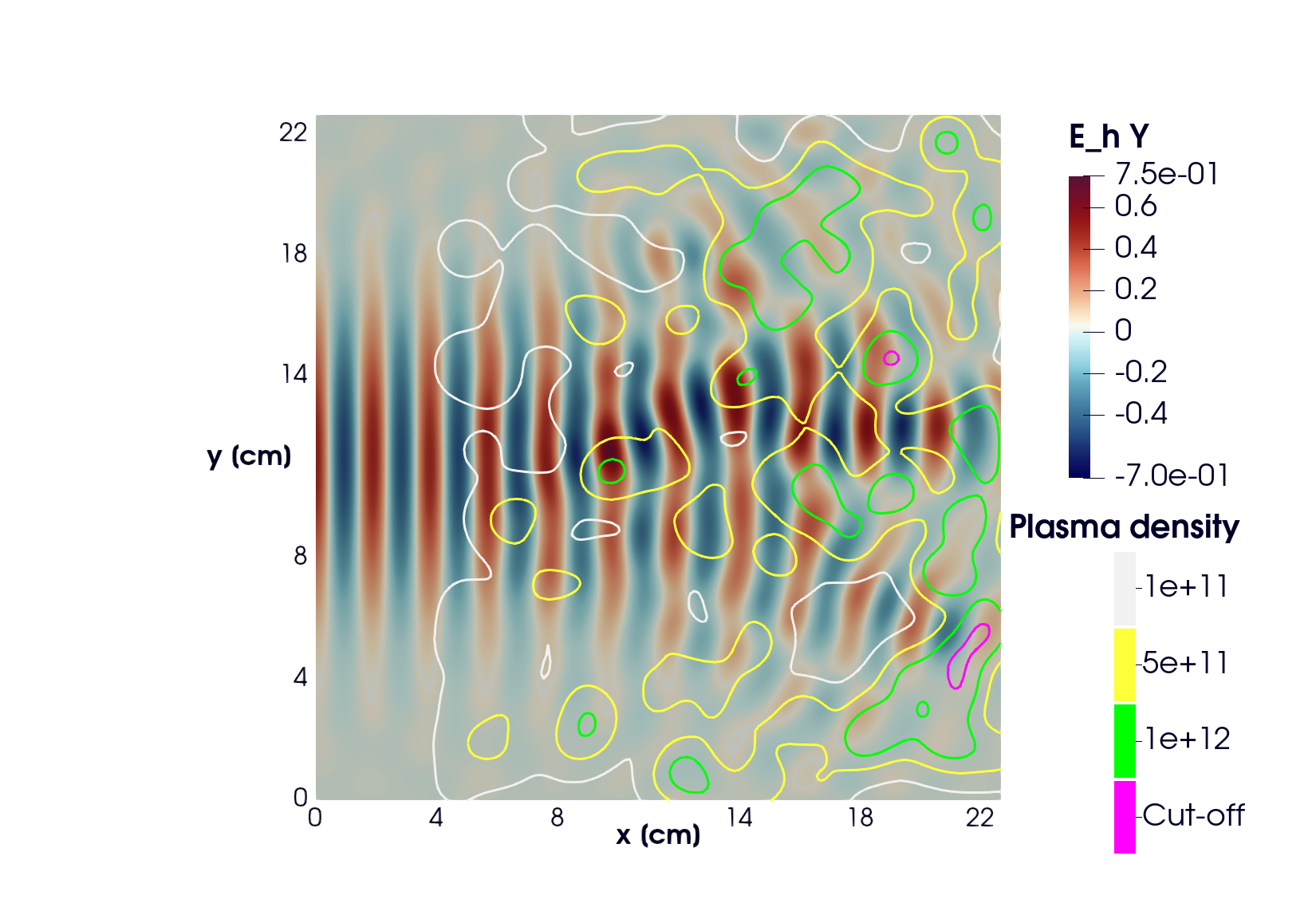}
          \includegraphics[width=0.194\textwidth, trim=14cm 4.8cm 13.9cm 3.5cm, clip]{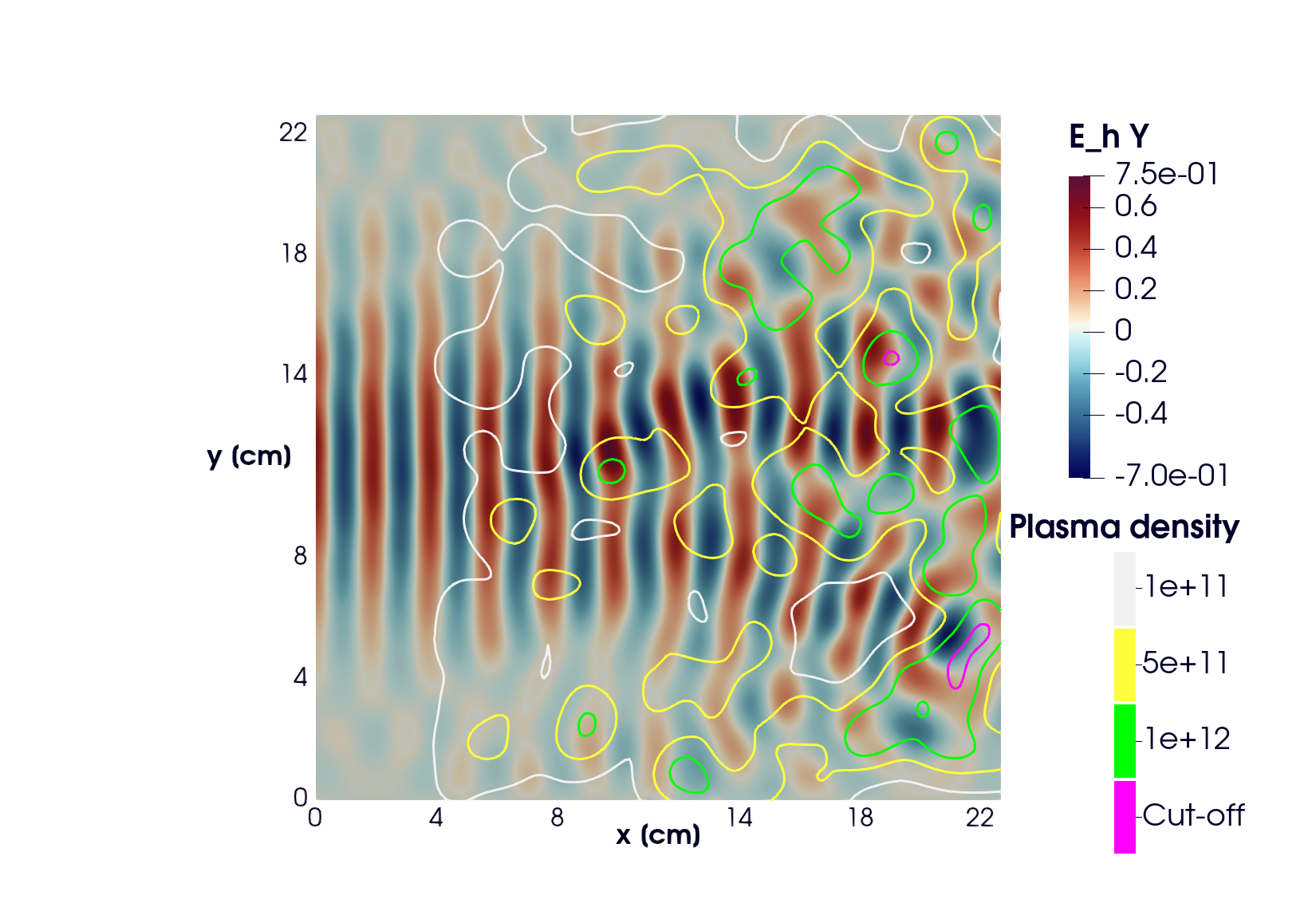}
          \includegraphics[width=0.194\textwidth, trim=14cm 4.8cm 13.9cm 3.5cm, clip]{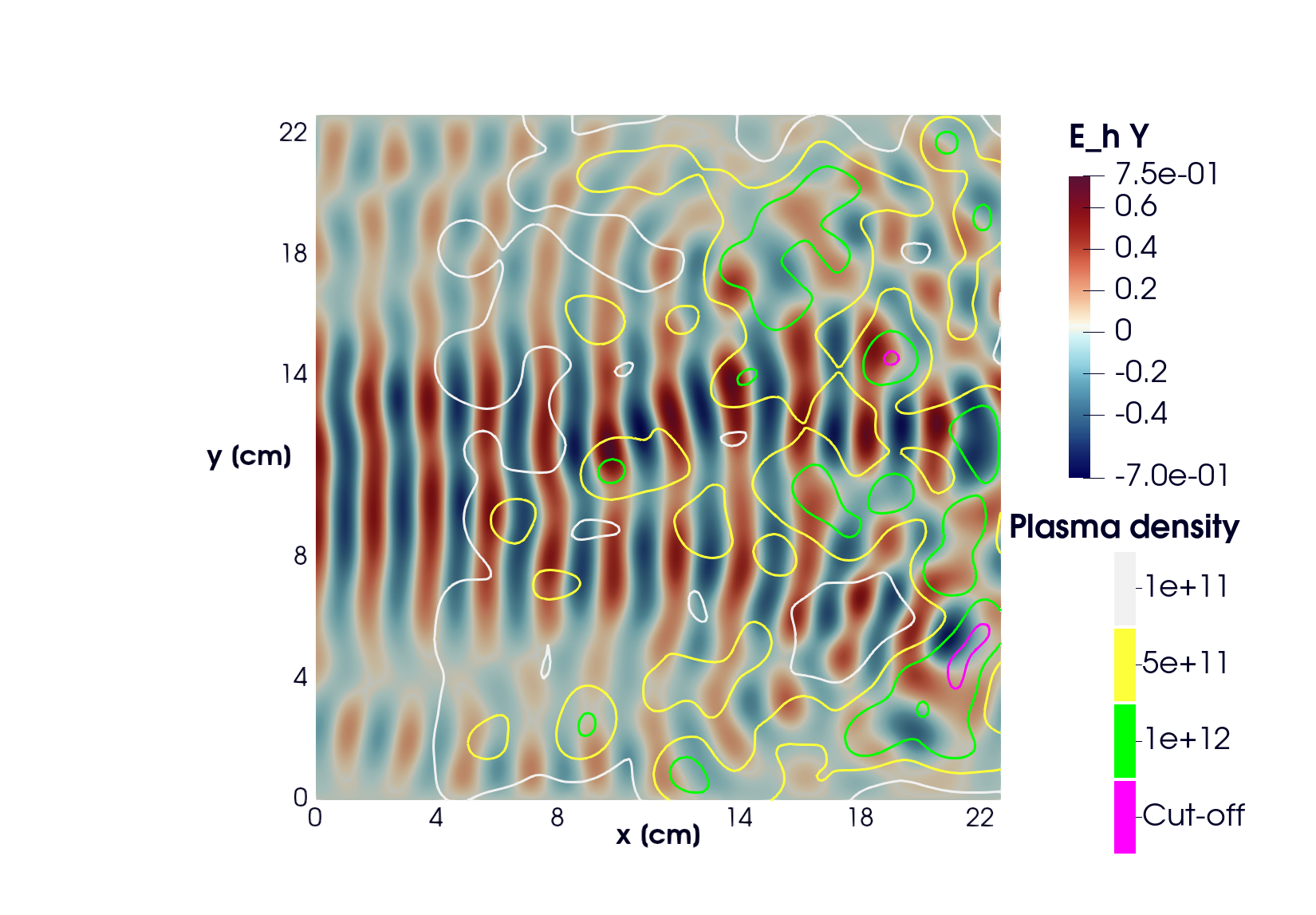}
         \caption{Top: real and imaginary part of $( \bE^{\mathrm{freq}} )_y$. Bottom: evolution of $( \bE_h )_y$ at $t = 5T,10T,15T,28T, 50T$ (increasing from left to right). The color map is fixed to the amplitude of $(\Re \{\bE^{\mathrm{freq}}\})_y$, which is $[-0.7, 0.75]$.}
         \label{pics_example_high_fluctuations_Xmode_solution_Y}
 \end{figure}

\begin{figure}[H]
          \centering %[trim=left bottom right top, clip]
         \includegraphics[width=0.49\textwidth, trim=9cm 2.5cm 3cm 4.5cm, clip]{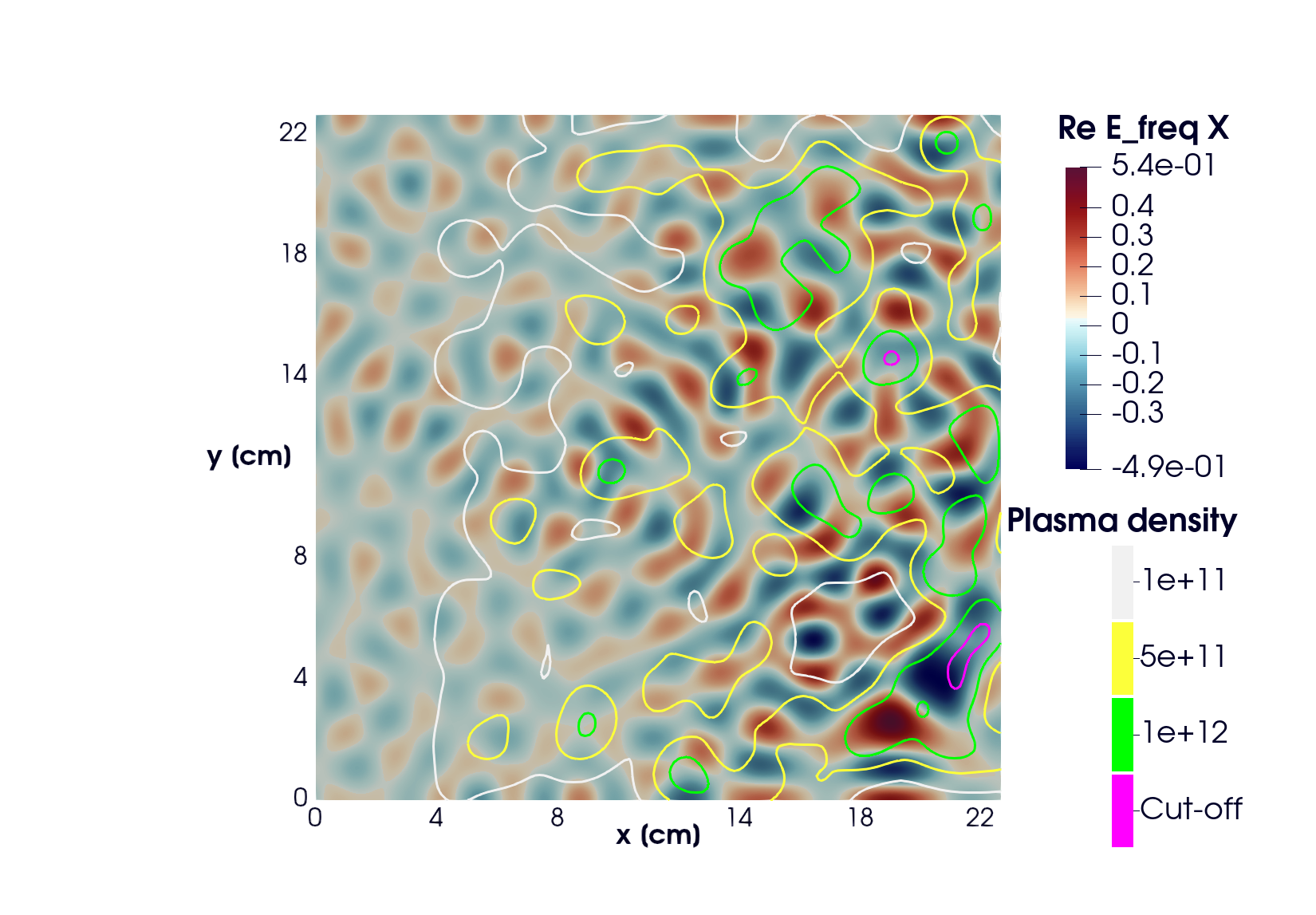}
        \includegraphics[width=0.49\textwidth, trim=9cm 2.5cm 3cm 4.5cm, clip]{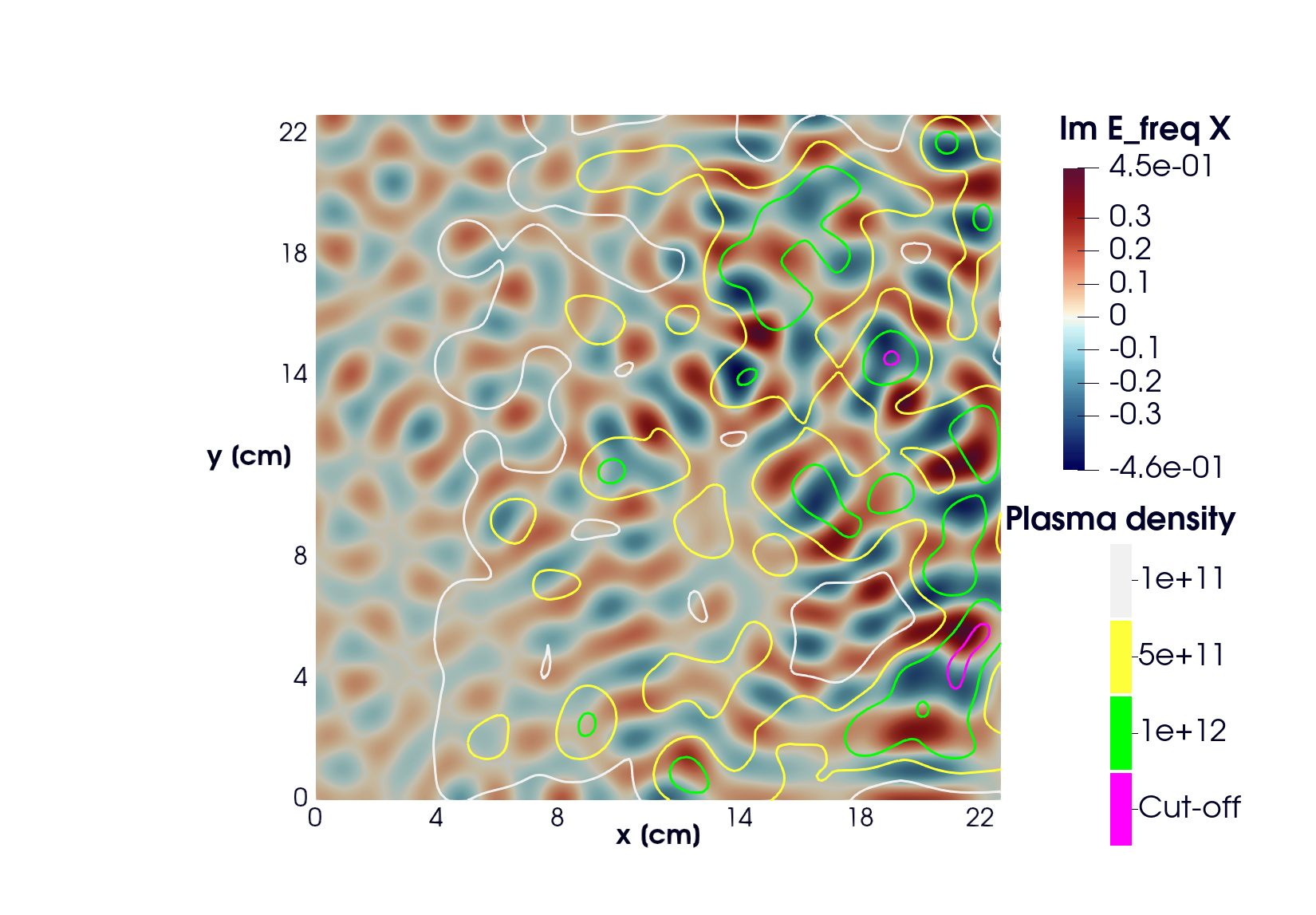}
          \includegraphics[width=0.194\textwidth, trim=14cm 4.8cm 13.9cm 3.5cm, clip]{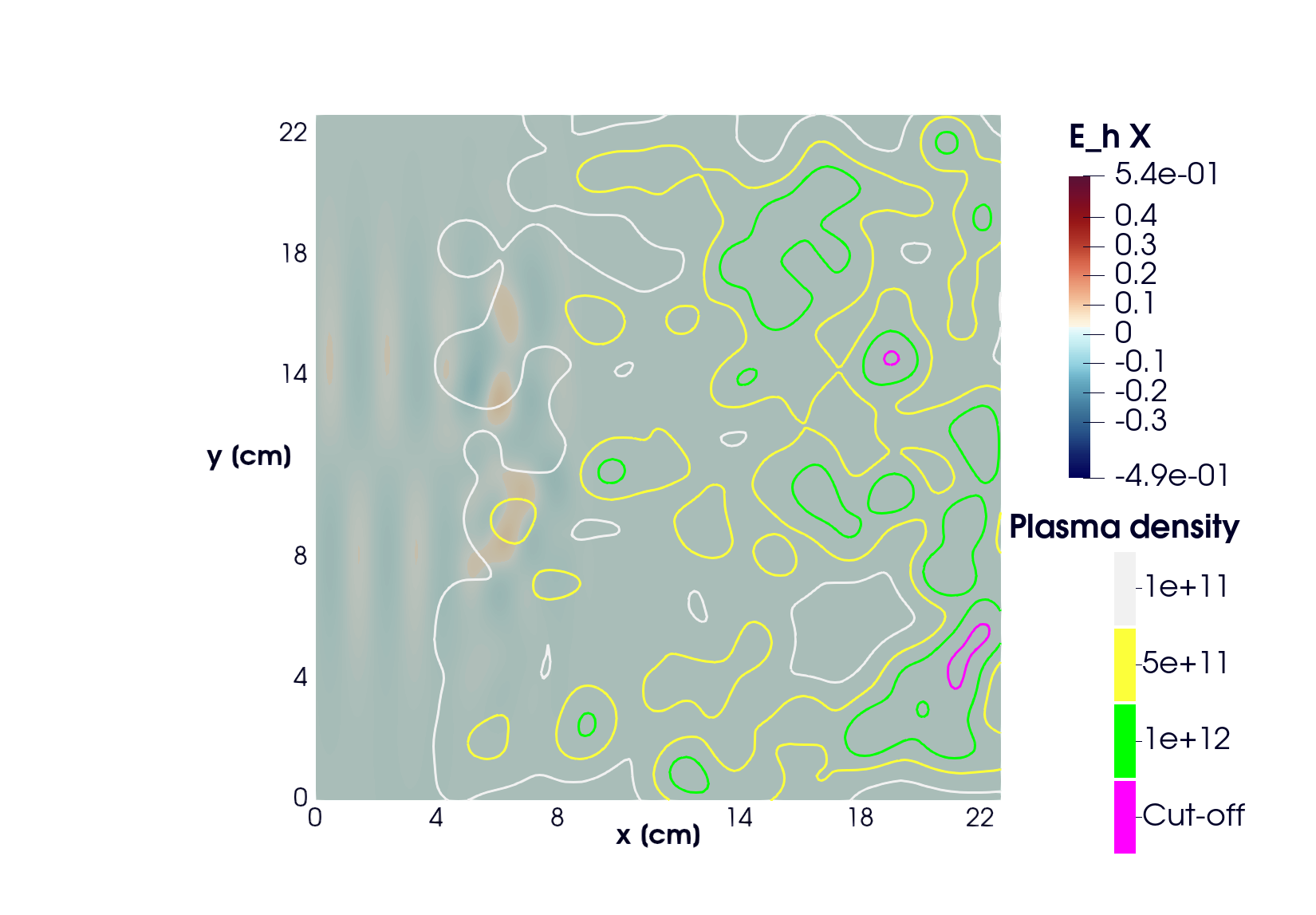}
          \includegraphics[width=0.194\textwidth, trim=14cm 4.8cm 13.9cm 3.5cm, clip]{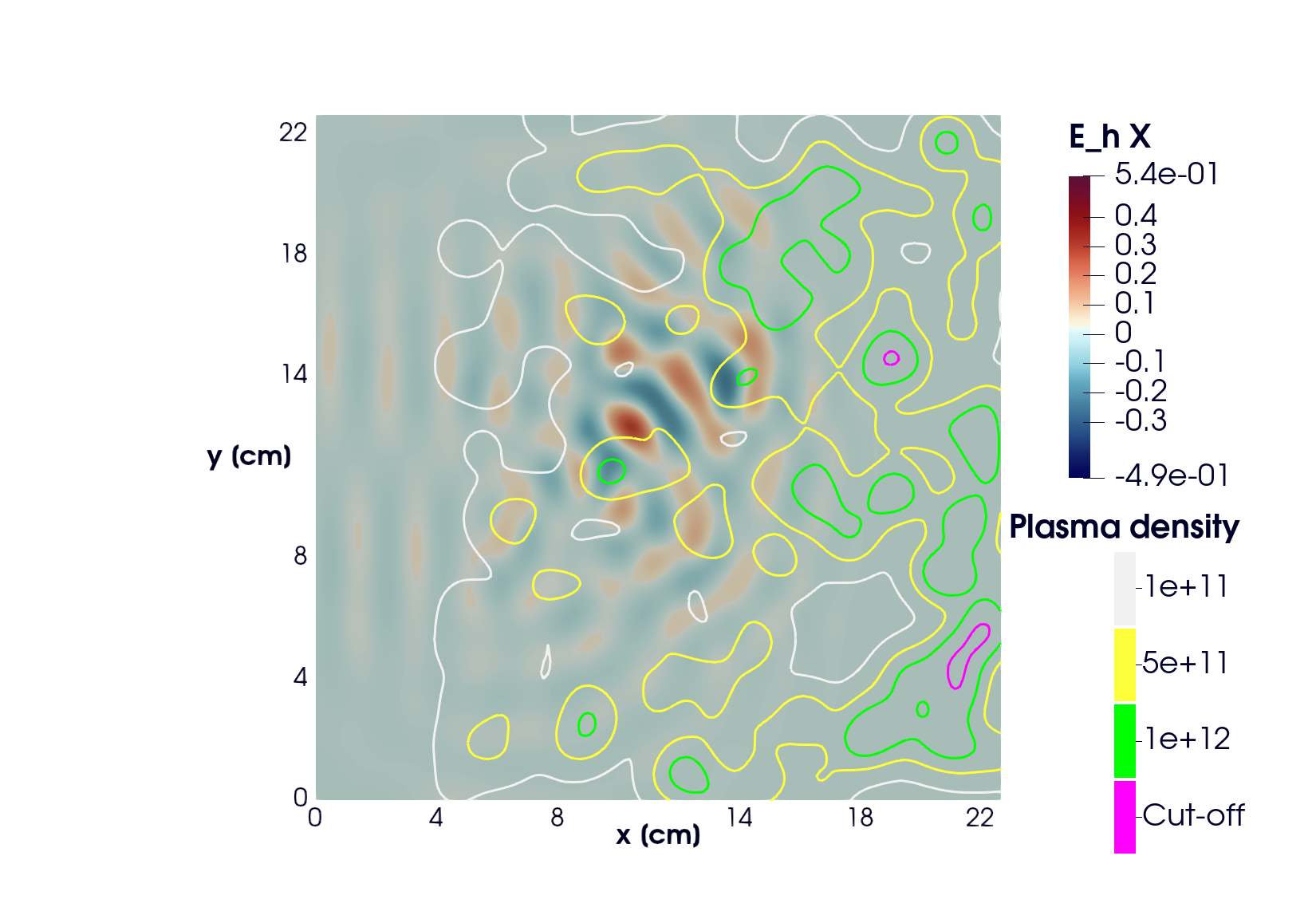}
          \includegraphics[width=0.194\textwidth, trim=14cm 4.8cm 13.9cm 3.5cm, clip]{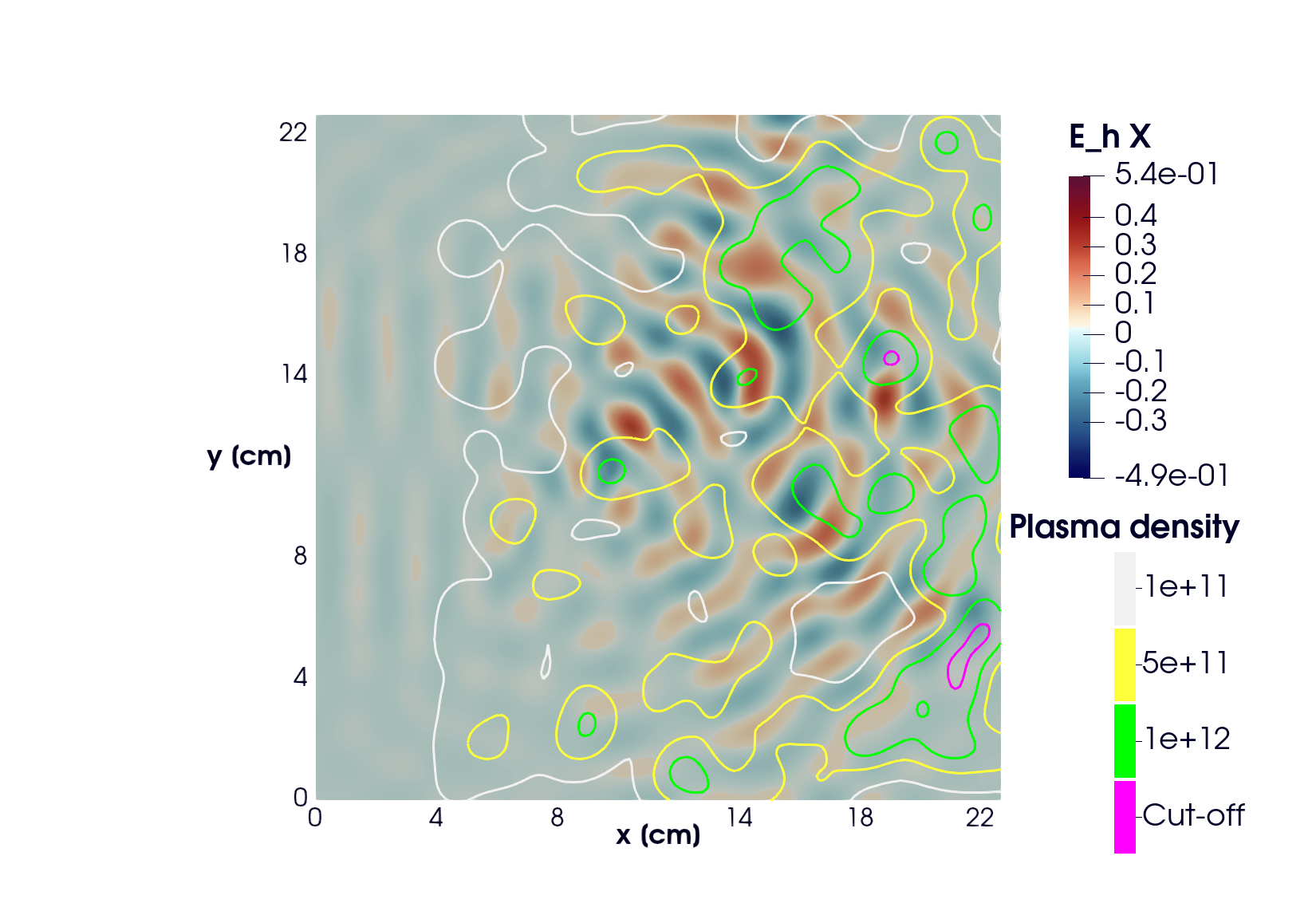}
          \includegraphics[width=0.194\textwidth, trim=14cm 4.8cm 13.9cm 3.5cm, clip]{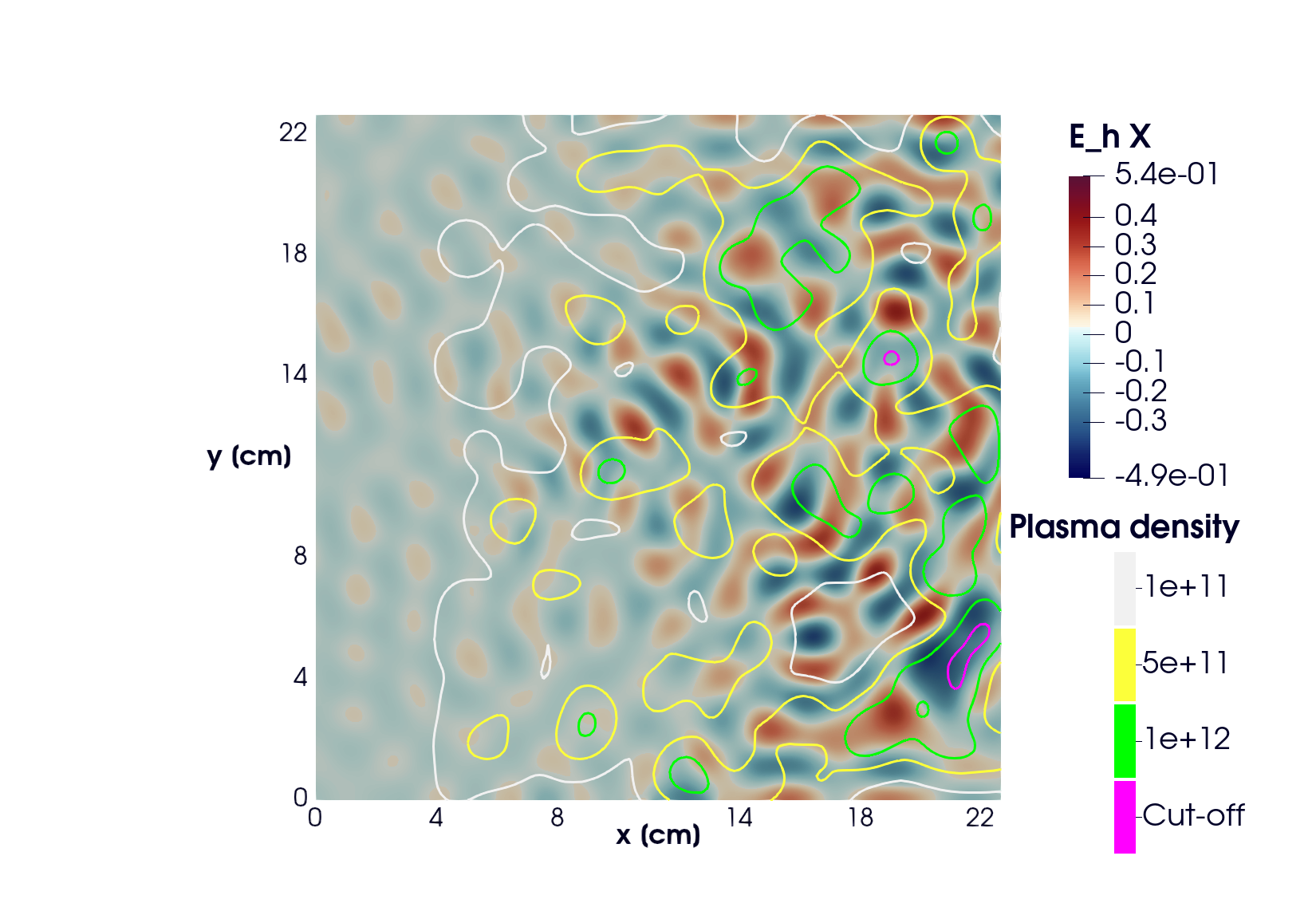}
          \includegraphics[width=0.194\textwidth, trim=14cm 4.8cm 13.9cm 3.5cm, clip]{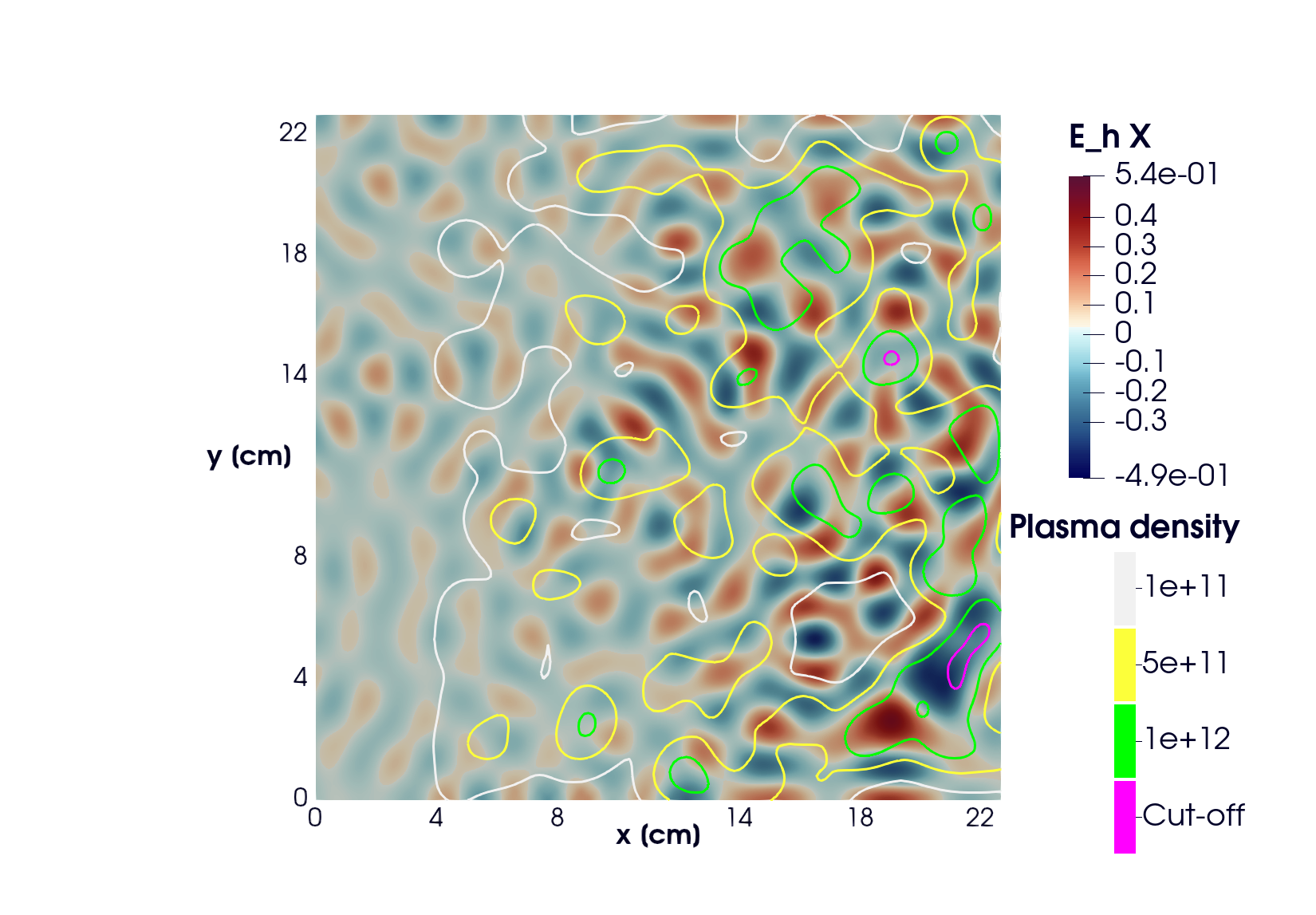}
         \caption{Top: real and imaginary part of $( \bE^{\mathrm{freq}} )_x$. Bottom: evolution of $( \bE_h )_y$ at $t = 5T,10T,15T,28T, 50T$ (increasing from left to right). The color map is fixed to the amplitude of $( \Re \{\bE^{\mathrm{freq}}\} )_x$, which is $[-0.49, 0.54]$.
         }
         \label{pics_example_high_fluctuations_Xmode_solution_X}
 \end{figure}

\begin{figure}[H]
          \centering %[trim=left bottom right top, clip]
          \includegraphics[width=0.4\textwidth, trim=0.8cm 1cm 0.8cm 4.6cm, clip]{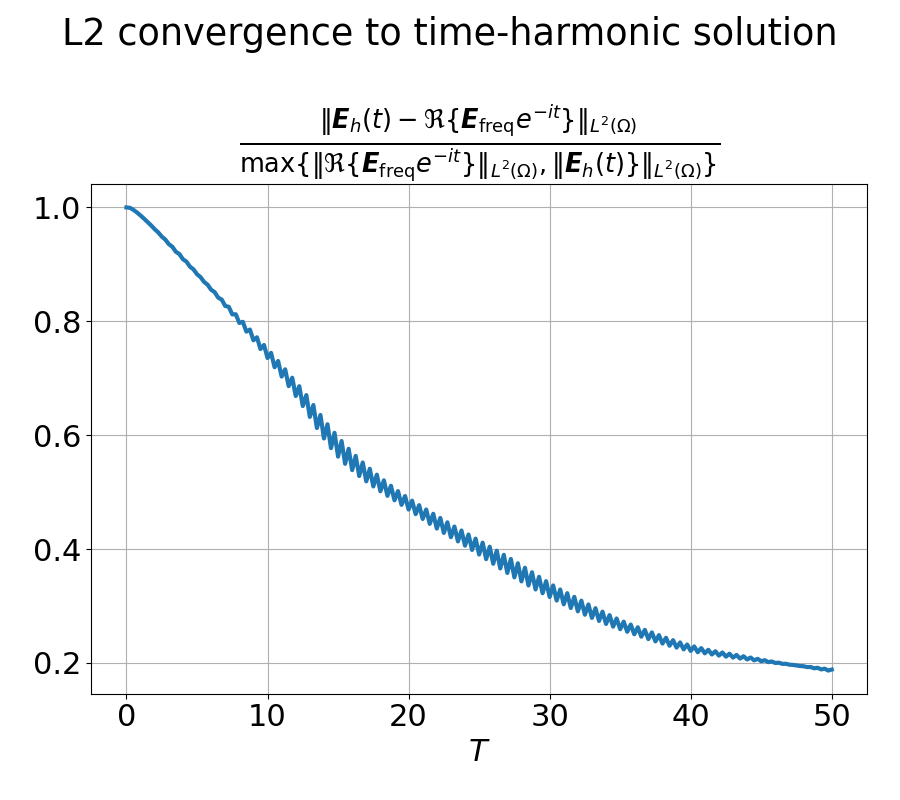}
          \includegraphics[width=0.59\textwidth, trim=9cm 2.5cm 2.3cm 4.5cm, clip]{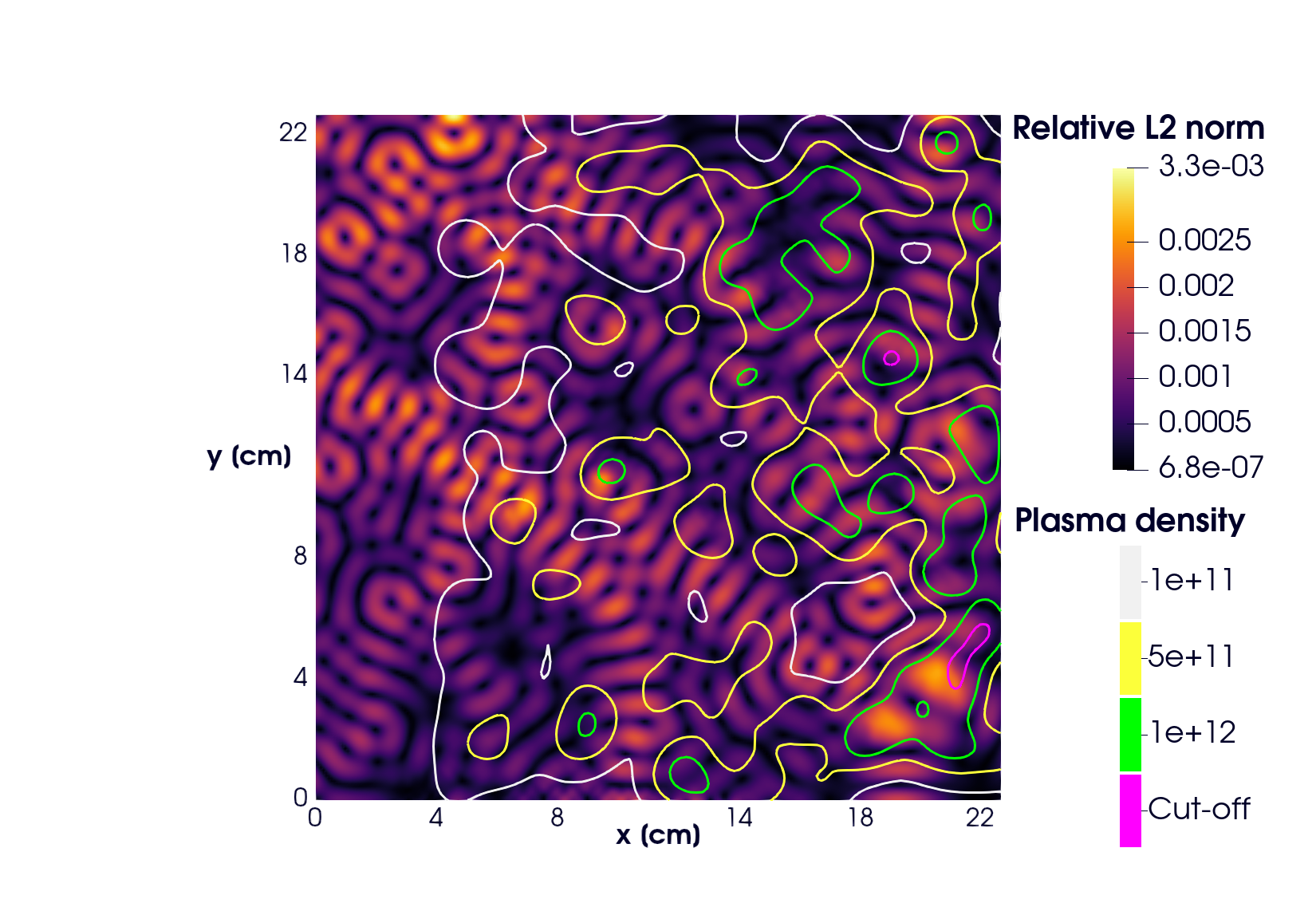}
         \caption{Time evolution of $\|\bs{R}\|_{L^2(\Omega)}$ (left) and pointwise Euclidian norm of $\bs{R}(\bx,50T)$ (right).
         }
         \label{pics_example_high_fluctuations_Xmode_conv_timeharmonic}
 \end{figure}

   In this experiment we consider an X-mode configuration, i.e. $\bs{e} = \bs{e}_y$ in \eqref{gaussian_beam_3d} and $(\bE_h)_z = \bs{0}$. Figure \ref{pics_example_high_fluctuations_Xmode_density} shows the electron-plasma density, which is the same as in figure \ref{pics_example_density}, but scaled to exclude the upper-hybrid resonance.
Figures \ref{pics_example_high_fluctuations_Xmode_solution_Y} ($y$-component) and \ref{pics_example_high_fluctuations_Xmode_solution_X} ($x$-component) show the real and imaginary parts of $\bE^{\mathrm{freq}}$ (top) and the evolution of $\bE_h$ (bottom).
  As in figures \ref{pics_example_Omode_conv_timeharmonic} and \ref{pics_example_1lobXmode_conv_timeharmonic}, the left-hand side panel in figure \ref{pics_example_high_fluctuations_Xmode_conv_timeharmonic} shows the time evolution of $\|\bs{R}\|_{L^2(\Omega)}$ and the right-hand side plot shows the pointwise Euclidian norm of $\bs{R}(\bx,50T)$.

 The left plot in figure \ref{pics_example_high_fluctuations_Xmode_conv_timeharmonic} shows that $\|\bs{R}\|_{L^2(\Omega)}$ decreases until a certain threshold.
 In this experiment, the oscillations are small and do not grow in time. As before, the right plot shows that the difference at the last time point is small and it is spread over the whole domain.

We have observed in all the experiments that $\|\bs{R}\|_{L^2(\Omega)}$ decreases fast during the first periods and it stagnates after a large time at a treshold of approximately $20\%$. Furthermore, in all the experiments the pointwise Euclidian norm of $\bs{R}$ at the last time point is small and it does not exhibit any clear pattern. The potential explanations given in remark \ref{remark_conv_timeharmonic} are coherent with the numerical observations.

To conclude, we can confirm that the solutions produced by the Poisson splitting method behave qualitatively as expected. Specifically, the solutions maintain the correct qualitative behaviour after large simulation times.

\section{Conclusions}
In the first part of this work, we have studied the Hamiltonian structure of the ideal problem and the effect of the boundary terms in the energy-balance. We derived two methods based on time-splitting geometric integrators for non-canonical Hamiltonian systems and proposed a treatment of the boundary terms, that guarantees the enforcement of the Silver-Müller boundary conditions.
In particular, one method uses Poisson splitting \eqref{scheme_Poisson} and the other method is based on Hamiltonian splitting \eqref{scheme_Hamiltonian}. For comparison in the analysis we included a method involving no splitting \eqref{scheme_Crank_Nicolson}, which amounts to the Crank-Nicolson method. For the Poisson splitting and Crank-Nicolson schemes we showed that under certain assumptions the methods are long-time stable.

In the second part of this work, the properties of the proposed schemes were tested in a series of numerical experiments. We checked the convergence order, the numerical stability and compared the performance in terms of cost and accuracy. Furthermore, we showed that the proposed methods yield accurate approximations of the energy, total charge and zero divergence of magnetic field. In the study we found that the Poisson splitting scheme has the best cost/accuracy ratio by one order of magnitude: for a given accuracy it is the cheapest,
and for a given cost it is the most accurate. 
Furthermore, it is unconditionally stable and our findings also hold with larger Courant numbers.

Finally, we used the Poisson splitting method to run three two-dimensional experiments of physical interest. 
We tracked the convergence of the solution to the time-harmonic regime,
by using the same spatial discretization of the frequency-domain problem. 
As expected, our experiments showed that in different scenarios the solutions converge towards the time-harmonic solution, with a residual error that is consistent with our 
long-time stability analysis.

\section*{Acknowledgements}
This work has been carried out within the framework of the EUROfusion Consortium, funded by the European Union via the Euratom Research and Training Programme (Grant Agreement No 101052200 — EUROfusion). Views and opinions expressed are however those of the author(s) only and do not necessarily reflect those of the European Union or the European Commission. Neither the European Union nor the European Commission can be held responsible for them.

Y.\ Güçlü is member of the \emph{Gruppo Nazionale per il Calcolo Scientifico} of the Italian \emph{Istituto Nazionale di Alta Matematica ``Francesco Severi''} (GNCS-INdAM).

\printbibliography

@article{dumont_heating_2012,
	title = {Heating and current drive by ion cyclotron waves in the activated phase of ITER},
	volume = {53},
	number = {1},
	journal = {Nuclear Fusion},
	author = {Dumont, R. J. and Zarzoso, D.},
	year = {2012},
	publisher = {IOP Publishing and International Atomic Energy Agency}
}

@article{jacquot_2d_2013,
	title = {{2D} and {3D} modeling of wave propagation in cold magnetized plasma near the {Tore} {Supra} {ICRH} antenna relying on the perfecly matched layer technique},
	volume = {55},
	number = {1},
	journal = {Plasma Physics and Controlled Fusion},
	author = {Jacquot, J. and Colas, L. and Clairet, F. and Goniche, M. and Heuraux, S. and Hillairet, J. and Lombard, G. and Milanesio, D.},
	year = {2013}
}

@article{tierens_unconditionally_2012,
	title = {An unconditionally stable time-domain discretization on cartesian meshes for the simulation of nonuniform magnetized cold plasma},
	volume = {231},
	number = {15},
	journal = {Journal of Computational Physics},
	author = {Tierens, W. and De Zutter, D.},
	year = {2012},
	pages = {5144--5156}
}

@book{de_boor_practical_2001,
	title = {A {Practical} {Guide} to {Splines}},
	publisher = {Springer New York, NY},
	author = {de Boor, C.},
	year = {2001}
}

@article{buffa_isogeometric_2011,
	title = {Isogeometric discrete differential forms in three dimensions},
	volume = {49},
	number = {2},
	journal = {SIAM Journal on Numerical Analysis},
	author = {Buffa, A. and Rivas, J. and Sangalli, G. and Vázquez, R.},
	year = {2011},
	pages = {818 -- 844}
}

@techreport{maj_hamiltonian_nodate,
	title = {Hamiltonian and metriplectic structures of the cold-plasma model and their discretization},
	url = {https://app.readcube.com/library/8de573ae-13bf-4887-8bf5-8d9936dd78d5/item/b920781d-8861-447f-a1fc-576234ac5ba5},
	author = {Maj, O. and Güçlü, Y and Lafitte, O. and Morrison, P.J. and Ratnani, A. and Sonnendrücker, E.},
}

@article{DASILVA201524,
title = {Stable explicit coupling of the Yee scheme with a linear current model in fluctuating magnetized plasmas},
journal = {Journal of Computational Physics},
volume = {295},
pages = {24-45},
year = {2015},
author = {da Silva, F. and Campos Pinto, M. and Després, B. and Heuraux S.},
}

@book{assous2018mathematical,
title={Mathematical Foundations of Computational Electromagnetism},
author={Assous, F. and Ciarlet, P. and Labrunie, S.},
series={Applied Mathematical Sciences},
year={2018},
publisher={Springer International Publishing},
pages={53-59}
}

@article{FEEC_YamanSaidMartin2022,
author = {Güçlü, Y. and Hadjout, S. and Campos Pinto, M.},
year = {2023},
title = {A Broken FEEC Framework for Electromagnetic Problems on Mapped Multipatch Domains},
volume = {97},
journal = {Journal of Scientific Computing}
}

@book{geomIntegration2006,
  title={Geometric Numerical Integration: Structure Preserving Algorithms for Ordinary Differential Equations},
  author={Hairer, E. and Wanner G. and Lubich C.},
  series={Computational Mathematics},
  year={2006},
  publisher={Springer Berlin, Heidelberg}
}

@article{bicgstab,
author = {van der Vorst, H. A.},
title = {Bi-CGSTAB: A Fast and Smoothly Converging Variant of Bi-CG for the Solution of Nonsymmetric Linear Systems},
journal = {SIAM Journal on Scientific and Statistical Computing},
volume = {13},
number = {2},
pages = {631-644},
year = {1992}
}

@article{iterative_methods,
author = {Barrett, R. and Berry, M. and Chan, T. and Donato, J. and Dongarra, J. and Eijkhout, V. and Pozo, R. and Romine, C. and Van der Vorst, H.},
year = {1996},
title = {Templates for the Solution of Linear Systems: Building Blocks for Iterative Methods},
volume = {64},
journal = {Mathematics of Computation}
}

@Inbook{why_difficult_Helmholtz_Ernst2012,
author={Ernst, O. G. and Gander, M. J.},
title={Why it is Difficult to Solve Helmholtz Problems with Classical Iterative Methods},
bookTitle={Numerical Analysis of Multiscale Problems},
year={2012},
publisher={Springer Berlin Heidelberg},
pages={325-363}
}

@book{chen2013introduction,
title={Introduction to Plasma Physics and Controlled Fusion: Volume 1: Plasma Physics},
author={Chen, F.F.},
year={2013},
publisher={Springer US}
}

@article{shiftedlaplacian_preconditioner,
author = {Gander, M. and Graham, I. and Spence, E.},
year = {2015},
title = {Applying GMRES to the Helmholtz equation with shifted Laplacian preconditioning: what is the largest shift for which wavenumber-independent convergence is guaranteed?},
volume = {131},
journal = {Numerische Mathematik}
}

@article{optimized_schwarz_maxwell_Dolean,
author = {Dolean, V. and Gander, M.J. and Gerardo-Giorda, L.},
title = {Optimized Schwarz Methods for Maxwell's Equations},
journal = {SIAM Journal on Scientific Computing},
volume = {31},
number = {3},
pages = {2193-2213},
year = {2009}
}

@book{bookDDM_Dolean,
author = {Dolean, V. and Jolivet, P. and Nataf, F.},
title = {An Introduction to Domain Decomposition Methods},
publisher = {Society for Industrial and Applied Mathematics},
year = {2015}
}

@article{reflectometry_FDTD_Silva_2019,
year = {2019},
volume = {14},
number = {08},
author = {F. da Silva and S. Heuraux and E. Ricardo and A. Silva and T. Ribeiro},
title = {Modelling reflectometry diagnostics: finite-difference time-domain simulation of reflectometry in fusion plasmas},
journal = {Journal of Instrumentation},
}

@article{2D_vs_3D_FDTD_daSilva_2022,
year = {2022},
publisher = {IOP Publishing},
volume = {17},
number = {01},
author = {F. da Silva and E. Ricardo and J. Ferreira and J. Santos and S. Heuraux and A. Silva and T. Ribeiro and G. De Masi and O. Tudisco and R. Cavazzana and O. D’Arcangelo},
title = {Benchmarking 2D against 3D FDTD codes for the assessment of the measurement performance of a low field side plasma position reflectometer applicable to IDTT},
journal = {Journal of Instrumentation}
}

@article{waveheating_Heuraux_2015,
author = {Heuraux, S. and da Silva, F. and Ribeiro, T. and Després, B. and Campos Pinto, M. and Jacquot, J. and Faudot, E. and Wengerowsky, S. and Colas, L. and Lu, L.},
year = {2015},
title = {Simulation as a tool to improve wave heating in fusion plasmas},
volume = {81},
journal = {Journal of Plasma Physics}
}

@article{McLachlan_Quispel_2002,
title={Splitting methods},
volume={11},
journal={Acta Numerica},
author={McLachlan, R. I. and Quispel, G. R. W.},
year={2002},
pages={341–434}
}

@article{ITER_position_reflectometers,
author = {da Silva, F. and Heuraux, S. and Manso, M.},
year = {2006},
pages = {1205 - 1208},
title = {Developments on reflectometry simulations for fusion plasmas: Application to ITER position reflectometers},
volume = {72},
journal = {Journal of Plasma Physics}
}

@book{abraham2012manifolds,
  title={Manifolds, Tensor Analysis, and Applications},
  author={Abraham, R. and Marsden, J.E. and Ratiu, T.},
  series={Applied Mathematical Sciences},
  year={2012},
  publisher={Springer New York}
}

@article{morrison1998,
author = {Morrison, P.},
year = {1998},
title = {Hamiltonian Description of the ideal fluid},
volume = {70},
journal = {Review of Modern Physics}
}

@article{Santos2021,
  title={A 3D CAD model input pipeline for REFMUL3 full-wave FDTD 3D simulator},
  author={Santos, J.M. and Ricardo, E. and da Silva, F. and Ribeiro, T. and Heuraux, S. and Silva, A.},
  journal={Journal of Instrumentation},
  year={2021},
  volume={16}
}

@article{DaSilva2024_overviewFDTD,
title = {An overview of the evolution of the modeling of reflectometry diagnostics in fusion plasmas using finite-difference time-domain codes},
journal = {Fusion Engineering and Design},
volume = {202},
year = {2024},
author = {da Silva, F. and Heuraux, S. and Ribeiro, T. and Ricardo, E. and Santos, J. and Silva, A. and Ferreira, J. and Vicente, J. and De Masi, G. and Tudisco, O. and Cavazzana, R. and Marchiori, G. and Luís, R. and Nietiadi, Y.}
}

@article{Heuraux2014_differentplasmas_FDTD,
title = {Study of wave propagation in various kinds of plasmas using adapted simulation methods, with illustrations on possible future applications},
journal = {Comptes Rendus Physique},
volume = {15},
number = {5},
pages = {421-429},
year = {2014},
author = {Heuraux, S. and Faudot, E. and da Silva, F. and Jacquot, J. and Colas, L. and Hacquin, S. and Teplova, N. and Syseova, K. and Gusakov, E.}
}

@book{FE_Maxwell_Monk,
    author = {Monk, P.},
    title = {Finite Element Methods for Maxwell's Equations},
    publisher = {Oxford University Press},
    year = {2003}
}

@book{Stix1992WavesIP,
  title={Waves in plasmas},
  author={Stix, T. H.},
  year={1992},
  publisher={American Institute of Physics Melville, NY}
}

@book{swanson2020plasma,
  title={Plasma Waves},
  author={Swanson, D.G.},
  series={Series in Plasma Physics and Fluid Dynamics},
  year={2020},
  publisher={CRC Press}
}

@book{Brambilla1998,
    author = {Brambilla, M.},
    title = {Kinetic Theory of Plasma Waves: Homogeneous Plasmas},
    publisher = {Oxford University Press},
    year = {1998}

}

@misc{petsc,
  author = {Balay, S. and Abhyankar, S. and Adams, M.F. and Benson, S. and Jed Brown and Peter Brune and Kris Buschelman and Emil~M. Constantinescu and Lisandro
  Dalcin and Alp Dener and Victor Eijkhout and Jacob Faibussowitsch and William~D.
  Gropp and V\'{a}clav Hapla and Tobin Isaac and Pierre Jolivet and Dmitry Karpeev
  and Dinesh Kaushik and Matthew~G. Knepley and Fande Kong and Scott Kruger and
  Dave~A. May and Lois Curfman McInnes and Richard Tran Mills and Lawrence Mitchell
  and Todd Munson and Jose~E. Roman and Karl Rupp and Patrick Sanan and Jason Sarich
  and Barry~F. Smith and Stefano Zampini and Hong Zhang and Hong Zhang and Junchao
  Zhang},
  title         = {{PETS}c {W}eb page},
  url           = {https://petsc.org/},
  year          = {2024}
}

@article{petsc4py,
  title         = {Parallel distributed computing using Python},
  author        = {Dalcin, L.D. and Paz, R.R. and Kler, P.A. and Cosimo, A.},
  journal       = {Advances in Water Resources},
  volume        = {34},
  number        = {9},
  pages         = {1124 - 1139},
  year          = {2011}
}

@misc{sebelin1997,
title = {Uniqueness and existence result around lax-milgram lemma: application to electromagnetic waves propagation in tokamak plasmas},
author = {Sebelin, E. and Peysson, Y and Litaudon, X. and Moreau, D. and Miellou, J.C. and Lafitte, O.},
year = {1997}
}

@phdthesis{LiseMarie_phdthesis,
	title = {Analyse mathématique et numérique de problèmes d'ondes apparaissant dans les plasmas magnétiques},
	school = {Université Pierre et Marie Curie},
	author = {Imbert-Gérard, L.M.},
	year = {2013},
}

@article{wellposedness_PEC_Back2015,
author = {Back, A. and Hattori, T. and Labrunie, S. and Roche, J. and Bertrand, P.},
year = {2015},
title = {Electromagnetic wave propagation and absorption in magnetised plasmas: Variational formulations and domain decomposition},
volume = {49},
journal = {ESAIM Mathematical Modelling and Numerical Analysis}
}

@mastersthesis{masterElena,
  author = {Moral Sanchez, E.},
  title = {A Time-Harmonic Numerical Model for Electromagnetic Wave Propagation in Plasmas},
  school = {Department of Mathematics, Technical University of Munich},
  year = {2022}
}

@article{MCDONALD1988337,
title = {Phase-space representations of wave equations with applications to the eikonal approximation for shor-wavelength waves},
journal = {Physics Reports},
volume = {158},
number = {6},
pages = {337-416},
year = {1988},
author = {McDonald, S.W.}
}

@article{Bernstein1975,
    author = {Bernstein, I.B.},
    title = {Geometric optics in space- and time-varying plasmas},
    journal = {The Physics of Fluids},
    volume = {18},
    number = {3},
    pages = {320-324},
    year = {1975}
}

@article{Kravtsov1969,
    author = {Kravtsov, Y.A.},
    title = {The geometric optics approximation in the general case of inhomogeneous and nonstationary media with frequency and spatial dispersion},
    journal = {Soviet Physics JETP},
    volume = {28},
    number = {4},
    year = {1969}
}
\end{document}